\documentclass[reqno]{amsart}
\usepackage{amscd,amssymb}

\newcommand{\splitpageyesno}[2]{#2}
\splitpageyesno{
\usepackage{splitpage}
}
{
\usepackage{fullpage}
}
\newcommand{\spmline}[2]{
\splitpageyesno{
\begin{multline*}
#1 \\
#2
\end{multline*}
}
{\[
#1 #2
\]
}
}
\newcommand{\spdisplay}[1]{
\splitpageyesno{\[ #1 \]}{$#1$}}
\usepackage[cmtip,arrow,curve,matrix]{xy}

\numberwithin{equation}{section}
\newtheorem{Theorem}[equation]{Theorem}

\newtheorem{Lemma}[equation]{Lemma}
\newtheorem{Proposition}[equation]{Proposition}
\newtheorem{Corollary}[equation]{Corollary}

\theoremstyle{definition}

\newtheorem{Definition}[equation]{Definition}
\newtheorem{Remark}[equation]{Remark}
\newtheorem{Example}[equation]{Example}

\newcounter{thmItem}

\newenvironment{thmList}{\begin{list}%
{\rm \roman{thmItem})}{\usecounter{thmItem}
\setlength{\labelwidth}{2em}
\setlength{\itemindent}{2em}
\setlength{\leftmargin}{0pt}
\setlength{\listparindent}{0pt}
\setlength{\parsep}{0pt}
\setlength{\partopsep}{0pt}
\setlength{\itemsep}{\medskipamount}
\setlength{\topsep}{\medskipamount}
}}{\end{list}}

\newcounter{textItem}

\newenvironment{textList}{\begin{list}%
{\rm (\arabic{textItem})}{\usecounter{textItem}
\setlength{\labelwidth}{2em}
\setlength{\itemindent}{2em}
\setlength{\leftmargin}{0pt}
\setlength{\listparindent}{0pt}
\setlength{\parsep}{0pt}
\setlength{\partopsep}{0pt}
\setlength{\itemsep}{\medskipamount}
\setlength{\topsep}{\medskipamount}
}}{\end{list}}

\newcounter{condItem}

\newcommand{\haff}{\Hat{\mathbb{A}}^{1}}
\newcommand{\artin}{\mathcal{A}}

\newcommand{\BU}[1]{BU\langle #1 \rangle}

\newcommand{\CatOf}[1]{(\text{#1})}
\newcommand{\cat}[1]{\mathcal{#1}}
\newcommand{\C}{\mathbb{C}}
\newcommand{\clockwise}{\mathrm{cl}}
\newcommand{\counterclockwise}{\mathrm{cc}}
\newcommand{\clchi}{\tilde{\chi}}
\DeclareMathOperator*{\colim}   {colim}
\newcommand{\cp}{{\C P^{\infty}}}
\newcommand{\cpplus}{\C P^{\infty}_{\plus}}
\newcommand{\CS}{\underline{C/S}}

\newcommand{\divisor}[1]{[#1]}
\newcommand{\divisorellA}{[\ell (A)]}

\DeclareMathOperator{\Defmtns}{Def}
\newcommand{\Deformations}[1]{\Defmtns (#1)}

\DeclareMathOperator{\Det}{Det}

\newcommand{\e}{0}
\newcommand{\einfty}{E_{\infty}}
\newcommand{\EllCurves}{\mathrm{Ell}}
\newcommand{\eps}{\epsilon}
\newcommand{\eqdef}{\overset{\text{def}}{=}}
\newcommand{\E}{\mathbf{E}}
\newcommand{\EGAI}[1]{\cite[#1]{EGAI}}

\newcommand{\funiv}{f_{\mathrm{univ}}}
\newcommand{\F}{\mathbf{F}}
\newcommand{\Fgps}{\mathrm{FGps}}
\newcommand{\fs}[1]{\underset{#1}{+}}
\newcommand{\fm}[1]{\underset{#1}{-}}
\newcommand{\fmlgpof}[1]{\widehat{#1}}
\newcommand{\frob}{\phi}

\newcommand{\G}{\mathbf{G}}
\newcommand{\genericunit}{\epsilon}

\newcommand{\Gmh}               {\widehat{\mathbb{G}}_m}
\newcommand{\GpOf}[1]{G_{#1}}
\newcommand{\GF}[1]{\mathbb{F}_{#1}}
\newcommand{\GS}{\underline{G/S}}

\newcommand{\hinfty}{H_{\infty}}

\newcommand{\hopgoe}{\mathrm{gh}}
\newcommand{\hot}{\widehat{\otimes}}
\newcommand{\HT}{\underline{H/T}}

\newcommand{\I}{\mathcal{I}}

\newcommand{\iso}{\cong}

\newcommand{\juniv}{j_{\mathrm{univ}}}

\DeclareMathOperator{\Ker}{ker}

\newcommand{\lvl}{\mathrm{level}}
\renewcommand{\L}{\mathcal{L}}
\newcommand{\Line}{\mathbb{L}}
\newcommand{\lie}{\mathcal{L}}
\newcommand{\lien}[1]{\lie({#1})}

\newcommand{\M}{\mathcal{M}}

\newcommand{\maxideal}{\mathfrak{m}}
\newcommand{\maxidof}[1]{\maxideal_{#1}}
\newcommand{\moeight}{MO\langle 8 \rangle}
\newcommand{\MU}[1]{MU\langle #1 \rangle}
\newcommand{\musix}{\MU{6}}

\newcommand{\noether}{\mathcal{N}}
\newcommand{\normff}{\mathbf{N}}
\newcommand{\normq}{N}
\DeclareMathOperator{\ob}{Ob}
\renewcommand{\O}{\mathcal{O}}

\newcommand{\pihat}{\widehat{\pi}_{0}}
\newcommand{\pointy}[1]{\langle #1 \rangle}
\newcommand{\plus}{+}

\newcommand{\Prob}{\mathcal{P}}
\newcommand{\psb}[1]{[ \! [#1] \! ]}
\newcommand{\psile}{{\psi_\ell^{\E}}}
\newcommand{\psilpe}{\psi_{\ell'}^{\E}}
\newcommand{\psilppe}{\psi_{\ell''}^{\E}}
\newcommand{\psilestar}{\left(\psi_\ell^{\E}\right)^\ast}
\newcommand{\psilpestar}{\left(\psi_{\ell'}^{\E}\right)^\ast}
\newcommand{\psilf}{{\psi_\ell^{\F}}}
\newcommand{\psilfstar}{\left(\psi_\ell^{\F}\right)^\ast}

\newcommand{\psilgstar}{\left(\psi_\ell^{\G}\right)^\ast}
\newcommand{\psilv}[1]{{\psi_\ell^{#1}}}

\newcommand{\psilge}[1]{\psi_\ell^{G/{#1}}}
\newcommand{\psilpge}[1]{\psi_{\ell'}^{G/{#1}}}
\newcommand{\psilppge}[1]{\psi_{\ell''}^{G/{#1}}}
\newcommand{\psilgestar}[1]{\left(\psi_\ell^{G/{#1}}\right)^*}

\newcommand{\ptdspaces}{\Spaces_{\plus}}

\newcommand{\ot}{\otimes}

\newcommand{\relpsilf}{{\psi_\ell^{\F/\E}}}

\newcommand{\restr}[1]{\vert_{#1}}

\newcommand{\quot}{q}

\DeclareMathOperator{\RingSpectra}{RingSpectra}
\newcommand{\rT}[1]{\tilde{T}_{#1}}

\newcommand{\slot}{\,-\,}
\DeclareMathOperator{\spec}{spec}

\DeclareMathOperator{\spf}{spf}
\newcommand{\specialpoint}[1]{#1_{0}}
\newcommand{\Spectra}[1]{\mathbb{S}_{#1}}
\DeclareMathOperator{\Spaces}{Spaces}
\newcommand{\struc}{\pi}
\newcommand{\suchthat}                {\;|\;}

\newcommand{\htensor}[1]{\underset{#1}{\hot}}
\newcommand{\TMF}{TMF}

\DeclareMathOperator{\tr}{Tr}
\DeclareMathOperator{\transfer}{\tau}

\newcommand{\un}[1]{\underline{#1}}
\newcommand{\uC}{\un{C}}
\newcommand{\uhom}{\un{\hom}}
\newcommand{\ulvl}{\un{\lvl}}

\newcommand{\Vreg}{V_{\mathrm{reg}}}
\newcommand{\Vstd}{V_{\mathrm{std}}}
\newcommand{\W}{\mathbb{W}}
\newcommand{\wreath}{\textstyle{\int}}

\newcommand{\xra}[1]{\xrightarrow{#1}}

\newcommand{\Z}{\mathbb{Z}}

\begin{document}

\title{The sigma orientation is an $\hinfty$ map} 

\author[Ando]{Matthew Ando}
\address{Department of Mathematics \\
The University of Illinois at Urbana-Champaign \\
Urbana IL 61801 \\
USA}
\email{mando@math.uiuc.edu}

\author[Hopkins]{Michael~J.~Hopkins}
\address{Department of Mathematics \\ 
         Massachusetts Institute of Technology\\ 
         Cambridge, MA 02139-4307\\ USA}
\email{mjh@math.mit.edu}

\author[Strickland]{Neil~P.~Strickland}
\address{Department of Pure Mathematics\\
         University of Sheffield\\
         Sheffield S3 7RH\\
         Great Britain}
\email{N.P.Strickland@sheffield.ac.uk}

\thanks{The authors were supported by NSF.  Ando's grant number was
DMS---0071482}

\date{Version 4.32, November 2002}

\begin{abstract}
In~\cite{AHS:ESWGTC} the authors constructed a natural map, called the
\emph{sigma orientation},  from the Thom spectrum $\musix$ to any
elliptic spectrum in the sense of 
\cite{ho:icm}.  $\musix$ is an 
$\hinfty$ ring spectrum, and in this paper we show that if $(\E,C,t)$
is the elliptic spectrum associated to the universal deformation of a
supersingular elliptic curve over a perfect field of characteristic
$p>0$, then the sigma orientation is a map of $\hinfty$ ring spectra.
\end{abstract}

\maketitle

\tableofcontents

\section{Introduction} \label{sec:introduction}

In \cite{ho:icm,AHS:ESWGTC}, we introduced the
notion of an \emph{elliptic spectrum} and showed that any elliptic
spectrum $(\E,C,t)$ admits a \emph{canonical} $\musix$ orientation
\[
     \musix \xra{\sigma (\E,C,t)} \E
\] 
called the
\emph{sigma orientation} (see also \S\ref{sec:sigma-orientation}).
We conjectured that  the spectrum $\TMF$ of ``topological modular
forms'' of Hopkins and Miller admits 
an $\moeight$ orientation, such that for any elliptic spectrum $(\E,C,t)$
the diagram 
\[
\begin{CD}
\moeight @>>> \TMF \\
@AAA @VVV \\
\musix @> \sigma (\E,C,t) >> \E
\end{CD}
\]
commutes.  For more about the conjecture, 
see~\cite{ho:icm} and the introduction to~\cite{AHS:ESWGTC}.

The conjecture seems now to be within reach, although
that is the subject of another paper in preparation.  The
proof depends on the following feature of the sigma orientation,
which was not proved in \cite{AHS:ESWGTC}.  
Let $C_{0}$ be a supersingular elliptic
curve over a perfect field $k$ of characteristic $p>0$, and let $\E$ be
the even periodic ring spectrum associated to the universal
deformation of the formal group of $C_{0}$ (so it is a form of
$E_{2}$).  The Serre-Tate theorem endows $\E$ with the structure of
an elliptic spectrum (see \S\ref{sec:serre-tate-theorem}), and so a
map of ring spectra 
\begin{equation} \label{oi-eq:17}
    \sigma: \musix \to \E.
\end{equation}
Goerss and Hopkins, building on work of Hopkins and Miller,  have
shown that $\E$ is an $\einfty$ ring spectrum \cite{gh:rcrs}; it is
classical that $\musix$ is.  We need to 
know that  the map \eqref{oi-eq:17} is an $\hinfty$ map;
we prove 
that in this paper. (The paper of Goerss and Hopkins has not
yet been published.  Our result depends only on the existence of the
$\hinfty$ structure, and so a cautious statement of the our result is
that if $\E$ is an $\hinfty$ ring spectrum, then the
map~\eqref{oi-eq:17} is $\hinfty$.  See Remark~\ref{rem-5}.)

In Part \ref{part:hinfty-orientations}, we study the general problem
of showing that an orientation 
\[
   \MU{2k} \xra{g} \E
\]
is $\hinfty$, i.e. that  
for each $n$ the
diagram 
\begin{equation} \label{oi-eq:21}
\begin{CD}
  D_{n}\MU{2k} @> D_{n} g >> D_{n}\E \\
  @VVV                      @VVV \\
  \MU{2k} @> g >> \E
\end{CD}  
\end{equation}
commutes up to homotopy.   
Our analysis is based
on~\cite{Ando:PowerOps}, which treats the case of $\MU{0}$, 
the Thom spectrum associated to
$\BU{0}= \Z\times BU$.  We review that case in \S\ref{sec:nec-mu0}, in
a form which generalizes to $\musix$.

Briefly, suppose that $\E$ is a homotopy
commutative ring spectrum with the property that
$\pi_{\text{odd}}\E = 0$ and $\pi_{2}\E$ contains a unit, so 
$\pi_{0}\E^{\cpplus}$ is the ring of
functions on a formal group $G = G_{\E}$ over $S = \pi_{0}\E$.  If 
\[
   \MU{0} \xra{g} \E
\]
is an orientation, then the composition 
\[
     (\cp)^{L} \xra{} \MU{0} \xra{g}  \E
\]
represents a trivialization $s_{g}$ of the ideal sheaf $\I_{G} (\e)$ of
functions on $G$ which vanish at the identity, that is, a coordinate
on the formal group $G$.  The association 
\[
    g \mapsto s_{g}
\]
gives a bijection between $\MU{0}$-orientations on $\E$ and
coordinates on $G$.  

Suppose in addition that $\pi_{0}\E$ is a complete 
local
ring with perfect residue field of characteristic $p>0$, and
the height of the formal group $G$ is finite.  In
\S\ref{sec:algebr-geom-hinfty}, following 
\cite{Ando:PowerOps},  we show that an $\hinfty$ structure on $\E$  
adds the following structure to the formal group $G$.
Given a map (of complete
local rings) $i: S\to R$, a finite abelian group $A$, and a
level structure  (in the sense of \cite{Drinfeld:EM}; see
\S\ref{sec:level-structures}) 
\begin{equation} \label{oi-eq:18}
\ell:A\to i^\ast G(R),
\end{equation}
there is a map
$\psi_\ell:S\to R$, and an isogeny $f_\ell:i^\ast
G\to\psi_\ell^\ast G$ with kernel $A$. (The behavior of this structure
with respect to variation 
in $A$ gives \emph{descent data for level structures} as described in
Definition \ref{def-1} or Proposition \ref{t-pr-descent-is-descent}.)

If $s$ is the coordinate on $G$ associated to an orientation
$g$, then the $\hinfty$ structure gives \emph{two} coordinates on 
$\psi_{\ell}^{*}G$: one
($\psi^{*}_{\ell} s$) comes 
from pulling back along $\psi_{\ell}$; the other ($\normq_{\ell}i^{*}s$)
is obtained from the invariant function 
\begin{equation} \label{eq:71}
        \prod_{a\in A}T_{a}^{*}i^{*}s
\end{equation}
on $i^{*} G$ by descent along the isogeny 
$f_{\ell}$ (see Proposition \ref{t-pr-OGp-NOG-new}).
In \S\ref{sec:comp-hinfty-struct-mu0}, we show that these two
coordinates arise from the two ways of navigating the
diagram 
\[
\begin{CD}
    (BA^{\ast}\times \cp)^{\Vreg\otimes L} 
    @>>>
    D_{n}\MU{0} @> D_{n}g>> D_{n} \E \\
  @. @VVV @VVV \\
 @.    \MU{0} @> g >> \E,
\end{CD}
\]
where $|A|=n$ and $\Vreg$ denotes the regular representation of $A^{\ast}$
(a key
point is that~\eqref{eq:71} is the Euler class of the bundle 
$\Vreg\otimes L$ associated to the orientation $g$).
It follows that if $g$ is an $\hinfty$ map, then 
\begin{equation} \label{eq:73}
      \psi_{\ell}^{*}s = \normq_{\ell}i^{*}s.
\end{equation}
This condition is equivalent to the condition in \cite{Ando:PowerOps};
see Remark \ref{rem-3}.

In \S\ref{sec:necess-cond-mu2k} we modify the 
discussion of \S\ref{sec:nec-mu0} to handle $\musix$-orientations. 
If
\[
   \musix \xra{g} \E 
\]
is an orientation, then the composition 
\[
    ((\cp)^{3})^{\prod_{i} (1-L_{i})}  \xra{} \musix \xra{g} \E
\]
represents a \emph{cubical structure} $s_{g}$ on the line bundle $\I_{G}
(\e)$ (Definition~\ref{defn-cubical-structure});   in
\cite{AHS:ESWGTC}, we showed that the assignment $g\mapsto 
s_{g}$ is a bijection between the set of $\musix$-orientations of $\E$
and the set $C^{3} (G;\I_{G} (\e))$ of cubical structures.

As before, a cubical structure $s$ on $\I_{G}
(\e)$ gives 
rise to \emph{two} cubical structures $\psi_{\ell}^{*}s$ and
$\Tilde{\normq}_{\ell} s$ on $\psi_{\ell}^{*}\I_{G} (\e)$.  If $s=s_{g}$ is
the cubical structure associated to an $\musix$-orientation $g$, 
then these two cubical structures correspond
to the two ways of navigating the diagram~\eqref{oi-eq:21}.  If $g$ is
an $\hinfty$ orientation, then the cubical structure $s$ must satisfy
the equation  
\begin{equation} \label{eq:72}
     \psi_{\ell}^{*} s = \Tilde{\normq}_{\ell}i^{*}s.
\end{equation}

In Proposition \ref{t-pr-norm-condition-suffices} we show that the
necessary conditions~\eqref{eq:73} and~\eqref{eq:72} are sufficient if
we suppose in addition that 
$p$ is not a zero divisor in $\E$.  We have given a direct proof, but
our argument amounts to showing that 
for $k\leq 3$, the character map of \cite{HKR:ggc} for
$\E^{0} ( D_{p} \BU{2k}_{\plus})$ is injective.

Thus we have reduced the problem of checking whether the orientation
\eqref{oi-eq:17} is $\hinfty$
to the problem of checking the equation
\eqref{eq:72}.  That problem is mostly a matter of 
recalling the construction of the sigma orientation; 
we do that in Part \ref{part:sigma-orientation}.  Here are the main points.

{\samepage
\begin{Definition}
An \emph{elliptic spectrum}  consists of 
\begin{enumerate}
  \item an even, periodic, homotopy commutative ring spectrum $\E$;
  \item an elliptic curve $C$ over $\spec \pi_{0}\E$; and
  \item an isomorphism of formal groups 
\[
t:G_{\E}\iso\fmlgpof{C}
\]
over $\spec \pi_{0}\E$.
\end{enumerate}
\end{Definition}
}
The Theorem of the Cube (or Abel's Theorem, for that matter) shows
that if $C$ is an elliptic curve, then $\I_{C} (\e)$ has a \emph{unique}
cubical structure $s (C/S)$.  If $C$ is the elliptic curve associated
to an elliptic spectrum $(\E,C,t)$, then 
\[
    (t^{3})^{*} \widehat{s} (C/S)
\]
is a cubical structure on $\I_{G_{\E}} (\e)$; the associated
$\musix$-orientation is the sigma orientation. (The name comes
from the fact that if $C$ is a complex elliptic curve, then there is a
simple formula for $s (C/S)$ in terms of the Weierstrass
$\sigma$-function, which shows that the sigma orientation for the
Tate curve is the Witten genus.  See \cite{AHS:ESWGTC})

Now suppose that $(\E,C,t)$ is an elliptic spectrum, that $\E$ is an
$\hinfty$ spectrum, and that $\pi_{0}\E$ is a complete 
local ring with perfect residue field of characteristic $p>0$.
Suppose that for each level 
structure
\[
   A \xra{\ell} i^{*} G (R)
\]
we are given 
an isogeny of elliptic curves 
\[
h_{\ell}: i^{*} C\xra{} \psi_{\ell}^{*}C
\]
with kernel $\divisorellA$, such that 
\[
       t^{*} \widehat{h}_{\ell} = f_{\ell}
\]
(Such structure, with compatibility with variation in $A$, is called
an \emph{$\hinfty$ elliptic spectrum} in Definition \ref{def-3}).
The 
uniqueness of the cubical structure $s (C/S)$ implies that 
\[
       \psi_{\ell}^{*} s (C/S) = \normq_{\ell} i^{*}s (C/S),
\]
which implies equation~\eqref{eq:72}.
Thus we have the following.

\begin{Proposition}[\ref{t-pr-sigma-hinfty-hinfty}]
If $(\E,C,t)$ is an $\hinfty$ elliptic spectrum, and 
$p$ is regular in $\pi_{0}\E$, then 
the sigma orientation 
\[
    \musix \xra{\sigma(\E,C,t)} \E
\]
is an $\hinfty$ map. \qed
\end{Proposition}

The Serre-Tate Theorem together with the result of Goerss and Hopkins
implies that the elliptic spectrum associated to the 
universal deformation of a supersingular elliptic curve over a perfect
field of characteristic $p$ is an $\hinfty$ elliptic spectrum
(Corollary \ref{t-co-hinfty-ec}),
and so the Proposition implies our result.

\begin{Corollary}[\ref{t-co-sigma-to-E-2-hinfty}]
If $(\E,C,t)$ is the elliptic
spectrum associated to the universal deformation of a supersingular
elliptic curve over a perfect field of characteristic $p>0$, then the
sigma orientation  
\[
      \musix \xra{\sigma (\E,C,t)} \E
\]
is $\hinfty$.
\end{Corollary}

We have analyzed $\hinfty$ ring spectra using the algebraic geometry of 
group schemes and in particular level structures, and we have analyzed
orientations (i.e. Thom isomorphisms) using the algebraic geometry of
line bundles.  Part \ref{part:even-peri-cohom} describes the
relationship to topology.  Section~\ref{sec:even-cohom-abel} discusses
the relationship between level structures and the cohomology of
abelian groups; this is a variation 
of \cite{HKR:ggc}.  Section~\ref{sec:cohom-thom-spectra} expresses some
familiar results about the even-periodic cohomology of Thom complexes
in the language of line bundles.

The construction of the homomorphism $\psi_{\ell}$ and the isogeny
$f_{\ell}$ and the proof of the sufficiency of the equations~\eqref{eq:73}
and~\eqref{eq:72} depend on two technical results (Propositions
\ref{t-pr-level-str-p-neq-0} and \ref{t-pr-char-map-inj-cyclic}) about
level structures.  We prove those results in Part
\ref{part:level-struct-isog}.  In order to make the discussion more
self-contained, we also recall there some results about level
structures, primarily from
\cite{Drinfeld:EM,KaMa:AMEC,Strickland:FiniteSubgps}.

\section{Notation}

\subsection{Groups}

If $X$ is an object in some category with products, and $J\subseteq I$
is an inclusion of sets, the projection map $X^I\to X^J$ will be
denoted $\pi_J$, while $\hat\pi_{J}$ will denote the projection map
$X^{I}\to X^{(I\backslash J)}$.    The set $J$ will often be 
indicated by the sequence of its elements.  For example, $\pi_{23}$
will denote projection to product of the $2^{\text{nd}}$ and
$3^{\text{rd}}$ factors, while $\hat\pi_{1}$ will denote projection
away from the first factor.  If $\sigma:I\to I$ is an
automorphism, the symbol $\pi_\sigma$ refers to the induced
automorphism of $X^I$. 

If $X$ is a commutative group object, then the symbol $\mu_J$ will
denote the map $X^J\to X$ obtained by composing $\pi_J$ with the
iterated multiplication, while $\hat{\mu}_{J}$ will denote the map
$\mu_{J}\times\hat\pi_{J}$.  
In punctual notation, 
\begin{align*}
      \mu_{23} (a_{1}, a_{2}, a_{3},\dots ) &= a_{2}  a_{3} \\
      \hat\mu_{23} (a_{1},a_{2},a_{3},\dotsc) & = (a_{1}, a_{2} a_{3}, \dotsc )
\end{align*}
and so forth.

If $X$ is a commutative group in a category of objects over a
base $S$, then the symbol $\e: S\to X$ will stand for the identity
section, and we shall generally abbreviate to $\struc$ the  symbol for
the structural map  
\[
     \pi_{\emptyset}: X\to S.
\]

\subsection{Formal schemes and formal groups}

As in \cite{AHS:ESWGTC},  we view affine schemes as
representable functors from rings to sets, and define a
\emph{formal scheme} to be a filtered colimit of affine schemes; the
value of the colimit is the colimit of the values
\begin{equation} \label{eq:55}
      (\colim_{\alpha} X_{\alpha}) (R) = \colim_{\alpha} X_{\alpha} (R).
\end{equation}

In this paper we make one important modification to the notation
\eqref{eq:55}. 
Recall from \EGAI{\textbf{0}, 7.1.2} that a \emph{preadmissible} ring
is a linearly topologized ring which contains an \emph{ideal of definition}: an open ideal $I$
such that, for all open neighborhoods $V$ of zero, $I^{n}\subset V$
for some $n>0$. An \emph{admissible} ring is preadmissible ring which is
complete and separated.  
If $R$ is an admissible ring, then the 
ideals of definition form a fundamental system of
neighborhoods of $0$.  The \emph{formal spectrum} of $R$ is the formal
scheme \EGAI{\textbf{I}, 10.1.2,10.6}
\[
    \spf R \eqdef \colim_{J} \spec R/J,
\]
where the colimit is over the poset of ideals of definition.  In fact
we shall only  need the case that $R$ is a local ring; a local ring is
admissible if it is complete and separated in its adic topology.  
If $R$ is an admissible ring and if $X$ is a formal scheme, then we
define $X (R)$ to 
be the set of \emph{natural transformations} 
\[
     \spf R \to X.
\]
Thus if $R$ is admissible then  $\haff (R)$ is the set of
\emph{topologically} nilpotent elements of $R$, rather than just the
set of nilpotent elements.  Similarly 
\[
      (\spf R') ( R )
\]
is the set of continuous ring homomorphisms from $R'$ to $R$.

A \emph{local scheme (of residue characteristic $p$)}  is a scheme of
the form $\spec R$, where $R$ is a local ring (of residue characteristic
$p$).  An \emph{local formal scheme (of residue characteristic $p$)}
is a formal scheme of the form 
$\spf R$, where $R$ is an admissible local ring (of residue
characteristic $p$).

Let $S$ be a formal scheme.  A \emph{formal group scheme} over $S$ is
a commutative group in the category of formal schemes over $S$.  
If $A$ is a finite abelian group, then $A_{S}$ will denote the
constant formal group scheme over $S$ given by $A$.
A \emph{formal group} over $S$ is a formal group scheme 
which is locally isomorphic to $S\times \haff$ as a
pointed formal scheme over $S$.  If $R$ is complete local ring, then a
formal group over $R$ means a formal group over $\spf R$.  If $G$ is a
formal group over $R$ and 
\[
    j: R\to R'
\]
is a map of complete local rings, then with the pull-back diagram 
\[
\begin{CD}
      j^{*} G @>>> G \\
       @VVV @VVV \\
      \spf R' @> j >> \spf R
\end{CD}
\]
in mind, we write $j^{*}G$ for the resulting formal group over $R'$.

We shall be primarily interested the case of a formal group $G$ of
finite height over a local formal scheme $S$ whose closed point
$\specialpoint{S}$ is the spectrum of a perfect field of
characteristic $p>0$.   Let 
$\specialpoint{G}$ be the fiber of $G$ over $\specialpoint{S}$,
i.e. the pull-back in the diagram
\[
\begin{CD}
   \specialpoint{G} @> i >> G \\
   @VVV     @VVV \\
   \specialpoint{S} @>>> S.
\end{CD}      
\]
By construction, $(G/S,i,\mathrm{id}_{S_{0}})$ is a
deformation of $\specialpoint{G}$ in the sense of Lubin and Tate
(see Definition~\ref{def-deformation}). 

This deformation is classified by a pull-back diagram 
\[
\begin{CD}
    G @>>> G' \\
    @VVV  @VVV \\
    S @>>> S',
\end{CD}
\]
where $(G'/S',\funiv,\juniv)$ is the universal deformation of
Lubin and Tate (\cite{LubinTate:FormalModuli}; see
\S\ref{sec:lubin-tate-groups}).
Various facts 
about $G/S$ then follow from facts about $G'/S'$ by change of base.

\subsection{Ideal sheaves associated to divisors}

If $G$ is a formal group or elliptic curve over a local formal scheme $S$,
then an \emph{(effective) divisor} on $G$ is a closed subscheme $D$ of
$G$ such
that the ideal sheaf $\I (D)$ is invertible, and $D$ is 
finite, free, and of finite presentation over $S$.
The sheaf $\I (D)$ is 
the inverse of the sheaf which is usually denoted $\O (D)$.
For example, if $w\in G (R)$ then 
$\I (w)$ is the ideal of functions on $G$ which \emph{vanish} at
$w$.  More generally, if
$W$ is a finite set and
\[
    \ell: W\to G (R)
\]
is a map of sets, then we will write $\I (\ell)$ for the ideal
associated to the divisor 
\[
 \divisor{\ell (W)} \eqdef \sum_{w\in W} \divisor{\ell (w)},
\]
so 
\[
    \I (\ell) \iso \bigotimes_{w\in W} \I (\ell (w)).
\]

\subsection{Spectra}

The category of spectra over a universe $U$ will be denoted
$\Spectra U$.  The category $\Spectra U$ is enriched over the category
$\ptdspaces$ of pointed topological spaces, and our notation will
reflect this.  Thus, the 
object $\Spectra U(E,F)$ will refer to a pointed topological space,
and for a pointed space $X$, $E^X$ is the function object (spectrum).
What would in category theory denoted $E\otimes X$ will in this case
be denoted $E\wedge X$.  There are natural homeomorphisms of pointed
spaces 
\begin{equation}\label{eq:44}
\ptdspaces(X,\Spectra U(E,F))\iso \Spectra U(E\wedge X,F)\iso
\Spectra U(E,F^X). 
\end{equation}

If $V$ is a vector bundle over a space $X$, then $X^V$ will refer to the
pointed space which is the Thom complex of $V$.  When $V$ is a virtual
bundle, then $X^V$ will refer to the Thom {\em spectrum} of $V$, arranged
so that the ``bottom cell'' is in the virtual (real)
dimension of $V$.  With this convention, the Thom spectrum of
an ``honest'' vector bundle is the suspension spectrum of the Thom
complex, so no real problem should come up when regarding an actual
vector bundle as a virtual.

We write $\Vstd$ for (the vector bundle over $B\Sigma_{n}$ associated
to) the standard complex representation of $\Sigma_{n}$, and if $A$ is an
abelian group, then $\Vreg$ will denote (the vector bundle 
over $BA$ associated to) the complex regular representation of $A$.

\subsection{Even periodic ring spectra}

A (homotopy commutative) ring spectrum $\E$ will be called \emph{even}
if $\pi_{\text{odd}}\E = 0$, and \emph{periodic} if $\pi_{2}\E$
contains a unit.   A ring spectrum $\E$ will be called
\emph{homogeneous} if it is a homotopy commutative algebra spectrum
over an even periodic ring spectrum. 
We will be particularly interested in homogeneous
spectra $\E$ in which the ring $\pi_{0}\E$ is preadmissible in some
natural topology (possibly discrete).  If $\E$
is a such a spectrum, then we write $\pihat \E$ for the separated
completion of $\pi_{0}\E$, and we define 
\[ 
  S_\E \eqdef \spf( \pihat \E )
\]
for the formal scheme defined by $\pihat \E$.

Let $\E$ be such a spectrum, and let $X$ be a space.  
If $\{X_\alpha\}$ is the set of compact subsets of 
$X$ and $\{I_{\beta}\}$ is the set of ideals of definition of $\pi_{0}\E$,
then $\pi_{0}\E^{X_{\plus}}$ is preadmissible in the 
topology defined by the kernels of the maps 
\[
 \pi_{0}\E^{X_{\plus}} \rightarrow 
\left(\pi_{0}\E^{(X_{\alpha})_{\plus}}\right)/I_{\beta},
\]
and we define $X_{\E}$ to be the formal scheme
\[
 X_{\E} = \spf \pihat \E^{X_{\plus}};
\]
this gives a covariant functor
from spaces to formal schemes over $S_{\E}$.  If $\F = \E^{X_{\plus}}$
then 
\[
   S_{\F} = X_{\E},
\]
and we shall use these notations interchangeably.

The most important example of these constructions is that $\E\iso
E_{n}$ is the spectrum associated to the universal deformation of a
formal group of height $n$ over a perfect field $k$ of characteristic
$p>0$, so
\[
    \pi_{0}\E\iso \W k\psb{u_{1},\dotsc ,u_{n-1}},
\]
and $X$ is a space with the property that $H_{*} (X,\Z)$ is 
concentrated in even degrees.  In that case, the natural map
of rings 
\[
     \pi_{0}\E^{X_{\plus}} \rightarrow \O (X_{\E}) = \pihat \E^{X_{\plus}}
\]
is an isomorphism, but in general all we have is a surjective map.

If $\E$ is a homogeneous ring spectrum, then it is complex
orientable, and 
\[
    G_{\E} = (\cp)_{\E} 
\]
is a formal group over $S_{\E}$.

\part{$\hinfty$ orientations}
\label{part:hinfty-orientations}
\section{Algebraic geometry of even $\hinfty$ ring spectra}

\label{sec:algebr-geom-hinfty}

\subsection{Descent data for level structures}

Let $\E$ be a homogeneous ring spectrum.
In this section we investigate the additional
structure which adheres to $\GpOf{\E} = (\cp)_{\E} / S_{\E}$ when
$\E$ is 
an $\hinfty$ spectrum (see \S\ref{sec:hinfty-ring-spectra}).  In order
to make precise 
statements, it is convenient to suppose that $\pi_{0}\E$ is 
a
complete
local ring with perfect residue field of characteristic $p>0$, and
that  $\GpOf{\E}$ is a 
formal group
of finite 
height.  In that case, we shall show that an $\hinfty$ structure on
$\E$ determines ``descent data for level structures'' on $\GpOf{\E}$.
In \S\ref{sec:desc-level-struct} we shall give a definition of this
notion in the usual language of descent; the definition we give there
is equivalent to the following.

\begin{Definition} \label{def-1} 
Let $G$ be a formal group 
over a formal scheme $S$.
\emph{Descent data for level
structures on} $G$ assign to every map of formal schemes $i:
T = \spf R \to S$, finite abelian group $A$, and level structure
(Definition \ref{def-level-structure})  
\begin{equation} \label{eq:18}
\ell:A_{T}\to i^\ast G,
\end{equation}
a map of formal schemes $\psi_{\ell}: T\to S$ and an isogeny
$f_\ell:i^\ast G\to\psi_\ell^\ast G$ with kernel 
\[
\divisorellA \eqdef \sum_{a\in A} \divisor{\ell (a)},
\]
satisfying the 
following. 
\begin{enumerate}
\item \label{descent-it-natural}
If $j: T\to T'$ is a map of formal schemes and
\[
   j^{*} \ell: A_{T'}\to j^{*}i^{*} G
\]
is the resulting level structure, then 
\[
     \psi_{j^{*}\ell} = j\circ \psi_{\ell},
\]
and 
\[
      f_{j^{*}\ell} = j^{*} f_{\ell}.
\]
\item \label{descent-it-2}
If $B\subseteq A$, then with
the notation
\begin{equation}\label{eq-hinfty-compat-diagram-intro}
\begin{CD}
B@>>> A @>>> A/B \\
@V \ell' VV @V \ell VV @VV \ell'' V \\
i^\ast G @= i^\ast G @>> f_{\ell'}> \psi_{\ell'}^{*} G,
\end{CD}
\end{equation}
we have 
\begin{equation}\label{eq-exact-equality-intro}
\begin{split}
\psi_{\ell''} &  = \psi_{\ell} : T \rightarrow S\\
f_\ell &  = f_{\ell''}\circ f_{\ell'}: i^{*}G \rightarrow
\psi_{\ell}^{*}G = \psi_{\ell''}^{*}G.
\end{split}
\end{equation}
\item \label{descent-it-3} 
If $\ell$ is the inclusion of the trivial subgroup, then $f_\ell$ and
$\psi_\ell$ are the identity maps.  Among other things this implies
that if $\ell$ and $\ell'$ differ by an automorphism of $A$, then
$f_\ell=f_{\ell'}$.
\end{enumerate}
We 
shall write $(\psi,f)$ for such descent data.
Formal groups with descent data for level structures form a category:
if $G/S$ and $G'/S'$ are two formal groups with descent data for level
structures, then a \emph{map} from $G'/S'$ to $G/S$ is a pull-back
\[
\begin{CD}
 G' @>>> G \\
 @VVV   @VVV \\
 S' @>>> S
\end{CD}
\]
in the category of formal schemes, such that the induced isomorphism
\[
     G'\rightarrow S'\times_{S} G
\]
is a group homomorphism, and such that 
the descent data for
$G$ pull back to the descent data for $G'$.
\end{Definition}

Let $\cat{C}$
be the category whose objects are homogeneous 
ring spectra 
$\E$ with the property that $\pi_{0}\E$ is 
a complete local 
ring with perfect residue field of 
positive characteristic and
$\GpOf{\E}$ is a 
formal group
of
finite height, and whose morphisms are maps $f: \E\to 
\F$ of ring spectra with 
the property that $\pi_{0}f$ is a map of local rings.
Let  $\hinfty\cat{C}$ be the subcategory of $\cat{C}$ consisting of
$\hinfty$ ring spectra and $\hinfty$ maps.
We shall
construct the dotted arrow in the diagram
\begin{equation} \label{eq:21}
\xymatrix{
{\hinfty\cat{C}} 
 \ar[d]
 \ar@{-->}[r]
&
{\CatOf{formal groups with descent data}}
 \ar[d] \\
{{\cat{C}}}
 \ar[r]
&
{\CatOf{formal groups}.}
}
\end{equation}
The main result is Theorem \ref{t-th-hinfty-adds}.

\subsection{Descent data from $\hinfty$ ring spectra}

The basic operation on the homotopy groups of an $H_\infty$ ring
spectrum is the transformation
\[
D_n:\pi_0\E\to\pi_0\Spectra U(D_n S^0,\E)=\pi_0 \E^{{B\Sigma_n}_\plus}.
\]
This map is multiplicative in the sense that $D_n(fg)=D_n(f)D_n(g)$,
but it is not quite additive.  In fact, it follows from
Proposition~\ref{t-pr-may-et-al} that
\begin{equation} \label{eq:17}
D_n(f+g)=\sum_{i+j=n}\tr_{ij}D_i(f)D_j(g),
\end{equation}
where $\tr_{ij}$ is the transfer map associated to the inclusion
$\Sigma_i\times\Sigma_j\subseteq\Sigma_n$.

If $\E$ is a spectrum such that $\pi_{0}\E$ is a 
complete
local ring with perfect residue field of 
characteristic $p>0$, and the formal group $\GpOf{\E}$ is
of finite height, then 
there is a slightly more convenient 
operation to work with. 
Suppose that $A$ is a finite abelian group, and let $A^{*}$ be its
group of complex characters.  With these hypotheses, 
Proposition \ref{t-E-BA-ast} says that the natural
map~\eqref{eq-hom-a-g-map}  
\begin{equation} \label{eq:56}
      (BA^{\ast})_{\E} \rightarrow \uhom (A,\GpOf{\E})
\end{equation}
is an isomorphism.  Define a functor $D_A:\Spectra 
U\to\Spectra U$ by 
\begin{equation}\label{eq-a-extended-power}
D_A(X) = \lie(U^{A^\ast},U) \underset{A^*}{\wedge}X^{(A^{\ast})},
\end{equation}
where
$X^{(A^{\ast})}$ denotes the external smash product
\[
\bigwedge_{\alpha\in A^\ast}X\in\ob\Spectra {U^{A^\ast}}
\]
(We will also have use for the functor on pointed spaces given by the
the analogue of~\eqref{eq-a-extended-power}).

\begin{Definition}\label{def-psile} 
Given a complete local ring $R$ and a level
structure~\eqref{def-level-structure}  
\[
   A_{\spf R} \xrightarrow{\ell} i^{*}G,
\]
we define $\psile: \pi_{0}\E \rightarrow
R$ to be the map given by the composition
\spmline{
\pi_0\E
 \xra{D_{A}}
\pi_0\Spectra U(D_A S^0,\E) =
\pi_0\E^{BA^\ast_\plus} \rightarrow
}
{
 \O ((BA^{\ast})_{\E}) 
\xra{\chi_{\ell}}
R,
}
where $\chi_{\ell}$ is the map classifying the
homomorphism $\ell$ as in~\eqref{eq:25}.
\end{Definition}

\begin{Lemma}[\cite{Ando:PowerOps}]\label{lem-psi-additive}
The map $\psile$ is a continuous ring homomorphism.
\end{Lemma}

\begin{proof}
$\psile$ is certainly multiplicative.  It's additive because 
equation \eqref{eq:17} and the double coset
formula imply that $\psile(x+y)-\psile(x)-\psile(y)$ is a sum
of elements in the image of the transfer map from proper subgroups of
$A^\ast$. The result therefore follows from
Proposition~\ref{t-pr-ra-kills-transfers}.

To see that $\psile$ is continuous, note that 
\begin{equation}\label{eq:67}
\O (( BA^\ast)_{\E}) \iso \O (\uhom (A,\GpOf{\E})) 
\end{equation}
is a local ring by Proposition
\ref{t-pr-hom-finite-flat}. 
It suffices to show that for $y$ in the maximal ideal of $\pi_{0}\E$,
$D_{A}y$ is in the 
maximal ideal of $\O (( BA^\ast)_{\E})$.  Since, modulo the
augmentation ideal of 
$\pi_{0}\E^{BA^{\ast}_{\plus}}$ we have 
\[
     D_{A} y = y^{|A|},
\]
it follows that $\psile$ is continuous.
\end{proof}

\begin{Remark}  
It is not necessary to use an abelian group  $A$.
The initial ring $R$ over which we may define a \emph{ring homomorphism}
\[
    \psi: \pi_{0}\E \xra{} \pi_{0} \E^{B (\Sigma_{k})_{\plus}} \xra{} R
\]
as above is the quotient of $\pi_{0}\E^{B (\Sigma_{k})_{\plus}}$ by
the ideal generated by the images of transfers from proper subgroups
of $\Sigma_{k}$.  Strickland \cite{Strickland:MorESym} shows that
\[
   \spf \left( \pi_{0}E_{n}^{B (\Sigma_{k})_{\plus}} / \text{proper transfers}\right) 
\]
is the scheme of ``subgroups of order $k$ of $\GpOf{{E_{n}}}$''.  
The analysis of this paper can be carried through in that setting.
\end{Remark}

The operation $\psi_\ell$ is clearly natural in the sense that given a
map
$f:\E\to\F$ of $H_\infty$ spectra of the indicated kind, with the
property that $\pi_{0}f$ is continuous, then 
the level structure $\ell$ gives a level structure 
\[
       A_{\spf R'} \xra{\ell} j^{*}\GpOf{\F},
\]
where $R'$
and $j$ are defined by  
\[
     \pi_{0}F \xra{j} R' = R \underset{i,\pihat \E,\pihat f}{\hot} \pi_{0}\F,
\]
and the
diagram
\[
\begin{CD}
R'  @< \psilf << \pi_0 \F  \\
@AAA  @AA \pi_0 f A  \\
R  @< \psile <<  \pi_0\E 
\end{CD}
\]
commutes.

In the language of algebraic geometry, let 
\[
   T = \spf R \xra{i} S_{\E}.
\]
The map $\psile$ is a map of formal schemes
\[
\psile: T \to S_{\E},
\]
and  the naturality is expressed in terms of the commutative diagram
\[
\begin{CD}
T\underset{i,S_{\E},S_{f}}{\times} S_{\F} @> \psilf >> S_{\F} \\
@VVV @VV S_{f} V \\
T @> \psile >> S_{\E}.
\end{CD}
\]
Making use of the isomorphism
\[
   T\times_{i,S_{\E},S_{f}} S_{\F} \cong i^{*} S_{\F}, 
\]
we find that the map $\psilf$ can be factored
through a \emph{relative} map
\begin{equation}\label{eq-psi-f}
\relpsilf: i^{*} S_{\F}\to\psile^\ast S_{\F}
\end{equation}
as in
\begin{equation} \label{eq:61}
\xymatrix{
{i^\ast S_{\F}}
 \ar@/^1pc/[drrr]^{\psilf}
 \ar[dr]^{\relpsilf}
 \ar@/_1pc/[ddr] \\
&
{\psilestar S_{\F}}
 \ar[rr]_-{S_{f}^{*}\psile}
 \ar[d]
& &
{S_{\F}}
 \ar[d]^{S_{f}}
\\
& 
{T}
 \ar[rr]_{\psile}
& & 
{S_{\E}.}
}
\end{equation}

For example, let $G=\GpOf{\E}$, and take  $\F=\E^{\cpplus}$ so that $G=S_{\F}$.
There results a map 
\begin{equation} \label{eq:5}
\psilge{\E}:i^\ast G\to\psilestar G.
\end{equation}
This map turns out to be a homomorphism of groups, as one can see by
considering the ($H_\infty$) map 
\[
\E^{\cpplus}\to\E^{(\cp\times\cp)_{\plus}}
\]
coming from $\mu:\cp^2\to\cp$.  We shall eventually show (Proposition \ref{t-pr-isogeny-degree})
that $\psilge{\E}$ is an isogeny, with kernel $\ell: A\to i^{\ast}G.$
In order to give the proof, it is essential to understand the effect
of the operation $\psile$ on the cohomology of Thom complexes.

Suppose that
$V$ is a virtual bundle over a space $X$,
and write
\[
\F=\E^{X_{\plus}}.
\]
If we start with an element 
$m\in \pi_0\Spectra U\left((X)^{V},\E\right)$
and follow the construction of the map $\psile$, we wind up with an
element 
$\psilv V(m)$ in 
\begin{equation} \label{eq:3}
R\underset{\chi_{\ell},\pihat \E^{BA^\ast_\plus}}{\hot} \pihat \Spectra
U\left((BA^\ast\times X)^{\Vreg\otimes V},\E\right), 
\end{equation}
where $\Vreg$ denotes the regular representation of $A^{\ast}$.
As before, this map is additive, and in fact $\psilf$-linear:
\begin{equation}\label{eq-psi-linearity}
\psilv V(xm)=\psilf(x)\psilv V(m).
\end{equation}

Let $T=\spf R$; then we have a commutative diagram 
\[
\xymatrix{
{i^{*}S_{\F}}
 \ar[r]^-{\chi_{\ell}}
 \ar[d] 
&
{\uhom (A,\GpOf{\F})}
 \ar[r]
 \ar[d]
&
{S_{\F}}
\ar[d]
\\
{T}
 \ar[r]^-{\chi_{\ell}}
 \ar@/_2pc/[rr]^{i}
&
{\uhom (A,\GpOf{\E})}
 \ar[r]
&
{S_{\E}}
}       
\]
in which all the squares are pull-backs.
In the language of \S\ref{sec:cohom-thom-spectra}, the element $m$ is
a section of 
the line bundle $\Line (V)$ over $S_{\F}$.  Elements of \eqref{eq:3}
are sections of 
\[
    \chi_{\ell}^{*}\Line (\Vreg\otimes V)
\]
over $i^{*}S_{\F}$.  
Taking into account the
linearity~\eqref{eq-psi-linearity} we find that the map $\psilv V$ can
be interpreted as a map
\begin{equation}\label{eq-psilf-line-bundle}
\psilv V:\psilfstar\Line(V)
\to 
\chi_{\ell}^{*}\Line(\Vreg\otimes V)
\end{equation}
of line bundles over $i^{*}S_{\F}$.

\begin{Lemma}\label{t-le-properties-psilv}
The map $\psilv V$ has the following properties
\begin{textList}
\item If $m$ is a trivialization of $\Line (V)$, then $\psilv{V} (m)$ is a
trivialization of $\chi_{\ell}^{*} \Line (\Vreg\otimes V)$.
\item With the obvious identifications
\[
\psilv{V_1\oplus V_2}=\psilv{V_1}\otimes\psilv{V_2}.
\]
\item If $f:Y\to X$ is a map, then
\[
\psilv{f^\ast V} = f^\ast\psilv{V}.
\]
\end{textList}
\end{Lemma}

The most important example is the case $X=\cp$ and $V=L$, so $S_{\F} =
G = \GpOf{\E}$ and $\Line (L) = \I_{G} (\e)$ \eqref{eq:16}.  Then
\eqref{eq:37} gives an isomorphism 
\[
 \chi_{\ell}^{*}\Line (\Vreg\otimes L) \cong \I_{i^{*}G} (\ell),
\]
and so we may think of $\psilv{L}$ as a map of line bundles over $i^{*}G$
\[
       \psilfstar \I_{G} (\e) \rightarrow \I_{i^{*}G} (\ell),
\]
or on sections a $\psilf$-linear map 
\begin{equation} \label{eq:4}
      \Gamma (\I_{G} (\e)) \rightarrow \Gamma (\I_{i^{*}G} (\ell)).
\end{equation}

\begin{Proposition}[\cite{Ando:PowerOps}]\label{t-pr-isogeny-degree}
\begin{textList}
\item The map
\[
\psilge{\E}: i^{*} G \to \psilestar G
\]
of~\eqref{eq:5} is an isogeny with kernel $\divisorellA$.
\item If the isogeny $\psilge{\E}$ is used to identify 
\[
      \psilgestar{\E} I_{\psilestar G} (\e) \iso
      I_{i^{*}G} (\ell)
\]
as in~\eqref{eq:84},
then the map $\psilv L$ \eqref{eq:4} sends a coordinate $x$ on $G$ to
the trivialization $\psilgestar{\E} \psilestar x$ of $I_{i^{*}G} (\ell)$.
\end{textList}
\end{Proposition}

\begin{proof}
First, observe that $\psilge{\E}$ is an isogeny of degree $|A|$.  This
follows from the Weierstrass Preparation Theorem
\cite[pp. 129--131]{Lang:cf}, because, after 
reducing modulo the maximal ideal in 
$R$, the ring homomorphism 
\[
    (\psilge{\E})^{*}: R\htensor{\psile} \pihat E^{\cp_{+}}
                      \rightarrow 
                     R\htensor{i} \pihat E^{\cp_{+}}
\]
sends a coordinate $x$ to $x^{|A|}$.  To see this, note that 
the composition 
\[
\pi_{0}\F
\xrightarrow{D_{A}}
\pi_{0}\F^{BA^{*}_{\plus}}
\xra{\pi_{0}\F^{S^{0}\rightarrow BA^{*}_{\plus}}}
\pi_{0}\F
\]
is the map $x\mapsto x^{|A|}$.

To prove the first part of the Proposition, it remains to show
that the $\divisorellA$ is contained 
in the kernel, i.e. that $\psilge{\E} \ell = 0$, i.e. that if 
$x$ is a coordinate on $G_{\E}$, then 
$(\psilge{\E})^{*} \psilestar x$ vanishes on $\divisorellA$.  A coordinate
on $G$ is a generator the ideal $\I (\e)$,
which the zero
section identifies with $\Line (L)$.
The commutativity of the diagram
\begin{equation} \label{eq:13}
\xymatrix{
{\pi_0\E^{\cpplus}}
 \ar[d]_{D_{A}}
 \ar@/_4pc/[dd]_-{\psilge{\E}}
&
{\pi_{0}\E^{\cp^{L}}}
  \ar[l]_-{\pi_{0}\E^{\zeta}}
  \ar[d]^-{D_{A}}
 \ar@/^4pc/[dd]^-{\psilv{L}}
\\
{\pi_0\E^{BA^\ast\times\cpplus}}
 \ar[d]
&
{\pi_0\E^{(BA^\ast\times\cp)^{\Vreg\otimes L}}}
 \ar[d]
 \ar[l]^-{\pi_{0}\E^{\zeta}}
\\
{R{\hot}
\pihat \E^{BA^\ast\times\cpplus}} 
&
{R{\hot}
\pihat \E^{(BA^\ast\times\cp)^{\Vreg\otimes L}}}
 \ar[l]
}
\end{equation}
(in which the tensor products are taken over the ring
$\pihat \E^{BA^\ast_\plus}$) shows that
$\psilge{\E}$ takes a section of $\I (\e)$
to a 
section of $\I_{i^{*}G} (\ell)$.  The claim about $\psilv{L} x$
also follows from inspection of the diagram~\eqref{eq:13}.
\end{proof}

In fact Proposition \ref{t-pr-isogeny-degree} gives a simple 
description of the map $\psilv{V}$ for a general
virtual bundle $V$, and in 
particular, shows that it is determined by maps which have already
been constructed.  We shall express the answer in the language of
\S\ref{sec:cohom-thom-spectra}, where line bundles of the form 
$\Line(V)$ are computed in terms of divisors. As in
\S\ref{sec:cohom-thom-spectra}, it is illuminating to 
work at the outset with $V\otimes L$ over $X\times\cp$ and then pull
back along the identity section of $G_{\F}$.

With this in mind, let $\F=\E^{X_{\plus}}$, let 
\[
\G=\F^{\cpplus}= \E^{(\cp \times X)_{\plus}},
\]
and let $G = G_{\F}= S_{\G}$.
Let $D=D_V$ be the divisor on $G$
corresponding to $V$ as in Proposition~\ref{t-pr-line-v},
so that there is an isomorphism
\[
t_{V\otimes L}:\Line(V\otimes L)\iso \I(D^{-1}).
\]
Thus we may replace the domain of 
\[
\psilv{V\otimes L}:\psilgstar \Line(V\otimes L)\to
\chi_{\ell}^{*}\Line(\Vreg\otimes V\otimes L).
\]
with 
\[
\psilgstar\I(D^{-1}) \cong \psilgestar{\F}\psilfstar I (D^{-1}).
\]
Using the analogous isomorphism~\eqref{eq:36}
\[
\chi_{\ell}^{\ast}\Line(\Vreg\otimes V\otimes L) 
     \iso
\I \left(\sum_a T_a^\ast D^{-1}\right)
\]
to interpret the range, 
we may think of $\psilv{V\otimes L}$ as a map 
\begin{equation}\label{eq-psilvl-map}
\psilv{V\otimes L}:\psilgestar{\F}\psilfstar \I(D^{-1})
\to
\I\left(\sum_a T_a^\ast D^{-1}\right).
\end{equation}

\begin{Proposition}\label{prop-that-damn-isogeny-lemma}
In the guise of~\eqref{eq-psilvl-map}, the map $\psilv{V\otimes L}$
is given by
\[
f\mapsto \psilgestar{\F}\psilfstar f.
\]
\end{Proposition}

\begin{proof}
There are actually two assertions.  One is that 
\[
\langle f\rangle\ge D^{-1}\implies
\langle (\psilfstar f)\circ\psilge{\F}\rangle \ge \sum T_a^\ast D^{-1}.
\]
The other is that this gives the map $\psilv{V\otimes L}$.  The
verification of both assertions follows the lines of the proof of
Proposition \ref{t-pr-isogeny-degree}.  Indeed, everything involved takes
Whitney sums in $V$ to tensor products, and commutes with base change
in $V$.  It suffices then to verify the case when $X$ is a single
point, and $V$ has dimension 1.  In this case, the isomorphism $t_L$
is given by the inclusion of the zero section $\cp\to\cp^L$, and the
result follows from naturality of the maps $\psi_\ell$, as in the
diagram \eqref{eq:13}.
\end{proof}

The results of this section assemble to give the following.

\begin{Theorem}\label{t-th-hinfty-adds} 
Let $\E$ be a homogeneous $\hinfty$ ring spectrum.
Suppose that $\pi_{0}\E$ is 
a 
local ring
with perfect residue field of characteristic $p>0$, and the formal group
$G=\GpOf{\E}$ is 
of finite height. The rule which
associates to  a 
level structure  
\begin{equation} 
\ell:A_{\spf R}\to i^\ast G
\end{equation}
the map of formal schemes
$\psile: \spf R\rightarrow S_{\E}$
and the isogeny 
\[
\psilge{\E} :i^\ast G\to\psi_\ell^\ast G
\]
is descent data for level structures on the formal group
$G/S_{\E}$, and gives the dotted arrow in the diagram \eqref{eq:21}.
\end{Theorem}

\begin{proof}
Lemma \ref{lem-psi-additive} and Proposition \ref{t-pr-isogeny-degree}
show that 
$\psile$ is a ring homomorphism and $\psilge{\E}$ is an isogeny with
kernel $\divisorellA.$  
The compatibility of $\psile$ and $\psilge{\E}$
with variation in $A$ as described in 
Definition \ref{def-1} follows from the commutativity of the
diagrams~\eqref{eq:42}; a proof is given in
Appendix~\ref{sec:comp-oper-deriv}. 
\end{proof}

\section{A necessary condition for an $\MU{0}$-orientation to be $\hinfty$}

\label{sec:nec-mu0}

Let $\MU{0}$ be the Thom spectrum of the tautological bundle over
$\Z\times BU$, and let $\E$ be a homogeneous
ring spectrum.  In \S\ref{sec:spectrum-mu0} we 
recall that  to give a map of
(homotopy commutative) ring spectra
\begin{equation} \label{eq:1}
   g: \MU{0}\to \E.
\end{equation}
is to give a coordinate $s$ on
$G = \GpOf{\E}$.  

In \S\ref{sec:comp-hinfty-struct-mu0} we give
a necessary condition  for the map $g$ to 
be a map of $\hinfty$ spectra,
in the case that $\E$
is an $\hinfty$ ring spectrum, that $\pi_{0}\E$ is a 
complete local
ring with perfect residue field of characteristic $p>0$, 
and that the formal group $G$ is 
of finite height.
The result
may described as follows.   

Let $s=s_{g}$ be the coordinate on
$G$ associated to the orientation~\eqref{eq:1}.  In
\S\ref{sec:algebr-geom-hinfty} we showed that the hypotheses on $\E$
give descent data for level structures on 
$G$.   Given 
a level structure~\eqref{def-level-structure} 
\begin{equation} \label{eq:2}
      A_{T} \xra{\ell} i^{*}G,
\end{equation}
we 
get \emph{two} coordinates on the formal group $\psilestar G$: one
is just $\psilestar s$, the other is the norm $\normq_{\ell} i^{*}s$
of the coordinate $i^{*}s$ with respect to the isogeny 
\[
   i^{*}G \xra{\psilge{\E}} \psi_{\ell}^{*}G
\]
as in Proposition \ref{t-pr-OGp-NOG-new}.  We show that these two
coordinates 
correspond to the two ways of going around the diagram 
\[
\begin{CD}
  D_{A}\MU{0} @>>> D_{A}\E \\
 @VVV @VVV \\
  \MU{0} @>>> \E;
\end{CD} 
\]
the main result is Proposition \ref{t-pr-mu-zero-hinfty}.

\subsection{$\hinfty$ structures on Thom spectra of infinite loop spaces}

Suppose that $B\to \Z\times BO$ is a homotopy multiplicative map, and let $M$
be the associated Thom spectrum.  The spectrum $M$ has a natural
multiplication.  If $W:X\to B$ is a vector bundle over $X$ with a
$B$-structure,
then the Thom complex $X^W$ comes equipped with a canonical $M$-Thom
class
\[
\Phi_M(W):X^W\to M
\]

\begin{Lemma}\label{lem-thom-class}
The Thom class $\Phi_M(W)$ has the following properties.
\begin{thmList}
\item  It is multiplicative:
$\Phi_M(W\oplus W')=\Phi_M(W)\Phi_M(W').$
\item  It is preserved under base change: given
$f:X\to Y$,
\[
\Phi_M(f^\ast W) = f^\ast\Phi_M(W).
\]
\end{thmList} \qed
\end{Lemma}

An infinite loop map 
\[
   B\to \Z\times BO
\]
gives, for every vector bundle $W: X\to B$ with a 
$B$-structure, a $B$-structure to the 
vector bundle $D_{n}W$ over $D_{n}X$, and
so also to its restriction 
$\Vreg\otimes W$ to
$B\Sigma_n\times X$.
The Thom spectrum spectrum $M$ is then an $E_\infty$ ring spectrum, 
whose underlying $H_\infty$-structure is such that if
$u_W:X^W\to M$ is the $M$--Thom class of the $B$-bundle $W$, then the
composition
\[
\left(B\Sigma_{n} \times X \right)^{\Vreg\otimes W}\to 
D_{n}M\to M
\]
is the $M$-Thom class of $\Vreg\otimes W$.

\subsection{The spectrum $\MU{0}$}
\label{sec:spectrum-mu0}

A $\BU{0} = \Z\times BU$ bundle over a space $X$ is just a virtual
complex vector bundle $W$, with rank given by the locally constant
function  
\[
    X \xrightarrow{W} \Z\times BU \rightarrow \Z. 
\]
The tautological line bundle $L$ over
$\cp$ gives rise to a natural map 
\[
     \Phi_{\MU{0}} (L):  \cp^{L} \rightarrow \MU{0}.
\]
If $\E$ is an even periodic ring spectrum with formal group
$G=\GpOf{\E}$ and 
\[
g:  \MU{0} \rightarrow \E
\]
is a homotopy multiplicative map, then by Proposition
\ref{t-pr-line-v-unique} the composition 
\[
    (\cp)^{L} \xra{\Phi_{\MU{0}} (L)} \MU{0} \xra{g} \E
\]
is a trivialization $s_{g}$ of the ideal sheaf $\Line (L)\cong \I
(\e)$ over $G$, that is, a coordinate on $G$.  The standard result
about $\MU{0}$-orientations is  

\begin{Lemma} \label{t-le-mu-0-coord}
The assignment $g\mapsto s_{g}$ is a bijection between the set of maps
of homotopy commutative ring spectra $\MU{0}  \rightarrow \E$ and
coordinates on $\GpOf{\E}$.
\end{Lemma}

\begin{proof}
For $MU=\MU{2}$ instead of $\MU{0}$ the standard reference is
\cite{Adams:BlueBook}.  The minor modifications for $\MU{0}$ may be
found in \cite{AHS:ESWGTC}.
\end{proof}

It is customary to express Lemma \ref{t-le-mu-0-coord} in terms of
formal group laws.  A formal group law is the same thing as a 
formal group together with a coordinate:  the equivalence sends a
formal group $G$ over $R$ with multiplication  
\[
    G\times G \xra{m} G,
\]
and coordinate $s\in \O (G)$ to the power series
\[
    m^{*}s \in \O (G\times G)\cong R\psb{s,t}.
\]
We shall write $(G,s)$ for this group law.

For example, the tautological map 
\[
   (\cp)^{L} \rightarrow \MU{0}
\]
gives a coordinate $s_{\MU{0}}$ on $G_{\MU{0}}$. (Quillen's
Theorem \cite{Quillen:fglocct} is that $(G_{\MU{0}}, s_{\MU{0}})$ is the
\emph{universal} formal group law.)

The commutative diagram 
\[
\xymatrix{
{\spf \pihat \E^{\cpplus}}
  \ar[rr]^-{\spf \pihat g^{\cpplus}}
  \ar[d]
& & 
{\spf \pihat \MU{0}^{\cpplus}}
  \ar[d]
\\
{\spf \pihat \E} 
 \ar[rr]^-{\spf \pihat g}
& & 
{\spec \pi_{0}\MU{0}}
}
\]
gives a \emph{relative} map 
\[
    \Bar{g}: G_{\E} = \spf \pihat \E^{\cpplus} \rightarrow (\spf
\pihat g)^{*} G_{\MU{0}}.
\]
Naturality together with the analogous diagram for 
\[
   (\cp \times \cp)_{\plus} \rightarrow \cpplus
\]
shows that $\Bar{g}$ is a homomorphism of formal groups over
$S_{\E}$.  The construction of $s_{g}$ shows that $\Bar{g}$ is an
isomorphism of formal groups, and
\[
     \Bar{g}^{*}s_{\MU{0}} = s_{g}.
\]
In particular, we have the following.

\begin{Lemma} \label{t-le-mu-gp-law} 
If 
\[
     g: \MU{0} \to \E
\]
is a homotopy multiplicative map, then 
\[
    \pi_{0}g: \pi_{0}\MU{0} \rightarrow \pi_{0}\E
\]
classifies the group law $(G_{\E},s_{g})$. \qed
\end{Lemma} 

To understand the Thom class associated to a general virtual complex vector
bundle $W$ over a pointed space $X$, it is convenient as in
\S\ref{sec:cohom-thom-spectra} 
to work first with the bundle $W\otimes L$ over $X\times \cp$, and
then pull back along the identity section.  So let 
$\F=\E^{X_{\plus}}$, and let $f: \E=\E^{S^{0}}\to \F$ be the map associated to
the map $X_{\plus}\rightarrow S^{0}$.

The map 
\[
    (X\times\cp)^{W\otimes L} \rightarrow  \MU{0} \xrightarrow{g} \E
\]
represents a trivialization $s_{W}$ of the line bundle $\Line
(W\otimes L) $ over $\GpOf{\F} = (\spec \pi_{0} f)^{*}\GpOf{\E}.$

\begin{Lemma} \label{t-le-MU-thom-class-translation}
Suppose that $W$ is a line bundle, and let $b\in G (\pi_{0}\F)$ be the
corresponding point.  Under the isomorphism~\eqref{eq:24}
\[
\Line (W\otimes L) \iso T_{b}^{*}\I (\e),
\]
$s_{W}$ is the section 
\[
     s_{W} = T_{b}^{\ast} (\pi_{0}f)^{*} s_{g}
\]
\end{Lemma}

\begin{proof}
The map $\cp\times\cp \rightarrow \cp$ which classifies the tensor
product of line bundles is responsible for the group structure of $G$.
\end{proof}

Now take $X=BA^{*}$ so that $S_{\F}  = (BA^{\ast})_{\E}$, and take
$W=\Vreg$.  Let  
\[
\eps: S_{\F}\rightarrow S_{\E}
\]
be the structural map.  We have a homomorphism 
\[
    A \xra{} \GpOf{\F} 
\]
as in \eqref{eq:10} and an isomorphism of line  bundles over
$S_{\F}$ 
\begin{equation} \label{eq:7}
  \Line (\Vreg\otimes L) \iso
  \bigotimes_{a\in  A} T_{a}^{\ast} \eps^{\ast}\I (\e),
\end{equation}
as in \eqref{eq:12}.  Lemma \ref{t-le-MU-thom-class-translation}
implies the following. 

\begin{Proposition}   \label{t-pr-mu-zero-and-norm}
Under the isomorphism \eqref{eq:7}, we have 
\[
    s_{\Vreg} = \bigotimes_{a\in  A} T_{a}^{\ast} \eps^{\ast} s_{g}.
\]  \qed
\end{Proposition}

\subsection{Comparing the $\hinfty$ structures}

\label{sec:comp-hinfty-struct-mu0}

We continue to suppose that $\E$ is a homogeneous ring spectrum, and
that 
\[
g:   \MU{0}\rightarrow \E
\]
is a map of homotopy commutative ring spectra.  
Now suppose in addition that 
$\pi_{0}\E$ is a 
complete local
ring with perfect residue field of characteristic $p>0$, 
and that the formal group $G=\GpOf{\E}$ is
of finite height.
Let $A$ be a finite
abelian group.  
Proposition~\ref{t-E-BA-ast} implies that, with our hypotheses on
$\E$, 
the natural map~\eqref{eq-hom-a-g-map}
\[
   (BA^{\ast})_{\E} \xra{}\uhom (A,G)
\]
is an isomorphism.  If
\[
   A_{T} \xra{\ell} i^{*}G
\]
is a level structure~\eqref{def-level-structure} with cokernel 
\[
    i^{*} G\xra{\quot} G',
\]
then the homomorphism $\ell$ is classified by a map $\chi_{\ell}$
making 
the diagram 
\[
\xymatrix{
{T}
 \ar[r]^-{\chi_{\ell}}
 \ar[dr]_-{i}
&
{\uhom (A,G)}
 \ar[d]^-{\eps}
\\
&
{S_{\E}}
}
\]
commute. After changing base along $\chi_{\ell}: T \times G
\to \uhom (A,G)\times G$, the isomorphism \eqref{eq:7}
becomes 
\begin{equation} 
  \chi_{\ell}^{\ast}\Line (\Vreg\otimes L) \iso
  \bigotimes_{a\in  A} T_{a}^{\ast} i^{\ast}\I (\e),
\end{equation}
and then Proposition \ref{t-pr-OGp-NOG-new} gives an isomorphism
\begin{equation} \label{eq:62}
      \chi_{\ell}^{\ast}\Line (\Vreg\otimes L)
      \iso
      \quot^{*}\normq_{\quot} i^{*}\I_{G} (\e) \iso \quot^{*}\I_{G'} (\e).
\end{equation}
Proposition \ref{t-pr-mu-zero-and-norm} and Proposition
\ref{t-pr-OGp-NOG-new} have the following 

\begin{Corollary}   \label{t-co-mu-zero-and-norm}
With respect to the isomorphism \eqref{eq:62}, we have 
\[
    \chi_{\ell}^{*} s_{\Vreg} = \quot^{*} \normq_{\quot} i^{\ast}s_{g}.
\]\qed
\end{Corollary}

Now suppose in addition that $\E$ is an $\hinfty$ ring
spectrum.  Using the isogeny  \eqref{t-pr-isogeny-degree}  
\[
    \psilge{\E}: i^{*} G \xra{} \psilestar G, 
\]
equation~\eqref{eq:62} becomes 
\begin{equation} \label{eq:9}
 \chi_{\ell}^{*}\Line(\Vreg\otimes L)
 \iso
 \psilgestar{\E} \I_{\psilestar G}(\e).
\end{equation}
The two ways of going around the diagram 
\[
\begin{CD}
    (BA^{*}\times \cp)^{\Vreg\otimes L} 
    @>>>
    D_{A}\MU{0} @> D_{A}g>> D_{A} \E \\
  @. @VVV @VVV \\
 @.    \MU{0} @> g >> \E
\end{CD}
\]
give two different trivializations $s_\clockwise$ and
$s_{\counterclockwise}$ of $\Line (\Vreg \otimes L)$
over 
\[
(BA^{*}\times\cp)_{\E} = \uhom (A,G)\times G,
\]
and Corollary \ref{t-co-mu-zero-and-norm} shows that 
\[
\chi_{\ell}^{*}s_{\counterclockwise} =  \chi_{\ell}^{*}s_{\Vreg}
= \psilgestar{\E} \normq_{\psilge{\E}} s_{g}.
\]
By definition, 
\[
    \chi_{\ell}^{*} s_{\clockwise} = 
    \psilv L (s_{g}), 
\]
and so with respect to the isomorphism  \eqref{eq:9}, 
Proposition \ref{t-pr-isogeny-degree} gives
\[
  \chi_{\ell}^{*}s_{\clockwise} = 
  \psilgestar{\E} \psilestar s_{g}.
\]
Thus we have the following

\begin{Proposition}\label{t-pr-mu-zero-hinfty}
Let $g:\MU{0}\to\E$ be a homotopy multiplicative map, and
let $s=s_{g}$ be corresponding trivialization of $\I_{G} (\e)$.
If the map $g$ is $H_\infty$, then
for any level structure 
\[
    A \xra{\ell} i^{*}G,
\]
the section $s$ satisfies the
identity
\begin{equation} \label{eq:70}
\normq_{\psilge{\E}} i^{*}s = \psilestar s,
\end{equation}
in which the isogeny $\psilge{\E}$ has been used make the identification
\[
  \normq_{\psilge{\E}}i^{*}\I_{G} (\e) 
  \iso
  \I_{\psilestar G} (\e).
\] \qed
\end{Proposition}

\begin{Remark} \label{rem-2} 
The Proposition can be stated in terms of formal group laws along the
lines of Lemma \ref{t-le-mu-gp-law}.  Given an orientation $g$
and a level 
structure $\ell$ as in the Proposition, we get two ring homomorphisms 
\[
    \alpha,\beta: \pi_{0}\MU{0} \rightarrow \E,
\]
namely 
\[
    \alpha: \pi_{0}\MU{0} \xra{P_{A}} \pi_{0}\MU{0}^{BA^{*}_{\plus}}
            \xrightarrow{\pi_{0}g^{BA^{*}}_{\plus}}
\pi_{0}\E^{BA^{*}_{\plus}} \xrightarrow{\chi_{\ell}} R
\]
and 
\[
    \beta: \pi_{0}\MU{0} \xra{\pi_{0}g} \pi_{0}\E \xra{\psile} R.
\]
The Proposition implies  that $\alpha$ classifies the group law
\spdisplay{(\psilestar G,\normq_{\ell}i^{*}s_{g}),}
while $\beta$ classifies the group law
\spdisplay{
(\psilestar  G,\psilestar s_{g}).
}
\end{Remark}

\begin{Remark} \label{rem-3}
The necessary condition of Proposition
\ref{t-pr-mu-zero-hinfty} was introduced in \cite{Ando:PowerOps},  in
the case that $\E$ is the spectrum associated to the universal
deformation of the  Honda formal group of height $n$.  In that case,
if one has a level structure 
\[
    (\Z/p)^{n} \xra{\ell} i^{*} \GpOf{\E},
\]
one finds that
\begin{align*}
    \psile & = i \\
    \psilge{\E} & = p: \GpOf{\E} \to \GpOf{\E},
\end{align*}
so that equation~\eqref{eq:70} becomes
\[
\normq_{p} i^{*}s =  i^{*}s,     
\]
or after pulling back along $p$, 
\[
    \prod_{a\in (\Z/p)^{n}} T_{a}^{*} s  = p^{*} s.
\]
\end{Remark} 

\section{A necessary condition for an $\MU{2k}$-orientation to be $\hinfty$}
\label{sec:necess-cond-mu2k}

In this section we describe the modifications to Proposition
\ref{t-pr-mu-zero-hinfty} needed in the case that $k\geq 1$ and 
\[
   g: \MU{2k} \to \E
\] 
is a homotopy multiplicative map
from the Thom spectrum of $\BU{2k}$ to $\E$.

Let $bu$ denote connective $K$-theory.  We recall that 
\[
     [X,\BU{2k}] \iso bu^{2k} (X).
\]
This makes it clear 
that if $V$ is a $\BU{2k}$-bundle over $X$ and $W$ is any
(virtual) vector bundle, then $W\otimes V$ has a canonical
$\BU{2k}$-structure.  
Also, as 
the map 
\[
     \cp \xra{} BU = \BU{2}
\]
classifying the reduced tautological bundle $1-L$ may be viewed as an
element of $bu^{2} (\cp)$, 
it follows that the bundle 
\[
 V=(1-L_{1})\otimes\dots\otimes (1-L_{k})
\]
over $(\cp)^k$ has a $\BU{2k}$-structure.

Suppose that $\E$ is an even periodic ring spectrum
and let $G=G_{\E}$.  

\begin{Lemma} \label{t-le-Theta-and-BUn} 
\begin{textList}
\item Proposition \ref{t-pr-line-v} gives an isomorphism
\[
    t_{V}: \Line (V)\iso \Theta^{k} (\I_{G} (\e)).
\]
\item For the bundle $L\otimes V$ over $\cp^{k+1}$,  Proposition
\ref{t-pr-line-v} gives an isomorphism 
\[
    t_{L\otimes V}: \Line (L\otimes V)\iso 
\frac{\hat\mu_{12}^\ast\Line (V)}
     {\hat\pi_{2}^{\ast}\Line (V)}
\]
of line bundles over $\GpOf{\E}^{k+1}$.
\item  In the notation of Lemma \ref{lem-thom-class}, the $\MU{2k}$-Thom class of the
bundle $L\otimes V$ over 
$(\cp)^{k+1}$ is given by
\[
\Phi_{\MU{2k}}(L\otimes V) = 
\frac
{\hat\mu_{12}^\ast\Phi_{\MU{2k}}(V)}
{\hat\pi_{2}^\ast\Phi_{\MU{2k}}(V)}
\]
\end{textList}
\end{Lemma}

\begin{proof} The first part follows from the discussion of the line
bundles $\Line (V)$ in \S\ref{sec:cohom-thom-spectra}.  For the second
two parts,  
simply write
\spmline{
L\otimes V = (1 - LL_1)\otimes(1-L_2)\otimes\dots\otimes (1-L_k)
}
{
-
(1-L)\otimes(1-L_2)\otimes\dots\otimes (1-L_k),
}
and use Lemma \ref{lem-thom-class}.
\end{proof}

The Lemma implies that if
$g:\MU{2k}\to\E$ is a homotopy multiplicative map,
then the composition
\[
(\cp^k)^V\to\MU{2k}\xra{g}\E
\]
represents a trivialization $s=s_{g}$ of $\Theta^{k} (\I_{G}
(\e))$ (In fact it is easily seen to be a $\Theta^{k}$-structure on
$\I_{G} (\e)$ in the sense of \eqref{defn-cubical-structure}).
If $\E$ is an $\hinfty$ ring spectrum, and $A$ is a finite
abelian 
group, then the two ways of going around the diagram 
\begin{equation} \label{eq:38}
\begin{CD}
    (BA^{*}\times (\cp)^{k})^{\Vreg\otimes V} 
    @>>>
    D_{A}\MU{2k} @> D_{A}g>> D_{A} \E \\
  @. @VVV @VVV \\
  @.   \MU{2k} @> g >> \E
\end{CD}
\end{equation}
give two different trivializations $s_\clockwise$ and
$s_{\counterclockwise}$ of 
$
\Line (\Vreg \otimes V)
$
over 
$
\spf \pihat \E^{(BA^{*}\times (\cp)^{k})_{\plus}} = \uhom (A,G)\times G^{k}.
$
The second part of Lemma \ref{t-le-Theta-and-BUn} implies that 
\begin{equation} \label{eq:6}
       \Line (\Vreg\otimes V) \iso \bigotimes_{a\in A}
\rT{a} \Line (V),
\end{equation}
where $\rT{a}$ is translation operation introduced in~\eqref{eq:53-new}.

If 
\[
    A_{T} \xra{\ell} i^{*}G
\]
is a level structure on $G$~\eqref{def-level-structure}, then 
after changing base along the map 
\[
      T\times G^{k} \xra{\chi_{\ell}\times G^{k}} \uhom
(A,G)\times G^{k}
\]
and using the isogeny 
\[
      \psilge{\E}:  i^{*} G \xra{} \psilestar  G
\]
of Proposition \ref{t-pr-isogeny-degree}, 
we have isomorphisms
\splitpageyesno{
\begin{multline} \label{eq:54}
 \chi_{\ell}^{*}\Line(\Vreg\otimes V) \iso
 \psilgestar{\E} \Tilde{\normq}_{\psilge{\E}} i^{*}  \Theta^{k} (
\I_{G}(\e)) \\
 \iso 
 \psilgestar{\E} \Theta^{k} (\I_{\psilestar G} (\e));
\end{multline}
}
{
\begin{equation} \label{eq:54}
 \chi_{\ell}^{*}\Line(\Vreg\otimes V) \iso
 \psilgestar{\E} \Tilde{\normq}_{\psilge{\E}} i^{*}  \Theta^{k} ( \I_{G}(\e))
 \iso 
 \psilgestar{\E} \Theta^{k} (\I_{\psilestar G} (\e));
\end{equation}
}
from~\eqref{eq:6}
and~\eqref{eq-rnorm-idents}.
The third part of Lemma \ref{t-le-Theta-and-BUn} and Definition
\ref{def-rnorm-section}  imply that  with
respect to this isomorphism we have 
\[
\chi_{\ell}^{*}s_{\counterclockwise} = \psilgestar{\E}
\Tilde{\normq}_{\psilge{\E}} i^{*} s_{g}.
\]
By definition, we have 
\[
\chi_{\ell}^{*}s_{\clockwise} = 
\psilv V (s_{g}), 
\]
and with respect to the isomorphism~\eqref{eq:54}, Proposition
\ref{prop-that-damn-isogeny-lemma} gives the equation
\[
 \chi_{\ell}^{*}s_{\clockwise} = 
 \psilgestar{\E} \psilestar s_{g}.
\]

The analogue of Proposition \ref{t-pr-mu-zero-hinfty} is

\begin{Proposition} \label{t-pr-mu-2k-hinfty}
Let $g:\MU{2k}\to\E$ be a homotopy multiplicative map,
and 
$s$ corresponding section of $\Theta^{k} (\I_{G} (\e))$. 
If the map $f$ is $H_\infty$, then
for each level structure
\[
   A \xra{\ell} i^{*}G
\]
the section $s$ satisfies the
identity 
\[
\Tilde{\normq}_{\psilge{\E}}i^{*} s = \psilestar s,
\]
in which the map $\psilge{\E}$ has been used make the identification
\[
  \Tilde{\normq}_{\psilge{\E}}\Theta^{k} (\I_{i^{*}G} (\e)) 
  \iso
  \Theta^{k} (\I_{\psilestar G} (\e))
\]
as in~\eqref{eq:54}. \qed
\end{Proposition}

\section{The necessary condition is sufficient for $k\leq 3$}
\label{sec:necess-cond-suff}

Suppose that $\E$ is an even periodic $\hinfty$ spectrum.  Suppose
that $\pi_{0}\E$ is a $p$-regular admissible local ring with
perfect residue field of characteristic $p$, and the formal group
$G=\GpOf{\E}$ is 
of 
finite height.  Suppose that $k\leq 3$, and 
let $g:\MU{2k}\to\E$ be a homotopy multiplicative map.  Let 
$s = s_{g}$ be the section of $\Theta^{k} (I_{G} (\e))$ as in
\S\ref{sec:necess-cond-mu2k}. 

\begin{Proposition}\label{t-pr-norm-condition-suffices}
The map $g$ is $\hinfty$ if and only if for
each  level structure 
\[
    A \xra{\ell} i^{*}G,
\]
the section $s$ satisfies the
identity
\[
\tilde \normq_{\psilge{\E}} s = \psilestar i^{*}s,
\]
in which, as in Proposition \ref{t-pr-mu-2k-hinfty}, the isogeny
$\psilge{\E}$ has been used make the identification
\[
\Tilde{\normq}_{\psilge{\E}}\Theta^{k}(\I_{i^{*}G} (\e))
\iso
\Theta^{k}(\I_{\psilestar G} (\e)).
\]
\end{Proposition}

\begin{proof}
We must show that, for all $n$, the diagram 
\[
\begin{CD}
    D_{n} \MU{2k} @> D_{n} g>> D_{n}\E \\
   @VVV  @VVV \\
   \MU{2k} @> g >> \E
\end{CD}
\]
commutes.  The hypotheses on $\pi_{0}\E$ 
and the algebra of the $D_{n}$'s together with
the Sylow structure of the symmetric groups  reduce us immediately to
checking that the diagram
\[
\begin{CD}
    D_{A}\MU{2k} @> D_{A}g>> D_{A} \E \\
  @VVV @VVV \\
  \MU{2k} @> g >> \E
\end{CD}
\]
commutes when $A$ is a
Sylow subgroup of $\Sigma_{p}$ \cite[\S 7]{McClure:pohrt}.

Let $g_{\clockwise}$ and $g_{\counterclockwise}$
be the two ways of navigating this diagram.  Each is a
generator of $\pi_{0}\Spectra{U} (D_{A}\MU{2k},\E)$, so by the Thom
isomorphism their ratio is a generator of 
\[
\pi_{0}\Spectra{U}
(D_{A}\BU{2k}_{\plus}, \E).
\]
For $k\leq 3$, the natural map 
\[
    \pi_{0}\Spectra{U} (D_{A}\BU{2k}_{\plus},\E)
    \xrightarrow{\Delta^{*}}
    \pi_{0}\Spectra{U} ((BA^{*}\times \BU{2k})_{\plus},\E)
\]
is injective (see e.g. \cite[7.3]{McClure:pohrt}).  Let
$\F=\E^{\BU{2k}_{\plus}}$. Our hypotheses on $\E$ and the fact that,
for $k\leq 3$, $H_{*} (\BU{2k},\Z)$ is concentrated in even degrees
(for $\BU{6}$ see \cite{Singer:ConnectiveBU} or \cite{AHS:ESWGTC})
imply that   the natural maps induce isomorphisms 
\begin{align*}
   \pi_{0}\E & \iso \O (S_{\E}) \\
   \pi_{0}\F & \iso \O (S_{\F}) \\
   \pi_{0}\Spectra{U} ((BA^{*}\times \BU{2k})_{\plus},\E) & 
    \iso \O (\uhom (A,G_{\F})).
\end{align*}

By Proposition \ref{t-pr-char-map-inj-cyclic}, it suffices to show
that $g_{\clockwise}/g_{\counterclockwise} = 1$  after changing base
along the two maps 
\begin{align*}
    \ulvl (A,G_{\F}) &\rightarrow \uhom (A,G_{\F}) \\
     S_{\F} & \rightarrow \uhom (A,G_{\F})
\end{align*}
classifying respectively the level structure and the zero homomorphism.

After changing base to $\ulvl (A,G)$,
$g_{\clockwise}/g_{\counterclockwise}$ becomes  
\[
\left(\psilfstar s \right) / \left(\tilde \normq_{\psilge{\F}}i^{*} s\right) ,
\]
as in Proposition \ref{t-pr-mu-2k-hinfty}.
The base change $S_{\F} \rightarrow \uhom (A,G)$ corresponds to the
augmentation 
\[
     \pi_{0}\F^{BA^{*}_{\plus}} \rightarrow \pi_{0}\F,
\]
under which each $g$ restricts to the Thom class
\[
      \BU{2k}^{V^{p}} \xra{\Phi_{\MU{2k}} (V^{p})} 
      \MU{2k} \xra{g} \E
\]
of $V^{p}$, where $V$ is the standard bundle over $\BU{2k}$.
\end{proof}

\part{Even periodic cohomology of abelian groups and Thom complexes}
\label{part:even-peri-cohom}

\section{Even cohomology of abelian groups}
\label{sec:even-cohom-abel}

Suppose that $\E$ is a homogeneous ring spectrum,
with formal group $G=G_{\E}$,
and let $A$ be a finite abelian group.  Let $A^{*}$ be
the character group $A^{*} = \hom (A,\C^{\times})$. 
An element $a$ of $A$ may be viewed as a character of $A^{*}$, giving
a line bundle $V_{a}$ over $BA^{*}$ and so a map 
$(BA^\ast)_{\E}\to (\cp)_{\E}= G$, i.e. a $(BA^{\ast})_{\E}$-valued
``point'' of $G$.     As $a$ varies we get a map of sets
\begin{equation}\label{eq:10}
    A \xrightarrow{\chi} G ((BA^{*})_{\E}).
\end{equation} 
Since 
\[
   V_{a+b} = V_{a}\otimes V_{b},
\]
and since the group structure of $G$ comes from the map 
\[
    \cp\times \cp \to \cp
\]
which classifies the tensor product of line bundles, the map $\chi$ is
a group homomorphism, and so it is classified by a map of of formal schemes
\begin{equation}\label{eq-hom-a-g-map}
(BA^{*})_{\E} \xrightarrow{\clchi} \uhom (A,G).
\end{equation}
This map is often an isomorphism.  For example, we have the
\begin{Proposition} \label{t-E-BA-ast}
If $\pi_{0}\E$ is a complete local ring of residue characteristic
$p>0$, and if the height of the formal group $G_{\E}$ is finite, then 
the map $\clchi$ is an isomorphism of formal
schemes over $S_{\E}$. 
\end{Proposition} 

\begin{proof}
This formulation of the $\E$-cohomology of abelian groups appeared in
\cite{HKR:ggc}.
\end{proof}

Suppose that $\clchi$ is an isomorphism, and suppose that we have a
level structure 
\[
     A_{T} \xra{\ell} i^{*}G
\]
over a formal scheme $T$.
The homomorphism $\ell$ is classified by a map $\chi_{\ell}$ making
the diagram 
\begin{equation} \label{eq:25}
\xymatrix{
{T} 
 \ar[r]^-{\chi_{\ell}}
 \ar[dr]_-{i}
&
{(BA^{*})_{\E}}
 \ar[d] \\
& 
{S_{\E}}
}
\end{equation}
commute.  

\begin{Proposition}\label{t-pr-ra-kills-transfers}
If $\pi_{0}\E$ is a complete local ring, and $G$ is of finite height, 
and if ${A'}^\ast\subset
A^\ast$ is a proper subgroup, then the composite 
map of $\pi_0\E$-modules 
\[
\pi_0\E^{{BA'}^\ast_\plus} \xra{\text{transfer}}
\pi_0\E^{{BA}^\ast_\plus}\xra{\chi_{\ell}} \O (T)
\]
is zero.
\end{Proposition}

\begin{proof}
It suffices to consider the case that 
\[
   \ell: A \to i^{*} G
\]
is the tautological level structure over $\ulvl (A,G)$.

If $A$ is not a $p$-group then $\ulvl (A,G)$ is empty and the result
is trivial. 

Suppose that $A'=0$ and $A=\Z/p$.  Let $t\in \pi_{0}\E^{\cpplus}$ be a
coordinate, and let $F$ be the resulting group law.  Then 
\[
  \pi_{0}\E^{BA^{\ast}_{\plus}} \iso \pi_{0}\E \psb{t}/ [p]_{F} (t)
\]
and $\transfer: \pi_{0}\E^{{BA'}^{\ast}_{\plus}}  = \pi_{0}\E
\rightarrow \pi_{0}\E^{BA^{\ast}}$ 
is given by 
\[
    \transfer (1) = \pointy{p} (t)
\]
(see e.g. \cite{Quillen:Ele}), where $\pointy{p} (t)$ is the power
series such that  
\[
      t \pointy{p} (t) = [p]_{F} (t).
\]
The result follows from the
isomorphism~\eqref{eq:41}
\[
    \O (\ulvl (\Z/p,\GpOf{\E})) \iso   \pi_{0}\E\psb{t}/ \pointy{p} (t).
\]

For the general case, we may suppose that ${A'}^{\ast}\subsetneq
A^{\ast}$ is maximal, and so we have a pull-back diagram 
\[
\begin{CD}
B{A'}^{\ast} @> i' >> B0 \\
@V j' VV       @VV j V \\
BA^{\ast}   @> i >> B C^{\ast}
\end{CD}
\]
where $C\subsetneq A$ is cyclic of order $p$.  The commutativity of the diagram
\[
\xymatrix{
{\pi_{0}\E^{B{A'}^{\ast}_{\plus}}}
 \ar[d]_{\transfer}
&
{\pi_{0}\E^{B0_{\plus}}}
 \ar[d]^{\transfer}
 \ar[l]_{{i'}^{*}}
 \ar@(r,r)[dd]^{0}
\\
{\pi_{0}\E^{BA^{\ast}_{\plus}}}
 \ar[d]_{\pi}
&
{\pi_{0}\E^{BC^{*}_{\plus}}}
 \ar[d]
 \ar[l]_{i^{*}}
\\
{\O ( \ulvl (A,G))}
&
{\O (\ulvl (C,G))}
 \ar[l]
}
\]
implies that 
\[
   \pi (\transfer (1)) = 0.
\]
The result follows, since $\pi_{0}\E^{B{A'}^{\ast}_{\plus}}$ is a cyclic
$\pi_{0}\E^{BA^{\ast}_{\plus}}$-module via ${j'}^{*}$, and $\transfer$ is a map of
$\pi_{0}\E^{BA^{\ast}_{\plus}}$-modules.
\end{proof}

\section{Algebraic geometry of the Thom isomorphism} \label{sec:cohom-thom-spectra}

Suppose that $X$ is a space, and that $V$ is a complex vector bundle
over $X$.  The $\pi_0 \E^{X_\plus}$-module $\pi_0 \Spectra{U}
(X^{V},\E)$ is free of rank one 
(since $\E$ is complex orientable) and so can be interpreted as the
module of sections of a line bundle $\Line(V)$ over $X_{\E}$.
The fact that the Thom complex of an external Whitney sum is the smash
product of the Thom complexes gives rise to a canonical isomorphism
\begin{equation}\label{eq-line-tensor} 
\Line(V\oplus W)\iso \Line(V)\otimes\Line(W)
\end{equation}
This property can then be used to extend the definition of
$\Line(V)$ to virtual bundles; we define 
\begin{equation} \label{eq:20}
\Line(V-W) = \Line(V)\otimes \Line(W)^{-1}.
\end{equation}
If $f:X\to Y$ is a map, and $V$ is a virtual bundle over $Y$, then
there is an isomorphism
\[
\pi_0\Spectra U(X^{f^\ast
V},\E)\iso\left(\pi_0\E^{X_\plus}\right)
\otimes
\pi_0\Spectra U(Y^V,\E).
\]
In terms of algebraic geometry, this means that there is a natural
isomorphism
\begin{equation}\label{eq-line-pullback}
\Line(f^\ast V)\iso \left(f_{\E}\right)^\ast\Line(V).
\end{equation}

Here is a series of examples which lead to a fairly complete
understanding of the functor $\Line.$  

\begin{textList}
\item If $L$ denotes the canonical line bundle
over $\cp$, then the zero
section identifies $\pi_0\Spectra{U} (\cp^{L},\E)$ with the augmentation ideal in
$\pi_0 \E^{\cpplus}$, and so we have an isomorphism 
\begin{equation} \label{eq:16}
\Line (L)\iso \I (\e).
\end{equation}

\item\label{item-example-2} Suppose that $V$ is a line bundle over $X$, classified by a map
$b:X\to\cp$.  In terms of algebraic geometry, the map $b$ defines an
$X_{E}$-value point $b=b_{\E}$ of $G$. It follows
from~\eqref{eq-line-pullback} that 
\begin{equation} \label{eq:19}
\Line(V) \iso b^\ast\I(\e)\iso \e^\ast \I (-b).
\end{equation}

\item Taking $X$ to be a point and $V$ to be the trivial complex line
bundle in (\ref{item-example-2}), we have 
\begin{equation} \label{eq:22}
\Line(V) \iso \e^\ast \I (\e).
\end{equation}
Now $\Line (V)$ is the sheaf associated to $\pi_{2}\E$, while 
$\e^{\ast} \I (\e)$ is the sheaf of cotangent vectors at the 
origin of $G$, isomorphic to the sheaf $\omega_{G}$ of invariant
differentials on $G$.  
\item If $V$ is the trivial bundle of dimension $k$, then by
\eqref{eq:22} and \eqref{eq-line-tensor}, $\Line(V)$ is just
$\omega_{G}^k$.   If $f:\E\to \F$ is an $\E$-algebra
(e.g. $\F=\E^{X_{\plus}}$), this gives an interpretation of the
homotopy group $\pi_{2k}\F$ as the sections of $f^{*}\omega_{G}^{k}$.
\item If $V= (1-L)$ is the reduced canonical line bundle over $\cp$, then
using \eqref{eq:20}, \eqref{eq:16}, and \eqref{eq:22} we have
\[
 \Line (V) \iso \struc^{*} \e^{*}\I (\e)\otimes \I (\e)^{-1} = \Theta^{1} (\I (\e)),
\]
where $\struc: \GpOf{\E}\to S_{\E}$ is the structural map and
$\Theta^{1}$ is defined in Definition~\ref{defn-thetat}.
\item With the notation of example~(\ref{item-example-2}) consider
the line bundle $V\otimes L$ over $X\times\cp$.  Then
$\Line(V\otimes L)$ is pulled back from $\I_{G}(\e)$ along the map
\[
X_{\E}\times G \xra{b \times 1} G\times G \xra{\mu} G.
\]
It follows that 
\begin{equation} \label{eq:24}
\Line(V\otimes L) \iso T_b^\ast\I(\e)=\I_{X_{\E}\times G}(-b).
\end{equation}
\item More generally, suppose that $V=\sum n_iL_i$ is a virtual sum of
line bundles over $X$.  The line bundles $L_i$ define points $b_i$ of
$G$ over $X_{\E}$, and the bundle $V$ determines the divisor $D=\sum n_i \divisor{b_i}$.
It follows using~\eqref{eq-line-tensor} that
\begin{equation} \label{eq:23}
\Line(V\otimes L) = \I_{X_{\E}\times G}(D^{-1}),
\end{equation}
where $D^{-1}=\sum n_i \divisor{b_i^{-1}}$.
\item In fact, by the splitting principle, the line bundle
$\Line(V\otimes L)$ can be computed in this manner even when $V$ is
not a virtual sum of line bundles.  Indeed,  by the splitting
principle, there is a map $f:F\to X$ with the properties that
$f_{\E}$ is finite and faithfully flat, and  
$f^\ast V$ is a virtual sum of line bundles.  The
line bundle $\Line(f^\ast(V)\otimes L)$ can then be computed as
$\O(D^{-1})$ as above.  But the divisor
$D$ descends to $X_{\E}\times G$, even if none of its points do.
\item Let $A$ be a finite abelian group.  An element $a\in A$ can be
regarded as a character of $A^{\ast}$.  Let $V_a$
be the associated line bundle over $BA^\ast$.
Recall~\eqref{eq:10}
that this construction defines a group homomorphism
\[
    \chi: A \rightarrow G (BA^{*}_{\E}).
\]
The line bundle
$\Line(V_a\otimes V\otimes L)$ over
$BA^{*}_{\E}\times X_{\E}\times G$ is
\[
\Line (V_{a}\otimes V\otimes L) \iso T_{a}^{\ast}\I(D^{-1});
\]
taking $V$ to be the trivial line bundle over a point gives
\[
  \Line (V_{a}\otimes L)\iso T_{a}^{*}\I (\e) = I (a^{-1})
\]
\item\label{item-o-a} Now let 
\[
V_{\text{reg}}=\bigoplus_{a\in A}V_a
\]
be the regular representation of $A^\ast$.  
Over the scheme
$(BA^{*})_{\E}\times G$, the line bundle associated to the Thom
complex of $\Vreg \otimes V \otimes L$ is  
\begin{equation} \label{eq:36}
\Line (\Vreg \otimes V \otimes L) \iso 
\bigotimes_{a\in A}T_{a}^{*}\I (D^{-1}) 
\iso \I \left( \sum_{a} T_{a}^{*} D^{-1}\right).
\end{equation}
In particular,
\begin{equation} \label{eq:12}
\Line (\Vreg \otimes L) \iso \bigotimes_{a\in A}T_{a}^{*}\I (\e)
                        \iso \I (\chi).
\end{equation}
\item Suppose that the map 
\[
    \clchi : (BA^{*})_{\E} \rightarrow \uhom (A,G) 
\]
of \eqref{eq-hom-a-g-map} is an isomorphism.  If 
\[
   A_{T} \xra{\ell} i^{*}G \xra{\quot} G'
\]
is a level structure with cokernel $\quot$ over $T$,
then changing base 
in \eqref{eq:12} along 
\[
    T \times G \xra{\chi_{\ell}} \uhom (A,G) \times G
\]
(where $\chi_{\ell}$ is the map classifying the homomorphism $\ell$;
see \eqref{eq:25}) and using \eqref{eq:84} gives
\begin{equation} \label{eq:37}
\chi_{\ell}^{*}\Line(V_{\text{reg}}\otimes L) \iso 
\quot^{*}\normq_{\quot} \I_{G} (\e)\iso
\quot^\ast\I_{G'}(\e) \iso
\I_{G} (\ell).
\end{equation}
\item Restricting the above example to $BA^\ast$ 
we find that
\begin{align*}
\chi_{\ell}^{*}\Line (\Vreg) =
\e_G^\ast\quot^\ast\I_{G'}(\e) & = \e_{G'}^\ast\I_{G'}(\e)
\\
 &= \omega_{G'}.
\end{align*}
\end{textList}

This series of examples establishes the following results:

\begin{Proposition}\label{t-pr-line-v}
For a pointed topological space $X$, let $\F$ be the spectrum $\E^{X_{\plus}}$,
and let 
\[
G =G_{\F} = (\cp)_{\F}
\]
be the associated formal group. 
Attached to each (virtual) complex vector bundle $V$ over  $X$ is a
divisor $D=D_V$ on $G$, and an isomorphism
\begin{equation} \label{eq-line-od-iso}
t_V:\Line(V\otimes L)\iso \I_{G}(D^{-1}).
\end{equation}
The map $t_V$ restricts to an isomorphism
\[
t_V:\Line(V)\iso \e^\ast\I(D^{-1}).\qed
\]
\end{Proposition}

\begin{Proposition}\label{t-pr-line-v-unique}
The correspondence $V\mapsto D_V$ and the
isomorphism~\eqref{eq-line-od-iso} are determined by the following
properties
\begin{thmList}
\item  If $V=V_1\oplus V_2$, then $D_V$=$D_{V_1}+D_{V_2}$, and with
the identifications
\begin{align*}
\Line(V_1)\otimes\Line(V_2) &\iso\Line(V) \\
\I(D_{V_1}^{-1})\otimes\I(D_{V_2}^{-1})& \iso\I(D_V^{-1}),
\end{align*}
 there is an equality
\[
t_{V}=t_{V_1}\otimes t_{V_2}.
\]
\item If $f:Y\to X$ is a map of pointed spaces, and if $W=f^\ast V$,
then $D_{W}=f^\ast D_V$, and $t_W=f^\ast t_V$.
\item If $X$ is a point, and $V$ has dimension 1, then $D=\divisor{\e}$, 
and the isomorphism
\begin{equation} \label{eq:26}
t_L:\Line(L)\iso\I (\e)
\end{equation}
is given by applying $\pi_0\E^{(\slot)}$ to the 
zero section $\cpplus\to\cp^L$.\qed
\end{thmList}
\end{Proposition}

\part{Level structures and isogenies of formal groups}
\label{part:level-struct-isog}

\section{Level structures} \label{sec:level-structures}

\subsection{Homomorphisms}

Suppose that $A$ is a finite abelian group and $G$ is a formal group
over a formal scheme $S$.

\begin{Definition} \label{def-uhom}
We write $\uhom (A,G)$ for the functor from formal schemes 
to groups defined by the formula 
\splitpageyesno{
\begin{align*}
 \uhom (A,G) (T) = 
\{ \text{pairs } (u,\ell) \suchthat 
   & u: T \to S\;, \\
   &  \ell \in \hom (A,u^{*}G (T))\}.
\end{align*}
}
{
\[
 \uhom (A,G) (T) = 
\{ \text{pairs } (u,\ell) \suchthat 
   u: T \to S\;,\; \ell \in \hom (A,u^{*}G (T))\}.
\]
}
\end{Definition}

\begin{Remark} \label{rem-4}
We shall use the notation 
\[ 
   A_{T} \xra{\ell} u^{*}G
\]
to indicate that $T$ is a formal scheme and $(u,\ell) \in \uhom (A,G) (T)$.
\end{Remark}

\begin{Example} \label{ex-hom-zn-G}
Let $G$ be a formal group over $R$, and suppose that $x$
is a coordinate on $G$.  Let $F$ be the resulting group law.  The
``$n$-series'' of $F$ is the  power series $[n] (t) \in R\psb{t}$
defined by the formula 
\[
     [n] (x) = n^{*}x, 
\]
where the right-hand-side refers to the pull-back of functions along
the homomorphism 
$
     n: G \xra{} G.
$
To give a homomorphism 
\[
     \ell: \Z/n \rightarrow G (T)
\]
is to give a topologically nilpotent element $x (\ell (1))$ of $\O (T)$,
with the property that 
\[
    [n] (x (\ell (1))) = 0;
\]
the homomorphism $\ell$ is then given by 
\[
     x (\ell (j)) = [j] (x (\ell (1))).
\]
It follows that 
\[
    \uhom (\Z/n,G) = \spf  
    \Bigl( R\psb{x (\ell (1)) }/ 
    \bigl( [n] (x (\ell (1))) \bigr) \Bigr).
\]
\end{Example}

It is clear from the definition that 
if $B\subseteq A$ then there is a
restriction map 
\[
    \uhom (A,G) \rightarrow \uhom (B,G), 
\]
and if $A=B\times C$ then the resulting map 
\begin{equation}  \label{eq-uhom-prod-prod}
    \uhom (A,G) \rightarrow \uhom (B,G) \times_{S} \uhom (C,G)
\end{equation}
is an isomorphism.  Also from the definition we see that 
if $j: S'\to S$ is a map of 
formal schemes, then the natural map 
\[
     \uhom (A,j^{*}G) \rightarrow j^{*}\uhom (A,G)
\]
is an isomorphism.  Combining these observations with Example
\ref{ex-hom-zn-G} and the structure of finite abelian groups gives the
following.

\begin{Lemma}\label{t-hom-rep}
The functor $\uhom (A,G)$  is represented by an affine formal scheme
over $S$.    If $j: S'\to S$ is a map of 
formal schemes, then the natural map 
\[
     \uhom (A,j^{*}G) \rightarrow j^{*}\uhom (A,G)
\]
is an isomorphism of formal schemes over $S'$. \qed
\end{Lemma} 

For formal groups over $p$-local rings, only the $p$-groups
give anything interesting.

\begin{Example} \label{ex-hom-zn-G-n-unit}
Returning to Example \ref{ex-hom-zn-G}, the $n$-series is easily seen
to be of the form  
\[
     [n] (t) = nt + o (2).
\]
If $n$ is a unit in $R$, then 
\[
    R\psb{x}/ ([n] (x)) \iso R
\]
so $\uhom (\Z/n,G)$ is the trivial group scheme over  $R$.
\end{Example}

\begin{Example}\label{ex-hom-zn-G-p-local}
If $R$ is a complete 
local ring of residue
characteristic $p$,
then there is an $h$ with $1\leq h\leq \infty$ such that 
\[
       [p^{m}] (t) \equiv \genericunit t^{p^{mh}} + o (t^{p^{h}}+1)
\mod \maxidof{R}. 
\]
This $h$ is
called the \emph{height} of $G$.  If $h$ is finite, then the
Weierstrass Preparation Theorem \cite[pp. 129--131]{Lang:cf} implies
that there are 
monic polynomials $g_{m} (t)$ of degree $p^{mh}$ such that
\[
      [p^{m}] (t) = g_{m} (t)\cdot \genericunit,
\]
where $\genericunit$ is a unit of $R\psb{t}$.
It follows that $\O (\uhom (\Z/p^{m},G))$ is finite and free of rank
$p^{hm}$ over $R$.   
\end{Example}

These examples generalize to give the following.

\begin{Proposition} \label{t-pr-hom-finite-flat}
Let $G$ be a formal group of finite height over a local 
formal 
scheme $S$.  Then $\uhom (A,G)$ is a
local formal scheme 
over $S$.  For $B\subseteq 
A$, the forgetful map 
\[
    \uhom (A,G) \rightarrow \uhom (B,G)
\]
is a map of formal schemes, finite and free of rank $d^{h}$, where $d$
is the order of the $p$-torsion subgroup of $A/B$. \qed
\end{Proposition}




\subsection{Level structures}

The scheme $\uhom (A,G)$ has an important closed subscheme $\ulvl
(A,G)$, which was introduced by Drinfel'd \cite{Drinfeld:EM}.  Suppose
that $G$ is a formal group 
over a formal scheme $S$.  For simplicity, we suppose that $S$ is a
local formal scheme of residue characteristic $p>0$.

\begin{Definition} \label{def-level-structure} 
Let $T$ be a formal scheme.  A $T$-valued point 
\[
     A_{T} \xra{\ell} i^{*}G
\]
of $\uhom (A,G)$ is a \emph{level structure} if for each prime $q$
dividing $|A|$,  the subgroup $i^{*}G[q] = \Ker (q: i^{*}G\to i^{*}G)$
is a divisor 
on $G/T$, and there is an inequality of 
divisors
\[
\sum_{\substack{a\in A \\ q a = 0}}\divisor{\ell (a)} \leq i^{*}G[q]
\]
in $i^{*}G$.  The subfunctor of
$\uhom (A,G)$ consisting of 
level structures will be denoted $\ulvl (A,G)$. 
\end{Definition}

\begin{Remark}\label{rem-1} 
If we say that 
\[
   A_{T} \xra{\ell} i^{*} G
\]
``is a level structure,'' we mean that $T$ is a formal scheme, and
$(i,\ell)$ is a $T$-valued point of $\ulvl (A,G)$. We may omit one of $T$
and $i$ if it is clear from the context.
\end{Remark}

Here are some examples to give a feel for level structures.  First of
all, only $p$-groups of small rank can produce level structures.  

\begin{Lemma}
If
$|A|$ is not a power of $p$, then 
\[
     \ulvl (A,G) = \emptyset.
\]
If the height of $G$ is $h$ and the $p$-rank of $A$ is greater
than $h$, then again $\ulvl (A,G)=\emptyset$.
\end{Lemma}

\begin{proof}
If $|A|$ is not a power of $p$, then there is a prime $q\neq p$
such that the divisor
\[
     \sum_{qa=0} \divisor{\ell (a)}
\]
has degree greater than $1$.  However, $q: G\to G$ is an
isomorphism, so $G[q] = \divisor{\e}$ has degree $1$.
Similarly, if the height of $G$ is $h$ then the degree of
$G[p]$ is $p^{h}$.
\end{proof}

A level structure is trying to be a monomorphism; for example if $R$
is a domain in which $|A|\neq 0$, then a 
homomorphism 
\[
    \ell: A\to G (R)
\]
is a level structure if and only if it is a monomorphism (Corollary
\ref{t-co-lvl-struc-mono-p-not-zero-divisor-domain}).
However, naive monomorphisms from $A$ to $G$ can't in general be
a representable functor.

\begin{Example}\label{ex-char-p-lvl-not-mono}
Let $\Gmh$ be the formal multiplicative group with coordinate $x$ so
that the group law is 
\[
     x\fs{F} y = x + y - xy.
\]
The $p$-series is 
\[
    [p] (x) = 1 - ( 1- x)^{p}.
\]
The monomorphism
\[
       \Z/p \xrightarrow{} \Gmh (\Z\psb{y} / [p] (y))
\]
given by $j\mapsto [j] (y)$ becomes the zero map under the base change
\begin{align*}
\Z\psb{y}/ ( [p] (y))& \rightarrow \Z/p\\
       y &\mapsto 0.
\end{align*}
\end{Example}

On the other hand, the functor $\ulvl (A,G)$ is representable.

\begin{Lemma} \label{t-le-ulvl-represented}
Let $G$ be a formal group of finite height over a local formal scheme
$S$, and let $A$ be a finite abelian group.  The functor $\ulvl (A,G)$
is a closed formal subscheme of $\uhom (A,G)$.
\end{Lemma}

\begin{proof}
See Katz and Mazur \cite[1.3.4]{KaMa:AMEC} or
\cite{Strickland:FiniteSubgps}.
\end{proof}

It is clear from the definition that if $j: S'\to S$ is a map of 
formal schemes, then the natural map 
\[
     \ulvl (A,j^{*}G) \rightarrow j^{*}\ulvl (A,G)
\]
is an isomorphism of formal schemes over $S'$.

\begin{Proposition} \label{t-pr-level-A-finite-flat} 
Suppose that $G$ is a formal group
of finite height $h$ over a local formal scheme $S$ with perfect
residue field of 
characteristic $p>0$, and suppose that $A$ is a finite abelian 
$p$-group with $|A[p]|\leq p^{h}$.  Then we have the following.
\begin{thmList}
\item \label{item-ulvl-fin-flat} The functor $\ulvl (A,G)$ is
represented by a local formal scheme which is finite and 
flat  over $S$: indeed $\O (\ulvl (A,G))$ is a finite free $\O (S)$-module.
\item \label{item-ulvl-reg-local-ring}
If $G$ is the universal deformation of a formal group over a perfect field
(see \S\ref{sec:lubin-tate-groups}) then $\ulvl (A,G)$ is the formal
spectrum of a Noetherian complete local 
domain
which is regular of
dimension $h$. 
\end{thmList}
\end{Proposition}

\begin{proof}
With our hypotheses, we may suppose that $G/S$ is the universal
deformation of a formal group of height $h$ over a perfect field $k$
of characteristic $p$; the general case follows by change of base.

If $A=A[p]$, then the result is precisely the Lemma of
\cite[p. 572, in proof of Prop. 4.3]{Drinfeld:EM}.
The proof in the general case follows similar lines and is given in 
\cite{Strickland:FiniteSubgps}.  

The proof of \ref{item-ulvl-fin-flat}) for a general $A$ 
can be given easily: by definition of $\ulvl (A,G)$, the diagram 
\[
\begin{CD}
\ulvl (A,G) @> j >> \uhom (A,G) \\
@V i VV              @VV k V \\
\ulvl (A[p],G) @> l >> \uhom (A[p],G).
\end{CD}
\]
is a pull-back.  Proposition \ref{t-pr-hom-finite-flat} implies that
$k$ is finite and free, and so $i$ is too.
\end{proof}

\subsection{Level structures over $p$-regular schemes}

In this section, we suppose that $G$ is a formal group of finite
height over  a complete
local ring $E$ of residue characteristic $p>0$.
The following description of
the subscheme $\ulvl (A,G)$ was found by Hopkins in the course of
his work on \cite{HKR:ggc}.

\begin{Proposition} \label{t-pr-level-str-p-neq-0}
Suppose that $p$ is not a zero
divisor in $E$.  Let $x$ be a 
coordinate on $G$.  The scheme $\ulvl (A,G)$ is
the closed subscheme of $\uhom (A,G)$ 
defined by the ideal of annihilators of $x (\ell (a))$,
where $a$ ranges over the non-zero elements of $A[p]$.
\end{Proposition}

The proof will be given at the end of this section.  Note that
the ideal  in the Proposition is independent of the coordinate used
to describe it.

For $n\geq 1$ let $A[n]$ denote the $n$-torsion in $A$.  
Let $R$ be a complete local $E$-algebra, and consider the following
conditions  on a homomorphism
\[
     \ell: A\to G (R).
\]
Again, they are phrased in terms
of a choice of a coordinate $x$ on $G$, but they are easily seen to be
independent of that choice. 

\begin{enumerate}
\item [(A)]
If $0\neq a\in A[p]$  then $x (\ell (a))$ is regular (i.e. not a
divisor of zero).
\item  [(B)]
If $0\neq a\in A[p]$ then $x (\ell (a))$ divides $p$.
\item [(C)] 
$\prod_{ a\in A[p]} 
        (x - x (\ell (a)))$ divides $[p] (x)$.
\item [(D)] The natural map 
\begin{equation} \label{eq-Vandermonde}
   R\psb{x}\left/ \left(\prod_{a\in A[p]} (x-x (\ell (a))) \right)\right.
   \rightarrow 
  \prod_{a\in A[p]} \left( R\psb{x}/ (x-x (\ell (a))) \right)
\end{equation}
is a monomorphism.
\end{enumerate}
Condition (C) says precisely that there is an inequality of
Cartier divisors 
\[
\sum_{pa = 0} \divisor{\ell (a)}
 \leq G[p].
\]
Thus condition (C) is that $\ell$ is a level structure.

\begin{Proposition} \label{t-pr-level-A-char-0}
If $R$ is $p$-torsion free, then these conditions are equivalent.
\end{Proposition}

First we prove the following result.  It will be convenient to use the
symbol $\epsilon$ to denote a generic unit.  Its value may change from
line to line.

\begin{Lemma} \label{t-le-discriminant}
Let $n
= |A[p^{m}]|$.  The discriminant of the set 
\[
\{ x (\ell (a))\suchthat a\in A[p^{m}] \}
\]
is
\[
    \Delta = \epsilon \prod_{0\neq a\in A[p^{m}]} x (\ell (a))^{n}.
\]
\end{Lemma}

\begin{proof}
Let $F$ be the group law associated to a coordinate on $G$.  The formula 
\[
   x \fm{F} y = (x-y) \epsilon (x,y),
\]
where $\epsilon (x,y)\in E\psb{x,y}^{\times}$, gives 
\begin{align*}
     \Delta & = \prod_{a\neq b\in A[p^{m}]} (x (\ell (a)) - x (\ell (b))) \\
            & = \epsilon \prod (x (\ell (a)) \fm{F} x (\ell (b))) \\
            & = \epsilon \prod x (\ell (a) - \ell (b)) \\
            & = \epsilon \prod_{c\neq 0} \prod_{a-b=c} x (\ell (c)) \\
            & = \epsilon \prod_{c\neq 0} x (\ell (c))^{n}.
\end{align*}
\end{proof}

\begin{proof}[Proof of Proposition \ref{t-pr-level-A-char-0}]
Under the hypothesis that $p$ is regular in $R$, it is clear that (B)
implies (A).  Let's check that (A) implies (B).  Note that 
\[
     [p] (x) = x ( p + x e (x))
\]
for some $e (x)\in E\psb{x}$.  For $a\in
A[p]$ we have 
\[
    0 = [p] (x (\ell (a)))  = x (\ell (a)) (p + x (\ell (a)) e (x
(\ell (a)))). 
\]
If $x (\ell (a))$ is not a zero-divisor
in $R$, then we must have 
\[
      p =  - x (\ell (a)) e (x (\ell (a))).
\]

Next, let check that (C) implies (B).  If (C) holds,
then there is a power series $e (x)\in E\psb{x}$ such that 
\[
e (x) \prod_{a\in A[p]} (x-x (\ell (a))) = [p] (x) = p x + o (x^{2})
\]
The coefficient of $x$ on the left is (up to a sign)
\[
    e (0) \prod_{0\neq a\in A[p]} x (\ell (a))
\]
so (B) holds.

Next let's check that (A) implies (D).   
With respect to the basis of powers of $x$ in the domain and the
obvious basis in the range, the matrix of \eqref{eq-Vandermonde}
is the Vandermonde matrix on the set $x (\ell (A[p]))$.  Condition
(A) and Lemma \ref{t-le-discriminant} together imply that
\eqref{eq-Vandermonde} is  a monomorphism.  

Finally, let's check that (D) implies (C).  Each $x (\ell (a))$ is a
root of $[p] (x)$, so the image of $[p] (x)$ in the range of
\eqref{eq-Vandermonde} is zero.  If (D) holds then $[p] (x)$  is zero
in the domain, which implies (C). 
\end{proof}

\begin{Lemma} \label{t-le-statement-a-p-a-p-m}
Condition (A) holds if and only if, for \emph{all} non-zero $a\in A$, $x (\ell
(a))$ is a regular element of $R$.  Condition (B) holds if and only
if, for \emph{all} non-zero $a\in A$, $x (\ell (a))$ divides a power
of $p$.  If $p$ is regular in $R$, then (C) holds if and
only if, for each $m$, $\prod_{a\in A[p^{m}]} (x - x (\ell (a)))$
divides $[p^{m}] (x)$, and (D) holds if and only if for each $m$, the natural
map 
\begin{equation} \label{eq:82}
   R\psb{x}\left/ \left(\prod_{a\in A[p^{m}]} (x-x (\ell (a))) \right)\right.
   \rightarrow 
  \prod_{a\in A[p^{m}]} \left( R\psb{x}/ (x-x (\ell (a))) \right)
\end{equation}
is a monomorphism.
\end{Lemma}

\begin{proof}
Recall that the $p$-series $[p] (x)$ is divisible by $x$: let
$\pointy{p} (x)$ be the power series such that 
\[
    [p] (x) = x \pointy{p} (x).
\]
Thus 
\[
      x (\ell (pa)) = x (\ell (a))\pointy{p} (x (\ell (a))).
\]
so if $x (\ell (pa))$ divides zero (resp. a power of $p$), then so
does $x (\ell (a))$; this proves the statement about (A) and (B).  For
the statement about (D), suppose that $p$ is not a zero divisor in
$R$, and condition (D) holds; by Proposition
\ref{t-pr-level-A-char-0}, condition (A) holds.
With respect to the basis of powers of $x$ in the domain and the
obvious basis in the range, the matrix of~\eqref{eq:82}
is the Vandermonde matrix on the set $x (\ell (A[p^{m}]))$.  
Lemma \ref{t-le-discriminant} and the statement about (A) prove the
statement about (D).  For the statement about (C), suppose that $p$ is
not a zero divisor in $R$ and condition (C) holds; by Proposition
\ref{t-pr-level-A-char-0}, condition (D) holds.  It follows that for
each $m$, the natural map \eqref{eq:82} is a monomorphism.  
If $p^{m} a = 0$, then $x (\ell (a))$ is a
root of $[p^{m}] (x)$, so the image of $[p^{m}] (x)$ in the range of
\eqref{eq:82} is zero, which implies the statement about (C).
\end{proof}

\begin{Corollary}
\label{t-co-lvl-struc-mono-p-not-zero-divisor-domain}
If $R$ is a domain of characteristic $0$, then the 
conditions (A)---(C) hold if and only  if
\spdisplay{
\ell: A\to G (R)
}
is a
monomorphism.   \qed
\end{Corollary}

\begin{proof}[Proof of Proposition \ref{t-pr-level-str-p-neq-0}]
By Proposition \ref{t-pr-level-A-finite-flat},
\spdisplay{
R=\O (\ulvl (A,G))
}
is a finite free $E$-module.
It follows that $p$ is not a zero divisor in $R$.
Proposition \ref{t-pr-level-A-char-0} implies that $R$ is initial
among complete local $E$-algebras satisfying (A).
\end{proof}

\begin{Example} \label{ex-lvl-zp-G}
Let $G$ be a 
formal group of finite height over a
$p$-regular complete 
local ring $R$ of residue characteristic $p$.  
Suppose that $x$ is a coordinate  
on $G$.  Let $\pointy{p} (t) \in R\psb{t}$ be
the power series such that 
\[
     t \pointy{p} (t)  = [p] (t).
\]
Proposition \ref{t-pr-level-str-p-neq-0} implies that 
\begin{equation} \label{eq:41}
     \ulvl (\Z/p,G) \iso 
\spf \left(  
     R\psb{x (\ell (1))} / \pointy{p} \bigl(x (\ell (1))\bigr) 
     \right).
\end{equation}
This calculation occurs as part of the proof of the Lemma in the proof
of Proposition 4.3 of \cite{Drinfeld:EM}.
\end{Example}

\subsection{Calculations in $\uhom (\Z/p,G)$ via level structures}

Let $G$ be a 
formal group of
finite height over a  
$p$-regular complete local ring $R$ with perfect residue field of
characteristic $p>0$. Let $A$ be a finite abelian group. 
By construction there is a natural map 
\[
      \O (\uhom (A,G)) \rightarrow \O (\ulvl (A,G)).
\]
There is also a ring homomorphism 
\[
     \O (\uhom (A,G)) \rightarrow R
\]
classifying the zero homomorphism.  The proof of Proposition
\ref{t-pr-norm-condition-suffices} uses the following result.

\begin{Proposition} \label{t-pr-char-map-inj-cyclic}
The natural map 
\begin{equation} \label{eq:15}
    \O (\uhom (\Z/p,G)) \rightarrow R \times \O (\ulvl (\Z/p,G))
\end{equation}
is injective.
\end{Proposition}

\begin{Remark}  This result is equivalent to the injectivity for the
group $\Z/p$ of the character map of Hopkins-Kuhn-Ravenel, and as such is
proved in \cite{HKR:ggc}.
\end{Remark}

\begin{proof}
Let $h$ be the height of $G$.  
Let $\Lambda= (\Z/p)^{h}$.  Let $g (x)$ be the monic polynomial of
degree $p^{h}$ such that 
\[
 [p] (x) = g (x) \genericunit 
\]
where $\genericunit\in R\psb{x}^{\times}$.    Then 
\begin{equation}
   \uhom (\Z/p,G) = \spf R\psb{x}/ [p] (x) \iso \spf R\psb{x}/ g (x), 
\end{equation}
and 
Proposition \ref{t-pr-level-str-p-neq-0} (see Example
\ref{ex-lvl-zp-G})  implies that 
\[
    \O (\ulvl (\Z/p,G))  =  R\psb{x} / \pointy{p} (x).
\]

Let
\[
     D = \O (\ulvl (\Lambda,G))
\]
and let 
\[
    \ell : \Lambda \rightarrow G (D)
\]
be the tautological homomorphism.
By definition, 
\begin{equation} 
     D = \O (\uhom (\Lambda,G)) / J,
\end{equation}
where $J$ is the ideal obtained by
equating coefficients in
\begin{equation} 
 \prod_{a\in \Lambda} (x- x (\ell (a))) = g (x).
\end{equation}
By Proposition \ref{t-pr-level-A-finite-flat}, $D$ is finite and
free over $R$.  Therefore, letting $u$ denote the
map \eqref{eq:15}, it suffices to show that $D\hot u$ is
injective.

Each non-zero $a\in \Lambda$ gives a monomorphism $\Z/p
\hookrightarrow \Lambda$ and so a homomorphism
\begin{align*}
    \O (\ulvl (\Z/p,G)) & \xra{a} D \\
               x& \mapsto x (\ell (a)).
\end{align*}
We may view these all together as a ring homomorphism 
\[
   D\hot \O (\ulvl (\Z/p,G)) \xra{M} \prod_{0\neq a\in \Lambda} D\psb{x}/
(x-x (\ell (a))). 
\]
Note that the identity map of $D$ may be written as 
\[
    D \xra{F} D\psb{x}/ (x-x (\ell (0)))  = D.
\]

With this notation, the diagram 
\[
\begin{CD}
D\hot \O (\uhom (\Z/p,G)) @> D\hot u >>
 D\hot (R\times \O (\ulvl (\Z/p,G))) \\
@V \iso VV                              @VV F \times M V \\
D \psb{x}
/ \left(\prod_{a\in \Lambda} (x - x (\ell (a)))\right)
     @>>>
\prod_{a\in \Lambda} D\psb{x} / (x-x (\ell (a)))
\end{CD}
\]
commutes, where the map across the bottom is the evident map
\eqref{eq-Vandermonde}.  It is a monomorphism by Proposition   
\ref{t-pr-level-A-char-0}.
\end{proof}

\section{Isogenies} \label{sec:isogenies}

Throughout this section, $G$ is a formal group of finite height over a
local formal scheme $S$ with perfect residue field of characteristic
$p>0$.  
If
\[
    A_{S} \xra{\ell} G
\]
is a level structure, then the Cartier divisor 
\[
    \divisorellA  \eqdef \sum_{a\in A} \divisor{\ell (a)}
\]
is a subgroup scheme of $G$, and the quotient $G/\divisorellA$ is a formal 
group.  A finite free map of formal groups $G \xra{} G'$ is called an
\emph{isogeny}.  In this section we recall the construction of the
isogeny $G \xra{} G/\divisorellA$.

\subsection{The norm}

\label{sec:norm}

An important ingredient in the construction of the quotient is the 
following (see for example \cite[III \S
12]{Mu:AV}, \cite[III \S 2 no. 3]{DeGa:GA},
\cite{Strickland:FiniteSubgps}).  
Let $X \xra{\pi} Y$ be a finite free map of local formal
schemes.  Multiplication by a section $f\in \O_{X}$ defines an
$\O_{Y}$-linear endomorphism $f\cdot$ of $\O_{X}$.  

\begin{Definition}\label{def-norm} 
The \emph{norm}
of $f$ is the determinant 
\[
    \normff_{\pi} f = \Det (f\cdot) \in \O_{Y}.
\]
This is a multiplicative (but not additive) map $\O_{X}\rightarrow
\O_{Y}$.
\end{Definition}

\begin{Definition} \label{def-full-sections}
A list $P_{j}: Y\to X$, $j=1,\dotsc ,k$ of not-necessarily distinct sections
of $\pi$ is called a \emph{full set of sections} if, for all $f\in
\O_{X}$, 
\[
     \normff_{\pi} f = \prod_{j} f (P_{j}),
\]
where we have viewed $P_{j}$ as a $Y$-valued point of $X$ and so
written $f (P_{j})$ for $P_{j}^{*}f$.
\end{Definition}

\begin{Definition}\label{def-norm-line-bundle} 
Let $\L$ be an invertible sheaf of ideals on $X$. 
The 
\emph{norm} of $\L$ is the ideal sheaf 
$\normff_{\pi}\L$ on $Y$ generated by $\normff_{\pi}t$, where $t$ is a
generator of $\L$.  
\end{Definition} 

The multiplicativity of $\normff_{\pi}$ implies that $\normff_{\pi}\L$ is
independent of the choice of generator $t$, and,  if $s$ is a section
of $\L$, then its norm $\normff_{\pi}s$ 
is a section of $\normff_{\pi}\L$.  The norm is not additive, but it is
multiplicative in the sense that  
\[
    \normff_{\pi} (fs) = \normff_{\pi}f\cdot \normff_{\pi}s
\]
if $f$ is a section of $\O_{X}$ and $s$ is a section of $\L$.  

If $\{P_{1},\dotsc ,P_{k} \}$ is a  full set of sections of $\pi$,
then the map 
\[
       \bigotimes_{j} f_{j}\otimes s_{j} \mapsto \prod_{j} f_{j} s_{j} (P_{j})
\]
defines an isomorphism of $\O_{Y}$-modules
\[
        \bigotimes_{j} P_{j}^{*} \L \cong \normff_{\pi}\L.  
\]

In particular, $\normff_{\pi} \L$ is a line bundle, and the 
norm extends to a multiplicative map from the group of invertible
fractional ideals on $X$ to the the group of invertible fractional
ideals on $Y$.

\subsection{Quotients by finite subgroups}

\begin{Definition} \label{def-finite-subgp} 
A \emph{finite subgroup} of $G$ is a divisor $K$ on $G$ which
is also a subgroup scheme.
\end{Definition}

Let us write $\pi$ and $\mu$ respectively for the projection and
multiplication maps
\[
   \pi,\mu: G \times K \to G,
\]
and let $\O_{G/K}$ be the equalizer
\begin{equation} \label{new-eq:28}
\xymatrix{
{\O_{G/K}}
 \ar[rr]
& & 
{\O_{G}}
 \ar@<1ex>[rr]^-{\pi^{*}}
 \ar@<-1ex>[rr]_-{\mu^{*}}
& & 
{\O_{G} \otimes \O_{K}.}
}
\end{equation}

\begin{Proposition}\label{t-pr-quotient-by-finite-subgroup}
\begin{thmList}
\item For $f\in \O_{G}$, we have $\normff_{\pi}\mu^{*}f \in \O_{G/K}$.
\item If $x$ is a coordinate on $G$ and $y = \normff_{\pi}\mu^{*}x$, then 
\[
     \O_{G/K} = \O_{S}\psb{y}.
\]
In particular, $\O_{G/K}$ is the ring of functions on a formal scheme
$G/K$ over $S$.
\item $G/K$ has naturally the structure of a formal group, and as such
is the categorical cokernel of the inclusion
\[
   K \rightarrow G.
\]
\item For any map $f: T\to S$ of local formal schemes, there is a
canonical isomorphism of formal groups
\[
      f^{*}G / f^{*}K \cong f^{*} (G/K).
\]
\end{thmList} 
\end{Proposition}

\begin{proof}
The construction of quotients of formal groups by finite subgroups
goes back to Lubin \cite{Lu:FSFG}.  For the construction of the
quotient in the generality considered here, see
\cite{Drinfeld:EM} and \cite{Strickland:FiniteSubgps}. 
\end{proof}

\begin{Definition}  \label{def-isogeny}  
An \emph{isogeny} is a finite free homomorphism
$q: G\rightarrow G'$
of formal groups, i.e. a homomorphism of formal groups such that $\Ker
q$ is a finite subgroup of $G.$
\end{Definition}

Proposition \ref{t-pr-quotient-by-finite-subgroup} implies that 
isogenies are epimorphisms.

\begin{Lemma} \label{t-isogenies-epi}
If $f: G\to G'$ is an isogeny, then there is a unique isomorphism
$G/\Ker f \cong G'$ making the diagram 
\[
\xymatrix{
{G}   
 \ar[d]
 \ar[dr] \\
{G/\Ker f}
 \ar[r]_-{\cong} 
&
{G'}}
\]
commute.  If $g: G \to G''$ is another isogeny, such that 
\[
\Ker f \subseteq \Ker g,
\]
then there is a unique isogeny
\[
     h: G' \to G''
\]
such that $g=hf.$ \qed
\end{Lemma}

Let $\quot: G\to G'$ be an isogeny, and let $K=\Ker \quot$.  By
Proposition \ref{t-pr-quotient-by-finite-subgroup}, the norm map 
\[
     \normff_{\pi}\mu^{*}: \O_{G} \rightarrow \O_{G}
\]
takes values in $\O_{G/K}\subseteq \O_{G}$.

\begin{Definition}  \label{def-norm-quot}
We write $\normq_{\quot}$ for the norm
\[
    \normq_{\quot}: \O_{G}\to \O_{G'}
\]
induced by $\normff_{\pi} \mu^{*}$ and the isomorphism 
$
    G/K\cong G'.
$
\end{Definition}


\subsection{Level structures, quotients, and the norm of an ideal}
\label{sec:norm-an-ideal-new}

Let $A$ be a finite abelian group, and let 
\[
   A_{S} \xra{\ell} G
\]
be a level structure.

\begin{Proposition} \label{t-pr-level-str-quotient-new}
\begin{thmList}
\item 
The Cartier divisor 
\[
   \divisorellA  = \sum_{a\in A}\divisor{\ell (a)} 
\]
is a subgroup scheme of $G.$
\item
The sections 
\[
    \ell (a): S\to \divisorellA
\]
for $a\in A$ are a full set of sections of $\divisorellA$.
\end{thmList}
\end{Proposition}

\begin{proof}
It suffices to prove the first part in the case that 
$G$ is the universal deformation a formal group 
of height $h$ over a field of characteristic $p$, and $R=\O (\ulvl
(A,G))$.  In that case the result is essentially Proposition 4.4 of
\cite{Drinfeld:EM}: Drinfel'd actually considers a level structure of
the form 
\[
    \Lambda = (\Z/p^{n})^{h} \xra{\ell} G (R)
\]
and a subgroup $A\subseteq \Lambda$, but his argument uses only
the $A$-structure and the fact that $R_{\Lambda} (G)$ is a Noetherian
$p$-regular complete local 
domain.  Strickland \cite{Strickland:FiniteSubgps} gives a complete
proof in the generality considered here.  Katz and Mazur prove the
second part as \cite[Thm. 1.10.1]{KaMa:AMEC}, in the case that $G$ is
smooth curve over a scheme $S$; their proof proceeds by reducing to
and proving the local case considered here.
\end{proof}

It follows from Proposition \ref{t-pr-level-str-quotient-new} that, if 
\[
     G \xra{\quot} G'
\]
is the cokernel of the inclusion $\divisorellA \rightarrow G$,
then the norm of Definition \ref{def-norm-quot} is given by the
formula 
\[
     \quot^{*}\normq_{\quot} f = \prod_{a\in A} T_{a}^{*}f.
\]
In particular, if $x$ is a coordinate on $G$, then 
\begin{equation} \label{eq:lubin-coord}
y = \prod_{a\in
A}T_{a}^{*}x
\end{equation}
is the coordinate on the quotient $G'$ given in Proposition
\ref{t-pr-quotient-by-finite-subgroup}.  Indeed this is the 
coordinate discovered by Lubin \cite{Lu:FSFG}.  A coordinate on $G$ is
a trivialization of the line bundle $\I_{G} (0)$, and it is useful to
interpret the norm in terms of line bundles. 

Let $\L$ be an invertible sheaf of ideals in $\O_{G}$ 
(or, more 
generally, an invertible fractional ideal on $G$).
Proposition \ref{t-pr-quotient-by-finite-subgroup} implies that there
is a line bundle $\normq\L=\normq_{\quot}\L$ on $G'$, 
characterized by the formula 
\[
\quot^{*}\normq\L = \normff_{\pi} \mu^{*}\L;
\]
and if $t$ is a
trivialization of $\L$, then $\normq t$ is a  trivialization of
$\normq\L$. 
The map 
\begin{equation} \label{eq:85}
       \bigotimes_{a\in A} s_{a}\mapsto \prod_{a\in A} s_{a}
\end{equation}
defines an isomorphism 
\begin{equation} \label{eq:83}
    \bigotimes_{a\in A}T_{a}^{*}\L \cong  \quot^{*} \normq\L.
\end{equation}
If $s$ is a section of $\L$, then under this isomorphism 
\[
      \bigotimes_{a\in A} T_{a}^{*}s = \quot^{*}\normq s.
\]
Thus in the presence of the level structure, the first two parts of
Proposition \ref{t-pr-quotient-by-finite-subgroup} take the following
form.

\begin{Proposition} \label{t-pr-OGp-NOG-new}
There are canonical natural isomorphisms
\begin{align}
   \normq\O_{G} & \cong \O_{G'} \notag \\
   \normq\I_{G} (0) & \cong \I_{G'} (0)  \notag \\
  \quot^{*}\I_{G'} (0) & \cong \bigotimes_{a\in A} T_{a}^{*}I_{G} (0)
\cong I_{G} (\ell).   \label{eq:84} 
\intertext{If $s$ is a coordinate on $G$, then $\normq s$ is a
coordinate on $G'$, and under the isomorphism \eqref{eq:84},}
    \quot^{*}\normq s & = \bigotimes_{a\in A} T_{a}^{*}s. \notag
\end{align} \qed
\end{Proposition}

\section{Descent for level structures}

\label{sec:desc-level-struct}

In Definition \ref{def-1} we described ``descent data for level
structures'' as they appear on the formal group of an $\hinfty$ ring
spectrum.  In this section, we give an equivalent description
(see Proposition \ref{t-pr-descent-is-descent}) which displays the
relationship to the usual notion of descent data.  In addition to
justifying the terminology, the new formulation simplifies the task of
showing that the Lubin-Tate formal groups have canonical descent data
for level structures (Proposition \ref{t-pr-lubin-tate-hinfty}).

\subsection{Composition of isogenies: the simplicial functor $\ulvl_{*}$}

Let $\Fgps$ be the functor from admissible local rings
$R$ to sets whose value on $R$ is the set of formal groups $G/\spf
R$.  If $f: R\to R'$ is a map of admissible local rings, then $\Fgps
(f)$ sends $G/\spf R$ to $f^{*}G/\spf R'$.

Let 
\[
\ulvl (A) \xra{}\Fgps
\]
be the functor \emph{over} $\Fgps$ whose value on $R$ is the
set of formal groups $G/\spf R$ equipped with a level structure 
\[
   A_{\spf R} \xra{} G.
\]
We define 
\[
   \ulvl_{1} \eqdef \coprod_{A_{0}} \ulvl (A_{0});
\]
the coproduct is over all finite abelian groups.  We have adorned the
$\ulvl$ and the $A$ with subscripts so that we 
can make the more general definition 
\[
   \ulvl_{n}  = \coprod_{0 = A_{n}\subseteq A_{n-1}\dotsb \subseteq
A_{0}} \ulvl (A_{0}).
\]
The coproduct is over all sequences of \emph{inclusions of finite
abelian groups with $A_{n} = 0$}.  With this convention we also
have 
\[
      \Fgps  = \ulvl_{0}.
\]
We write 
\begin{equation} \label{eq:57}
  d_{0} : \ulvl_{1}\rightarrow \Fgps  
\end{equation}
for the structural map.  

Over $\ulvl_{1} (A)$ we have a level structure 
\[
     A \xra{\ell} d_{0}^{*} G 
\]
and an isogeny 
\[
       d_{0}^{*} G \xra{\quot_{A}} G/\ell (A)
\]
with kernel $A$.  These assemble to give a group $G/\ell$ and an 
isogeny 
\[
       d_{0}^{*} G \xra{\quot} G/\ell
\]
over $\ulvl_{1}$.  We write 
\begin{equation} \label{eq:58}
   d_{1}: \ulvl_{1} \xra{}\Fgps 
\end{equation}
for the map classifying $G/\ell$.

\begin{Lemma} \label{t-le-level-restricts-passes-to-quot}
Let 
\[
   A\xra{\ell} G
\] 
be a level structure.  If $B\subseteq A$, then the induced map 
\[
  B \xra{\ell\restr{B}} G
\]
is a level structure.  If $\quot: G \xra{} G'$ is an isogeny with kernel
$\ell\restr{B}$, then the induced map 
\[
    \ell': A/B \xra{} G'
\]
is a level structure.
\end{Lemma}

\begin{proof}
The first part is clear from the definition of a level
structure~\eqref{def-level-structure}.  For the second part, consider
the diagram 
\[
\begin{CD}
A @> \ell >> G \\
@VVV         @VV \quot V \\
A/B @> \ell' >> G'.
\end{CD}
\]
Let $D$ be the divisor 
\[
    D = \sum_{\substack{a\in A \\ p a = 0}}\divisor{\ell (a)}
\]
on $G$; by hypothesis we have an inequality of Cartier divisors 
\[
    D \leq G[p].
\]
It follows that 
\[
     \sum_{b\in B} T_{b}^{*}D \leq \sum_{b\in B} T_{b}^{*} G[p].
\]
The formula~\eqref{eq:lubin-coord} for the coordinate on the quotient $G'$
shows that the left side descends to the divisor 
\[
\sum_{\substack{c\in ( A/ B) \\ p c = 0}}\divisor{\ell' (c)},
\]
while the right side descends to the divisor $G'[p]$.
\end{proof}

The Lemma gives maps
\[
   d_{j}:  \ulvl_{n}  \rightarrow \ulvl_{n-1}
\]
for $0\leq j\leq n$ as follows.  For $0\leq j\leq n-1$, the map
$d_{j}$ sends a point
\begin{equation} \label{eq:30}
    0 = A_{n}  \subseteq \dotsb \subseteq  A_{j} \subseteq \ldots \subseteq A_{0} \xra{} G
\end{equation}
of $\ulvl_{n}$
to the point
\[
    0 = A_{n}  \subseteq\dotsb   \widehat{A_{j}} 
    \dotsb  \subseteq A_{0} \xra{} G
\]
of $\ulvl_{n-1}$ obtained by omitting $A_{j}$.  The map $d_{n}$
sends~\eqref{eq:30} 
to 
\[
      0 = A_{n-1}/ A_{n-1} \subseteq \dotsb  \subseteq A_{0}/A_{n-1}
\xra{}  G/\ell (A_{n-1}).
\]
In the case $n=1$ these are just the maps~\eqref{eq:57} and~\eqref{eq:58}.
We also have for $0\leq j \leq n$ a map 
\[
    s_{j}: \ulvl_{n}  \rightarrow \ulvl_{n+1}
\]
which sends the sequence~\eqref{eq:30}
to the sequence
\[
A_{n}\subseteq \dotsb \subseteq A_{j} \subseteq A_{j} 
\subseteq \ldots \subseteq A_{0} \xra{}G
\]
obtained by repeating $A_{j}$.  It is easy to check that

\begin{Lemma}\label{t-le-ulvl-simplicial}
 $(\ulvl_{*},d_{*},s_{*})$ is a simplicial functor. \qed
\end{Lemma}

\subsection{Descent data for functors over formal groups}

Now suppose that 
\[
\Prob \xra{}\Fgps
\]
is a functor over $\Fgps$, and if $x\in \Prob (R)$ is an
$R$-valued point, let's write $G_{x}$ for the resulting formal group
over $\spf R$.  As in the previous section, we define 
\begin{align*}
   \ulvl (A,\Prob) & = \ulvl (A)\times_{\Fgps} \Prob \\
   \ulvl_{n} (\Prob) & = \ulvl_{n}\times_{\Fgps} \Prob
\end{align*}
and so on.  A point $(\ell,x)\in\ulvl (A,\Prob) (R)$ is a point $x$ of
$\Prob (R)$ and a level structure 
\[
    A \xra{\ell} G_{x}.
\]
We write
\begin{equation} \label{eq:51}
    d_{0}: \ulvl_{1} (\Prob) \rightarrow \ulvl_{0} (\Prob) = \Prob.
\end{equation}
for the forgetful map 
\[
    d_{0} (\ell,x) = x.
\]
We also always have degeneracies 
\[
       s_{j}: \ulvl_{n} (\Prob) \rightarrow \ulvl_{n+1} (\Prob)
\]
for $0\leq j\leq n$.  

If $(\ell, x)$ is an $R$-valued point of $\ulvl (A,\Prob)$, then we
get an isogeny 
\[
       G_{x} \xra{} G_{x} / \ell.
\]
Suppose that we have a natural transformation 
\begin{equation} \label{eq:32}
       d_{1}: \ulvl_{1} (\Prob) \xra{} \Prob
\end{equation}
such that 
\begin{equation} \label{eq:33}
       G_{d_{1} (\ell,x)} = G_{x} / \ell,
\end{equation}
or equivalently that the diagram 
\begin{equation} \label{eq:52}
\begin{CD}
  \ulvl_{1} (\Prob) @>>> \ulvl_{1}  \\
   @V d_{1} VV @VV d_{1} V \\
  \Prob  @>>> \Fgps 
\end{CD}
\end{equation}
commutes.  Lemma \ref{t-le-level-restricts-passes-to-quot} then 
gives maps
\[
       d_{j}: \ulvl_{n} (\Prob) \rightarrow \ulvl_{n-1} (\Prob)
\]
for $0\leq j\leq n$.  

\begin{Definition}  \label{def-descent-functors}
\emph{Descent data for level structures on the
functor} $\Prob$ consist of a natural transformation~\eqref{eq:32}
such that 
\begin{enumerate}
\item the diagram~\eqref{eq:52} commutes, and 
\item $(\ulvl_{*} (\Prob),d_{*},s_{*})$ is a simplicial functor.
\end{enumerate}
\end{Definition}

\begin{Remark}
It is equivalent to ask for natural transformations
\[
       d_{j}: \ulvl_{n} (\Prob) \rightarrow \ulvl_{n-1} (\Prob)
\]
for $n\geq 1$ and $0\leq j\leq n$, such that $(\ulvl_{*}
(\Prob),d_{*},s_{*})$ is a simplicial functor, and the levelwise
natural transformation  
\[
    \ulvl_{*} (\Prob) \xra{} \ulvl_{*}
\]
is a map of simplicial functors.
\end{Remark}

For example, let $G$ be a formal group of finite height over a
$p$-local formal scheme $S$.  The formal scheme $S$ has the structure
of a functor over $\Fgps$: if $x: \spf R\to S$ is a point of $S$, then 
\[
    G_{x} = x^{*}G;
\]
We briefly write $G/S$ for $S$, considered as a functor over $\Fgps$ in this
way.  The functor $\ulvl (A,G/S)$ is just the functor called $\ulvl
(A,G)$ in \S~\ref{sec:level-structures}; in particular it is
represented by the $S$-scheme $\ulvl (A,G)$ of Lemma
\ref{t-le-ulvl-represented}. To give maps $\psi_{\ell}$ and $f_{\ell}$ 
which satisfy condition \eqref{descent-it-natural} of 
Definition \ref{def-1} amounts
to giving a map 
\[
   d_{1}: \ulvl_{1} (G/S) \rightarrow S
\]
and an isogeny 
\[
     d_{0}^{*} G \xra{\quot} d_{1}^{*}G
\]
whose kernel on $\ulvl (A,G/S)$ is $A$.  Lemma
\ref{t-le-level-restricts-passes-to-quot} gives maps 
\[
   d_{j}: \ulvl_{n} (G/S) \xra{} \ulvl_{n-1} (G/S)
\]
for $0\leq j\leq n$ as explained above.  With these definitions, 
parts \eqref{descent-it-2} and
\eqref{descent-it-3} of Definition \ref{def-1}
are equivalent to asserting that
\[
(\ulvl_{*}(G/S),d_{*},s_{*})
\]
is a simplicial 
functor, and over $\ulvl_{2} (G)$ the diagram 
\begin{equation} \label{eq:35}
\xymatrix{
&
{d_{0}^{*}d_{0}^{*} G}
\ar@{=}[dl] 
\ar[dr]^{d_{0}^{*}\quot}
\\
{d_{1}^{*} d_{0}^{*} G}
 \ar[d]_{d_{1}^{*}\quot}
& 
&
{d_{0}^{*} d_{1}^{*} G }
 \ar@{=}[d]
\\
{d_{1}^{*}d_{1}^{*}G}
 \ar@{=}[dr]
& 
&
{d_{2}^{*} d_{0}^{*} G}
 \ar[dl]^{d_{2}^{*}\quot}
\\
&
{d_{2}^{*} d_{1}^{*}G}
}
\end{equation}
commutes.

A more convenient formulation of Definition~\ref{def-1} is the
following. Let $\GS$ to be the functor over $\Fgps$
whose value on $R$  is the set of  
pull-back diagrams 
\[
\begin{CD}
        G' @> f >> G \\  
       @VVV    @VVV \\   
        \spf R @> i >> S.
\end{CD}
\]
such that the map 
\[
    G' \to i^{*} G
\]
induced by $f$ is a homomorphism (hence isomorphism) of formal groups
over $\spf R$.  
For a finite abelian group $A$, 
$\ulvl (A,\GS) (R)$ is the set of diagrams 
\[
\xymatrix{
{A_{\spf R}} 
 \ar[r]^-{\ell}
 \ar[dr]
&
{G'}
 \ar[r]^{f}
 \ar[d]
&
{G}
 \ar[d]
\\
&
{\spf R}
 \ar[r]^{i}
&
{S,}
}
\]
where the square part is a point of $\GS (R)$ and $\ell$ is a level
structure.  To give a map of functors 
\begin{equation} \label{eq:11}
     \ulvl_{1} (\GS) \xra{d_{1}} \GS
\end{equation}
making the diagram 
\[
\begin{CD}
   \ulvl_{1} (\GS) @> d_{1} >> \GS \\
   @VVV                        @VVV \\
   \ulvl_{1}  @> d_{1} >> \Fgps 
\end{CD}
\]
commute is to
give a pull-back diagram 
\[
\begin{CD}
    G/\ell @>>> G \\
    @VVV       @VVV \\
     \ulvl_{1}(G/S) @>>> S;
\end{CD}
\]
it is equivalent to give a map of formal schemes
\[
   d_{1}: \ulvl_{1} (G/S) \rightarrow S
\]
and an isogeny 
\[
     d_{0}^{*} G \xra{\quot} d_{1}^{*}G
\]
whose kernel on $\ulvl (A,G/S)$ is $A$. 

\begin{Proposition} \label{t-pr-descent-is-descent}
Let $G$ be a formal group over an admissible local ring $R$, and let
$S=\spf R$.  Descent data for level structures on 
the group $G/S$ are equivalent to descent data for level structures on 
the functor $\GS$. 
\end{Proposition}

\begin{proof}
One checks that the commutativity of the diagram~\eqref{eq:35} has
been incorporated in the structure of the functor $\GS$.
\end{proof}

\subsection{Noetherian rings and Artin rings} \label{sec:noeth-rings-artin}

Suppose that $\cat{D}$ is a subcategory of the category of admissible
local rings.  If $\Prob$
is a functor from complete local rings to sets, let
$\Prob^{\cat{D}}$ denote its restriction
to $\cat{D}$.

\begin{Definition} \label{def-desc-artin} 
\emph{Descent data for level structures on $\Prob^{\cat{D}}$}
consists of a natural transformation
\[
    d_{1}: \ulvl_{1} (\Prob)^{\cat{D}} \xra{d_{1}}
    \Prob^{\cat{D}},
\]
such that the restriction to $\cat{D}$ of the diagram 
\eqref{eq:52} commutes, and such that the $(\ulvl_{*} (\Prob))^{\cat{D}},
d_{*},s_{*})$ is a simplicial functor.
\end{Definition}

For example, let $\noether$ be the category of Noetherian complete
local rings, and let $\artin$ be the category of Artin local rings.
If $S$ and $T$ are Noetherian local formal
schemes, then the natural maps
\splitpageyesno{
\begin{multline} \label{eq:63}
     \CatOf{formal schemes} (S,T) \rightarrow 
     \CatOf{functors} (S^{\noether}, T^{\noether})  \\ 
     \rightarrow 
     \CatOf{functors} (S^{\artin},T^{\artin})
\end{multline}
}{
\begin{equation} \label{eq:63}
     \CatOf{formal schemes} (S,T) \rightarrow 
     \CatOf{functors} (S^{\noether}, T^{\noether}) \rightarrow 
     \CatOf{functors} (S^{\artin},T^{\artin})
\end{equation}
}
are isomorphisms.

\begin{Proposition} \label{t-pr-restr-to-Artin-OK}
If $G$ is a formal group over a
Noetherian
local formal scheme $S$ with perfect residue field of characteristic
$p>0$, then the forgetful maps, from the set of descent data for level
structures on $\GS$ to 
the set of descent data for level structures on $\GS^{\noether}$ and
on $\GS^{\artin}$, are isomorphisms. 
\end{Proposition}

\begin{proof}
If $G$ is a formal group over a Noetherian local formal scheme, then
by Proposition \ref{t-pr-level-A-finite-flat}, 
$\ulvl (A,G)$ is also a Noetherian local formal scheme.  The result
follows easily from the isomorphism \eqref{eq:63}.
\end{proof}

\section{Lubin-Tate groups}
\label{sec:lubin-tate-groups}

Let $k$ be a perfect field of characteristic $p>0$, and let $\Gamma$ be
a formal group of finite height over $k$.  In this section we shall
prove that the universal deformation of $\Gamma$ has descent for level
structures. 

\subsection{Frobenius}

Let $k$ be a perfect field of characteristic $p>0$, and let $\Gamma$ be
a formal group of finite height over $k$.  The Frobenius map $\frob$
gives rise to a relative Frobenius $F$ map as in the diagram 
\[
\xymatrix{
{\Gamma} 
 \ar[r]^{F}
 \ar[dr]
 \ar@/^2pc/[rr]^{\frob_{\Gamma}}
&
{\frob_{k}^{*}\Gamma}
 \ar[r]
 \ar[d]
&
{\Gamma}
 \ar[d]
\\
&
{\spec k}
 \ar[r]^{\frob_{k}}
&
{\spec k.}
}
\]
The Frobenius map $F$ is an isogeny of degree $p$, with kernel the
divisor $p\divisor{\e}$.

\subsection{Deformations}

If $T$ is a local formal scheme, then we write $T_{0}$ for its closed point.

\begin{Definition}  \label{def-deformation}
Let $T$ be a local formal scheme.  
A \emph{deformation of $\Gamma$ to $T$}
is a triple $(H/T,f,j)$ consisting of
a formal group $H$ over $T$ and a pull-back  
diagram 
\[
\begin{CD}
      H_{\specialpoint{T}} @> f    >> \Gamma \\
      @VVV              @VVV \\
      \specialpoint{T} @> j    >> \spec k,
\end{CD}
\]
such that the induced map $H_{\specialpoint{T}} \rightarrow j^{*} \Gamma$ is a
homomorphism (and so isomorphism) of formal groups over $\specialpoint{T}$.  
The functor from complete local rings to sets which assigns to $R$ the
set of deformations of $\Gamma$ to $\spf R$ will be denoted 
$\Deformations{\Gamma}$.
\end{Definition}

From the definition it is clear that if $(H/T,f,j)$ is a deformation
of $\Gamma$, then there is a natural transformation 
\[
   \HT  \rightarrow \Deformations{\Gamma}.
\]
Lubin and Tate \cite{LubinTate:FormalModuli} construct a 
deformation $(G/S,\funiv,\juniv)$ with an isomorphism 
\begin{equation} \label{eq:60}
   S \iso \spf \W k \psb{u_{1},\dots,u_{h-1}}
\end{equation}
inducing $\juniv: S_{0}\iso \spec k$
such that the natural transformation 
\begin{equation} \label{eq:74}
   \GS \rightarrow \Deformations{\Gamma}
\end{equation}
is an isomorphism of functors over $\Fgps$.  

\begin{Remark} 
Lubin and Tate claim only that the restriction to Noetherian complete
local rings  
\[
    \GS^{\noether} \rightarrow  \Deformations{\Gamma}^{\noether}
\]
is an isomorphism.  In fact, their argument proves the stronger
statement.  The main point is that, if $\Phi$ is a formal group law over
a field $k$ of characteristic $p>0$ and if $M$ is \emph{any}
$k$-vector space (not necessarily finite-dimensional), then the
natural map 
\[
    H^{2}_{k} (\Phi) \otimes M \rightarrow H^{2}_{M} (\Phi)
\]
is an isomorphism, a result which Lubin and Tate assert at the
beginning of 2.3 for $M$ finite dimensional.  Indeed, the 
argument of their Proposition 2.6 may be applied to give this
calculation of  $H^{2}_{M} (\Phi)$.
\end{Remark}

\subsection{Descent for level structures on deformations}

We continue to fix a formal group $\Gamma$ of finite
height over a perfect field $k$ of characteristic $p>0$.

Let $A$ be a finite abelian group.  If $R$ is a complete local ring,
then a point of $\ulvl (A,\Deformations{\Gamma})$ is a commutative diagram 
\begin{equation} \label{eq:31}
\begin{CD}
A @> \ell >> H @<<< H_{T_{0}} @> f >> \Gamma \\
@.          @VVV      @VVV             @VVV \\
@.           T  @<<< T_{0} @> j >> \spec k,
\end{CD}
\end{equation}
consisting of a deformation $(H,f,j)$ of $\Gamma$ to $R$ and a level
structure $\ell$ on $H$.  
The level structure in~\eqref{eq:31} of $\Deformations{A,\Gamma} (R)$
has a cokernel 
\[
     H \xra{\quot} H'
\]
If $x$ is 
a coordinate on $H$, then $x (\ell (a))$ is topologically nilpotent in
$\O (T)$.  It follows that $x (\ell (a)) = 0$ on $\specialpoint{T}$,
and so there is a canonical isomorphism $\overline{\quot}$
making the diagram 
\[
\xymatrix{
{H_{\specialpoint{T}}}
 \ar[d]_-{\quot_{\specialpoint{T}}}
 \ar[r]^-{F^{r}}
&
{(\frob^{r})^{*} H_{\specialpoint{T}}}
\ar[d] 
\ar[r]^-{\text{can}}
&
{H_{\specialpoint{T}}}
 \ar[r]^-{f}
 \ar[d]
&
{\Gamma}
 \ar[d] \\
{H'_{{\specialpoint{T}}}}
 \ar[ur]_-{\overline{\quot}}
 \ar[r]
&
{\specialpoint{T}}
 \ar[r]^-{\frob^{r}}
&
{\specialpoint{T}}
 \ar[r]^-{j}
&
{\spec k}
}
\]
commute.   In other words $(H',f \text{can} \overline{\quot},j
\frob^{r})$ is a point of $\Deformations{\Gamma} (T)$, and we have
constructed a natural transformation 
\begin{equation} \label{eq:50}
    \ulvl_{1} (\Deformations{\Gamma})  \xra{d_{1}} \Deformations{\Gamma}.
\end{equation}
satisfying \eqref{eq:33}.  The fact that $\frob^{r}\frob^{s} =
\frob^{r+s}$ then implies

\begin{Lemma}
The map $d_{1}$ is descent data for level structures on the functor
$\Deformations{\Gamma}$. \qed
\end{Lemma}

Now let $(G/S,\funiv,\juniv)$ be Lubin and Tate's universal deformation of
$\Gamma/\spec k$.  Using the isomorphism~\eqref{eq:74}
and Proposition \ref{t-pr-restr-to-Artin-OK} 
we may trade $\GS$ for $\Deformations{\Gamma}$ in \eqref{eq:50} to get
a map 
\begin{equation} \label{eq:34}
    d_{1} : \ulvl_{1} (\GS) \xra{} \GS
\end{equation}
such that 

\begin{Proposition}\label{t-pr-lubin-tate-hinfty} 
The natural transformation $d_{1}$ is descent data for level
structures on the functor $\GS$, and so gives descent data for level
structures on the formal group $G/S$. \qed
\end{Proposition}

\subsection{Comparison to the descent data coming from the
$E_{\infty}$ structure of Goerss and Hopkins}

The construction of the descent data in Proposition
\ref{t-pr-lubin-tate-hinfty} uses 
the equality 
\begin{equation} \label{eq:40}
 p^{r} \divisor{\e} =\Ker F^{r} : \Gamma
\rightarrow 
(\frob^{r})^{*}\Gamma
\end{equation}
of divisors on $\Gamma$  and
the equation 
\begin{equation} \label{eq:59}
    F^{r+s} = ( (\frob^{r})^{*} F^{s} ) F^{r}: 
   \Gamma  \rightarrow 
   (\frob^{r+s})^{*} \Gamma
\end{equation}
More generally, to give descent data for level structures on $\GS$ is
equivalent to giving a collection of isogenies 
\[
      F_{r}: \Gamma
             \rightarrow  (\frob^{r})^{*}\Gamma
\]
for $r\geq 1$ satisfying the analogues of~\eqref{eq:40}
and~\eqref{eq:59}.  The descent data in the Proposition are uniquely
determined by the choice $F_{r} = F^{r}$.  

Now let $\E$ be the homogeneous ring spectrum such that $G=\GpOf{\E}$ is
Lubin and Tate's universal deformation of $\Gamma$, so 
\[
    S_{\E} = S = \spf  \W k \psb{u_{1},\dots ,u_{h-1}}.
\]
In work in preparation, Goerss and Hopkins \cite{gh:rcrs} have shown
that $\E$ is an $\einfty$ ring spectrum; by Theorem
\ref{t-th-hinfty-adds} it follows that there is a map 
\[
        \ulvl_{1} (\GS) \xra{d_{1}^{\hopgoe}} \GS
\]
giving descent data for level structures on $\GS$. 

Let $A$ be a finite group of order $p^{r}$, and let
\[
     A_{\spf R} \xra{\ell} i^{*} G
\]
be a level structure on $G$.  
Reducing modulo the maximal ideal in the construction of $\psile$
\eqref{def-psile}, one sees that 
\[
     \psile = \frob^{r} : S_{\E} \rightarrow S_{\E}.
\]
Examination of the construction~\eqref{eq:61} of $\psilge{\E}$ shows
that 
\[
    (\psilge{\E})_{\specialpoint{S}} = F^{r} : 
    G_{\specialpoint{S}} \rightarrow 
    (\frob^{r})^{*}G_{\specialpoint{S}}.
\]
Thus we have the following result.

\begin{Proposition}\label{t-pr-d-1-top-is-d-1-frob}
If $\E$ is the spectrum associated to the universal
deformation of a formal group $\Gamma$ of finite height over a perfect
field $k$, then the descent data for level structures on $\GpOf{\E}$
provided by the $\einfty$ structure of Goerss and Hopkins coincide with
the descent data in Proposition \ref{t-pr-lubin-tate-hinfty}. \qed
\end{Proposition}

\begin{Remark} \label{rem-5} 
At the time of the writing of this paper, the result of Goerss and 
Hopkins is not published.  The arguments of this section do
not depend on their result beyond the existence of the $\hinfty$
structure, so a cautious statement of the Proposition is that the
descent data for 
level structures on $\GpOf{\E}$ provided by \emph{any} $\hinfty$ structure on
$\E$ coincide with the descent data in Proposition
\ref{t-pr-lubin-tate-hinfty}.
\end{Remark}

\part{The sigma orientation}
\label{part:sigma-orientation}
\section{$\Theta^{k}$-structures}
\label{sec:thetak-structures}

\subsection{The functors $\Theta^{k}$}

Suppose that $G$ is a formal group 
over a formal scheme $S$, and suppose that 
$\L$ is a line bundle over $G$.

\begin{Definition}
 A \emph{rigid} line bundle over $G$ is a line bundle $\L$ equipped
 with a specified trivialization of $\e^{*}\L$.
 A \emph{rigid section} of such a line bundle is a section $s$ which
 extends the specified section at the identity.  A \emph{rigid
   isomorphism} between two rigid line bundles is an isomorphism which
 preserves the specified trivializations.
\end{Definition}

\begin{Definition} \label{defn-thetat}
 Suppose that $k\geq 1$.   We
 define the line bundle $\Theta^k(\L)$ over $G^k$ by the formula
 \begin{equation}\label{eq-Definition-of-Theta-L}
    \Theta^k(\L) \eqdef
     \bigotimes_{I\subseteq \{1,\dots,k\}}
      (\mu_{I}^{*}\L)^{(-1)^{|I|}}.
 \end{equation}
If $s$ is a section of $\L$, then we
write $\Theta^{k}s$ for the section 
\[
         \Theta^{k} s = \bigotimes_{I\subseteq \{1,\dots,k\}}
      (\mu_{I}^{*}s)^{(-1)^{|I|}}.
\]
of $\Theta^{k} (\L)$.  
We define $\Theta^0(\L)=\L$ and $\Theta^{0} s = s$.
\end{Definition}

For example we have 
\begin{align*}
 \Theta^{1} (\L) & = \frac{\struc^{*}\e^{*}\L}{\L} \\
 \Theta^1(\L)_{a} & = \frac{\L_0}{\L_a} \\
 \Theta^2(\L)_{a,b} & = \frac{\L_0\ot\L_{a+b}}{\L_a\ot\L_b} \\
 \Theta^3(\L)_{a,b,c} & = 
 \frac{\L_0 \ot \L_{a+b} \ot \L_{a+c} \ot \L_{b+c}}
       {\L_a \ot \L_b \ot \L_c \ot \L_{a+b+c}}.
\end{align*}

We observe three facts about these bundles.
\begin{enumerate}
 \item $\Theta^k(\L)$ has a natural rigid structure for $k>0$.
 \item For each permutation $\sigma\in\Sigma_k$, there is a canonical
  isomorphism
  \[ \xi_\sigma:\pi_\sigma^*\Theta^k(\L)\iso\Theta^k(\L). \]
  Moreover,
  these isomorphisms compose in the obvious way.
 \item There is a canonical identification (of rigid line bundles over
 $X^{k+1}$)
\splitpageyesno{
\begin{multline}
\label{eq-Theta-cocycle-iso}
  \Theta^k(\L)_{a_1,a_2,\ldots} \ot
  \Theta^k(\L)_{a_0+a_1,a_2,\ldots}^{-1} \ot \\
  \Theta^k(\L)_{a_0,a_1+a_2,\ldots} \ot
  \Theta^k(\L)_{a_0,a_1,\ldots}^{-1} \iso 1.
\end{multline}
}
{
\begin{equation}
\label{eq-Theta-cocycle-iso}
  \Theta^k(\L)_{a_1,a_2,\ldots} \ot
  \Theta^k(\L)_{a_0+a_1,a_2,\ldots}^{-1} \ot
  \Theta^k(\L)_{a_0,a_1+a_2,\ldots} \ot
  \Theta^k(\L)_{a_0,a_1,\ldots}^{-1} \iso 1.
\end{equation}
}
\end{enumerate}

\begin{Definition} \label{defn-cubical-structure}
 A $\Theta^k$--structure on a line bundle $\L$ over $G$ is a
 trivialization $s$ of the line bundle $\Theta^k(\L)$ such that
 \begin{enumerate}
  \item for $k>0$, $s$ is a rigid section; \label{item-rigid}
  \item $s$ is symmetric in the sense that for each $\sigma\in\Sigma_k$, we
   have $\xi_\sigma\pi_\sigma^*s=s$;
  \item we have
\splitpageyesno{
\begin{multline}
\label{was-eq:8}
s(a_1,a_2,\ldots) \ot s(a_0+a_1,a_2,\ldots)^{-1} \ot \\
    s(a_0,a_1+a_2,\ldots) \ot s(a_0,a_1,\ldots)^{-1} = 1
\end{multline} 
}
{
\begin{equation} \label{was-eq:8}
s(a_1,a_2,\ldots) \ot s(a_0+a_1,a_2,\ldots)^{-1} \ot
    s(a_0,a_1+a_2,\ldots) \ot s(a_0,a_1,\ldots)^{-1} = 1
\end{equation}
}
   under the isomorphism~\eqref{eq-Theta-cocycle-iso}.
 \end{enumerate}
 A $\Theta^3$--structure is known as a \emph{cubical structure}
 \cite{Br:FTTC}.  We write $C^k(G;\L)$ for the set of
 $\Theta^k$-structures on $\L$ over $G$.  Note that $C^0(G;\L)$ is
 just the set of trivializations of $\L$, and $C^1(G;\L)$ is the set
 of rigid trivializations of $\Theta^1(\L)$.  We also define a functor
 from rings to sets by
\splitpageyesno{
 \begin{align*} \uC^k(G;\L)(R) = 
     \{(u,f) \suchthat & u:\spec(R)\xra{}S\;,\\
                       & f\in C^k_{\spec(R)}(u^*G;u^*\L)\}, 
 \end{align*}
}
{
 \[ \uC^k(G;\L)(R) = 
     \{(u,f) \suchthat u:\spec(R)\xra{}S\;,\;
        f\in C^k_{\spec(R)}(u^*G;u^*\L)\}, 
 \]
}
and we recall the following.
\end{Definition}

\begin{Proposition}[\cite{AHS:ESWGTC}]
Let $G$ be a formal group over a scheme $S$, and let $\L$ be a line
bundle over $G$.  The functor $\uC^{k} (G;\L)$ is represented by an
affine scheme over $S$, and for $j: S'\to S$, the natural map 
\[
       \uC^{k} (j^{*}G;j^{*}\L) \rightarrow j^{*}\uC^{k} (G;\L)
\]
is an isomorphism. \qed
\end{Proposition}

\subsection{Relations among the $\Theta^{k}$: the functor $\Delta$}

\begin{Definition} \label{def-Delta} 
If $\M$ is a line bundle over $G^{n}$, then we define $\Delta \M$ to
be the rigid line bundle over $G^{n+1}$ given fiberwise by the formula 
\[
   \Delta \M_{a_{1},a_{2},\dots ,a_{n+1}} = 
   \frac{\M_{a_{1},a_{3},\dots ,a_{n+1}}\ot \M_{a_{2},\dots ,a_{n+1}}}
        {\M_{a_{1}+a_{2},a_{3},\dots,a_{n+1}}\ot\M_{0,a_{3},\dots ,a_{n+1}}}.
\]
If $s$ is a section of $\M$ then we write $\Delta s$ for the rigid section 
\[
   \Delta s (a_{1},\dots ,a_{n+1}) =
           \frac{s (a_{1},\dots )\ot s(a_{2,\dots })}
                {s (a_{1} + a_{2},\dots )\ot s (0,a_{3},\dots)}
\]
of $\Delta \M$.
\end{Definition}

The following can be checked directly from the definitions.

\begin{Lemma} \label{t-Delta-facts}
\begin{thmList}
\item $\Delta$ is multiplicative: if $\M$ is a line bundle over
$G^{n}$ then there is a canonical isomorphism of
rigid line bundles  
\begin{equation} \label{eq-delta-mult}
     \Delta (\M_{1}\ot \M_{2}) \iso \Delta (\M_{1})\ot \Delta (\M_{2}).
\end{equation}
\item Under the identification \eqref{eq-delta-mult}, one has 
\[
     \Delta ( s_{1}\ot s_{2} ) = \Delta (s_{1}) \ot \Delta (s_{2}).
\]
\item If $\L$ is a line bundle over $G$ then for $k\geq 2$ there is a
canonical isomorphism of rigid line 
bundles 
\begin{equation}\label{was-eq:1}
     \Theta^{k}\L \iso \Delta \Theta^{k-1} \L.
\end{equation}
\item If $s$ is a section of $\L$ then under the isomorphism
\eqref{was-eq:1}, one has  
\[
     \Theta^{k}s = \Delta \Theta^{k-1}s
\]
\end{thmList} \qed
\end{Lemma}

\section{The norm map for
$\Theta^{k}$-structures}
\label{sec:norm-map-thetak-new}

Suppose that $A$ is a finite group, and let 
\[
    A_{S}\xra{\ell} G
\]
be a level structure on a formal group $G$ of finite height over a
local formal scheme with perfect residue field of characteristic
$p>0$.  Let
\[
    K = \divisorellA,
\]
and let 
\[
   G \xra{\quot} G'
\]
be the quotient of $G$ by $K$.
	
Fix $k\geq 1$, and, for $1\leq i\leq k$, let 
\[
   \mu_{i}, w_{i}: G^{k}\times K  \rightarrow G^{k} 
\]
be the maps given in punctual notation by 
\begin{align*}
 \mu_{i} (g_{1},\dotsc ,g_{k},a) = (g_{1},\dotsc ,g_{i}+a,\dotsc,g_{k}) \\
  w_{i} (g_{1},\dotsc ,g_{k},a) = (g_{1},\dotsc ,a,\dotsc,g_{k});
\end{align*}
that is, $w_{i}$ replaces $g_{i}$ with $a$.  
Let
\[
    \pi: G^{k}\times K \rightarrow G^{k}
\]
the projection onto the first $k$ factors.  

Let 
\[
      A_{S} \xra{\ell_{i}} G^{k}
\]
be the homomorphism to the $i$ factor, let 
\[
    G_{i} = G^{i-1}\times G' \times G^{k-i}, 
\]
and let 
\[
     \quot_{i}: G^{k}\to G_{i}
\]
be the projection.

If $\M$ is  an invertible fractional ideal on $G^{k}$, then by
Proposition \ref{t-pr-quotient-by-finite-subgroup}, we may define
$\normq_{i}\M$ to be the invertible fractional ideal on $G_{i}$ such that 
\[
  \quot_{i}^{*}\normq_{i}\M  = \normff_{\pi} (\mu_{i}^{*}\M)  \otimes \normff_{\pi} (
w_{i}^{*}\M)^{-1}.
\]
For $a\in A$, let $\tilde T_a \M$ be the line bundle whose fiber over
$(g_{1},\dots ,g_{k})$ is
\begin{equation} \label{eq:53-new}
\rT{a}\M_{g_{1},\dotsc ,g_{k}} \eqdef \frac{\M_{(a+g_{1},g_{2},\dots ,g_{k})}}{\M_{(a,g_{2},\dots)}}.
\end{equation}
If $s$ is a section of $\M$, define $\rT{a}s$ by
\[
\rT{a} s(g_{1},\dots ,g_{k}) = 
\frac{s(a+g_{1},\dots ,g_{k})}
     {s(a,g_{2},\dots ,g_{k})}.
\]
Proposition \ref{t-pr-level-str-quotient-new} implies that there is a
canonical isomorphism 
\begin{equation} \label{eq:86-new}
   \quot_{1}^{*} \normq_{1}\M \cong \bigotimes_{a\in A}\rT{a}\M, 
\end{equation}
which after base change along $\ell_{1}$ is equivariant with respect to
the evident action of $A$ on the right.

By construction, if $s$ is a section of $\M$ then we get a section
$\normq_{i}\M$ of $\normq_{i}\M$ by the formula 
\[
   \quot_{i}^{*}\normq_{i}s = \normff_{\pi} (\mu_{i}^{*}s) \otimes \normff_{\pi} (w_{i}^{*}s)^{-1}.
\]
Under the isomorphism \eqref{eq:86-new}, we have 
\[
     \quot_{1}^{*}\normq_{1}s = \bigotimes_{a\in A} \rT{a}s.
\]

Now suppose that $\L$ is an invertible sheaf of ideals on $G$.

\begin{Lemma}\label{t-le-N-Theta} There is a canonical isomorphism of
rigid line bundles 
\begin{equation} \label{eq:87-new}
  \quot_{i}^{*}\normq_{i}\Theta^{k} \L \cong (\quot^{k})^{*} \Theta^{k}\normq_{\quot}\L.
\end{equation}
\end{Lemma}

\begin{proof}
We prove the statement for $i=1$.  One checks directly that there are
natural isomorphisms of rigid line bundles
\begin{align*}
\rT{a}\Theta^{1}\L & \iso 
                     \Theta^{1}\rT{a}\L \iso
                     \Theta^{1}T_{a}^{*}\L    \\
   \Delta \rT{a}\M & \iso \rT{a}\Delta \M.
\end{align*}
The proof follows by induction, using the
Lemma \ref{t-Delta-facts} and the isomorphism~\eqref{eq:86-new}.
\end{proof}

In view of the Lemma, it is convenient to write
$\tilde{\normq}\Theta^{k}\L$ for the line bundle $\normq_{i}\Theta^{k}\L$,
considered as a line bundle on $(G')^{k}$; of course there is a 
canonical isomorphism of rigid line bundles
\begin{equation} \label{eq:r-norm-theta}
    \tilde{\normq}\Theta^{k}\L\cong \Theta^{k}\normq\L.
\end{equation}

\begin{Proposition}  \label{t-pr-N-i-s-descends-to-theta}
If $s$ is a $\Theta^{k}$-structure on $\L$, then
\[
  \quot_{i}^{*}\normq_{i}s = \quot_{j}^{*}\normq_{j}s, 
\]
and under the isomorphism~\eqref{eq:87-new}, $\normq_{i}s$ descends to a
$\Theta^{k}$-structure on $\normq\L$. 
\end{Proposition}

\begin{proof}
The first point is that $\normq_{i}s$ is independent of $i$.  If
$\M=\Theta^{k}\L$, then there is a canonical isomorphism 
\[
  \frac{\mu_{i}^{*}\M}{w_{i}^{*}\M}\cong 
  \frac{\mu_{j}^{*}\M}{w_{j}^{*}\M}.
\]
The cocycle condition \eqref{was-eq:8} for $s$ implies that, under this
isomorphism,  
\[
  \frac{\mu_{i}^{*}s}{w_{i}^{*}s}\cong 
  \frac{\mu_{j}^{*}s}{w_{j}^{*}s},
\]
so $\quot_{i}^{*}\normq_{i}s = \quot_{j}^{*}\normq_{j}s$ as required.
It follows from Proposition \ref{t-pr-quotient-by-finite-subgroup} that 
$\normq_{i}s$ 
is invariant under the action of the $i$-factor of $K^{k}$ on $G^{k}$.
The question of whether $\normq_{i}s$ descends to a section of
$\Theta^{k}\normq\L$ amounts to the question of whether a ratio of 
sections of $\O_{G^{k}}$ in the equalizer of
\[
\xymatrix{
{\O_{G^{k}}}
 \ar@<1ex>[rr]^-{\pi^{*}}
 \ar@<-1ex>[rr]_-{\mu^{*}}
& & 
{\O_{G^{k}\times K^{k}}},
}
\]
and it suffices to check that it is in the equalizer of 
\[
\xymatrix{
{\O_{G^{k}}}
 \ar@<1ex>[rr]^-{\pi^{*}}
 \ar@<-1ex>[rr]_-{\mu_{i}^{*}}
& & 
{\O_{G^{k}\times K}},
}
\]
for each $i$.
It follows 
that $\normq_{i}s$ descends
to a section of $\Theta^{k}\normq\L$.
It is then straightforward if tedious to check that $\normq_{i}s$ is a
$\Theta^{k}$-structure.
\end{proof}

Lemma \ref{t-le-N-Theta} and Proposition
\ref{t-pr-N-i-s-descends-to-theta} permit us to make the following

\begin{Definition} \label{def-rnorm-section}
If $s$ is a $\Theta^{k}$-structure on $\L$, then let $\tilde{\normq}s$ be the
$\Theta^{k}$-structure on $\normq\L$ such that 
\begin{align}
   (\quot^{k})^{*} \tilde{\normq}s & = \quot_{1}^{*}\normq_{1} s =
   \bigotimes_{a\in A}\rT{a}s \\
\intertext{under the isomorphisms~\eqref{eq:86-new} and~\eqref{eq:87-new}:}
\label{eq-rnorm-idents}
    (\quot^{k})^{*}\Theta^{k} \normq \L        \cong 
    (\quot^{k})^{*}\tilde{\normq}\Theta^{k} \L & \cong 
    \quot_{1}^{*}\normq_{1}\Theta^{k}\L \cong 
    \bigotimes_{a\in A} \rT{a}\Theta^{k}\L.
\end{align}
\end{Definition}

\section{Elliptic curves}

\begin{Definition} \label{def-elliptic-curve}
An \emph{elliptic curve} is a pointed proper smooth curve 
\[
\xymatrix{
{C}
 \ar[r]
&
{S}
 \ar@/_1pc/[l]_{\e}
}      
\]
whose geometric fibers are connected and of genus $1$.
\end{Definition}

Much of the theory of level structures, isogenies,
$\Theta^{k}$-structures, which we have described in detail in this
paper for formal groups, is well-known in the case of elliptic curves. In this
section we briefly recall some results which we need.  Details may be
found in \cite{DeligneRapoport,KaMa:AMEC,Mu:AV,Si:AEC95}.

\subsection{Abel's Theorem} 
\label{sec:abels-theorem}

Note that the discussion of the line bundles $\Theta^{k}\L$ in
\S\ref{sec:thetak-structures} applies to abelian groups in any category where
the notion of line bundle makes sense.  
The first result about elliptic curves is that they are group schemes.

\begin{Theorem}[Abel] \label{t-th-abel}
An elliptic curve $C/S$ has a unique structure of abelian group scheme such
that the  rigid line bundle 
\spdisplay{
\Theta^{3} (\I_{C}(\e))
}
is trivial.  
The (necessarily unique) rigid trivialization $s (C/S)$
of $\Theta^{3} (\I (\e))$ is a cubical structure.   
\end{Theorem}

\begin{proof}
See for example \cite[p. 63]{KaMa:AMEC} or \cite{DeligneRapoport}.
\end{proof}

\begin{Remark}  The \emph{theorem of the cube} says that any line
bundle over an abelian variety has a unique cubical structure.  A
general enough statement of the theorem of the cube, together with the
group structure on elliptic curves,
implies Theorem~\ref{t-th-abel}.  We have stated
Theorem~\ref{t-th-abel} 
to emphasize that the group structure on an elliptic curve is
\emph{constructed} to trivialize $\Theta^{3} (\I (\e))$, so that by
the time you get around to applying the theorem of the cube, you
already know the conclusion for $\I (\e)$.
\end{Remark}

\subsection{Level structures on elliptic curves}
\label{sec:level-struct-ellipt} 

The study of level structures on elliptic curves is due to Katz and
Mazur \cite{KaMa:AMEC}.  
Let $C$ be an elliptic curve over a scheme $S$, and let $A$ be an
abelian group.

\begin{Definition} \label{def-level-structure-ec} 
A homomorphism 
\[
    \ell: A_{S} \xra{} C
\]
is a \emph{level structure} if the Cartier divisor  
\[
\divisorellA =  \sum_{a\in A} \divisor{\ell (a)}
\]
is a
sub-group scheme. 
\end{Definition}

\begin{Lemma}  Let $C$ be an elliptic curve over 
an local formal scheme $S$ with perfect residue field of
characteristic $p>0$.
If 
\[
   \ell: A_{S}\to \fmlgpof{C} 
\] 
is a level structure on the formal group of $C$, then 
\[
   A_{S}\to \fmlgpof{C}  \rightarrow C 
\]
is a level structure on $C$.  
\end{Lemma}

\begin{proof}
This follows from the definition and Proposition \ref{t-pr-level-str-quotient-new}.
\end{proof}

Let 
\[
   A \xra{\ell} C
\]
be a level structure.  The inclusion 
\[
      \divisorellA  \rightarrow C
\]
has a cokernel $C/\divisorellA$ which is an elliptic curve.  If 
\[
\quot:    C\to C' = C/\divisorellA
\]
denotes the projection, then $\quot^{*}$ identifies $\O_{C'}$ with the
equalizer 
\[
\xymatrix{
  {\O_{C'}}
 \ar[rr]
& & 
{\O_{C}}
 \ar@<1ex>[rr]^-{\pi^{*}}
 \ar@<-1ex>[rr]_-{\mu^{*}}
& & 
{\O_{C} \otimes \O_{\divisorellA}.}
}
\]
If $f\in \O_{C},$ then $\normff_{\pi}\mu^{*}f$ is in this equalizer, and we
write 
\[
   \normq = \normq_{\quot}: \O_{C} \rightarrow  \O_{C'}
\]
for the resulting map; explicitly, 
\[
   \quot^{*} \normq f = \prod_{a\in A} T_{a}^{*}f.
\]

It is easy to check that, if $t$ is a coordinate on $C$, then $\normq t$ is
a coordinate on $C'$; in particular, if 
\[
    \ell: A_{S} \rightarrow \fmlgpof{C}
\]
is a level structure, then the natural map of formal groups
\[
     \fmlgpof{C}/\divisorellA \rightarrow \fmlgpof{C/\divisorellA}
\]
is an isomorphism. 

As in the case of a formal group, we have 
\[
    \normq\I_{C} (\e) = \I_{C'} (\e).
\]
After pulling back along the level structure $\ell$, we have 
\[
     \quot^{*}\normq\L \cong \bigotimes_{a\in A} T_{a}^{*}\L, 
\]
and this isomorphism is equivariant with respect to the standard
action of $A$ on the right.  

The discussion of the reduced norm $\Tilde{\normq}$ of
\S\ref{sec:norm-map-thetak-new} applies to elliptic curves as well.  The
main point is that, if $\L$ is a fractional ideal on the elliptic
curve $C$, then the isogeny $\quot$ gives isomorphisms of rigid line bundles
\[
  \Tilde{\normq}\Theta^{k}\L \iso \Theta^{k}\normq\L
\]
over $(C')^{k}$ as in~\eqref{eq:r-norm-theta}, and if $s$ is a 
$\Theta^{k}$ structure on $\L$, then $\Tilde{\normq} s$ is a
$\Theta^{k}$-structure on $\Tilde{\normq}\L$, as in Proposition
\ref{t-pr-N-i-s-descends-to-theta} and Definition \ref{def-rnorm-section}.

\subsection{The Serre-Tate theorem}
\label{sec:serre-tate-theorem}
Let $C_{0}$ be an elliptic curve over a field $k$ of characteristic
$p>0$.

\begin{Definition}  \label{def-deformation-ec}
A \emph{deformation of} $C_{0}$
is a triple $(D/T,f,j)$ consisting of
an elliptic curve $D$ over a local formal scheme $T$ of residue
characteristic $p>0$ and a pull-back diagram 
\[
\begin{CD}
      D_{\specialpoint{T}} @> f    >> C_{0} \\
      @VVV              @VVV \\
      \specialpoint{T} @> j    >> \spec k
\end{CD}
\]
of elliptic curves.  A \emph{map} deformations 
\[
(\alpha,\beta):
(D,f,j)\to (D',f',j')
\]
is 
a pull-back square 
\[
\begin{CD}
      D @> \alpha >> D' \\
   @VVV       @VVV \\
    T  @> \beta >> T'
\end{CD}  
\]
such that the diagram 
\[
\xymatrix{
{D_{\specialpoint{T}}} 
 \ar[r]^{\alpha_{\specialpoint{T}}}
 \ar@/^2pc/[rr]^{f}
 \ar[d]
&
{D_{\specialpoint{T'}}'}
 \ar[r]^{f'}
 \ar[d]
&
{C_{0}}
 \ar[d]
\\
{\specialpoint{T}}
 \ar[r]^{\specialpoint{\beta}}
 \ar@/_2pc/[rr]_{j}
&
{\specialpoint{T'}}
 \ar[r]^{j'}
&
{\spec k}
}
\]
commutes.  
\end{Definition}

Let $C_{0}$ be a supersingular elliptic curve over a 
perfect field $k$ of characteristic $p>0$.

\begin{Theorem}[Serre-Tate] \label{t-th-serre-tate}
The natural transformation 
\[
   \Deformations{C_{0}}^{\noether} \rightarrow 
   \Deformations{\fmlgpof{C_{0}}}^{\noether}
\]
is an isomorphism of functors over $\Fgps^{\noether}$.  
Let $G/S$ be the universal deformation of the formal group
$\fmlgpof{C_{0}}$.  Then there is a deformation $(C/S,\funiv,\juniv)$
of $C_{0}$ to $S$ such that the natural maps 
\[
    \CS^{\noether} \rightarrow  
    \Deformations{C_{0}}^{\noether}  \rightarrow 
    \Deformations{\fmlgpof{C_{0}}}^{\noether} \leftarrow
    \GS^{\noether}
\]
are isomorphism of functors over $\Fgps^{\noether}$.
\end{Theorem}

\begin{proof}
The Serre-Tate Theorem as stated in \cite{Katz:SerreTate} proves that
the forgetful natural transformation induces an isomorphism 
\begin{equation} \label{eq:64}
   \Deformations{C_{0}}^{\artin} \rightarrow 
   \Deformations{\fmlgpof{C_{0}}}^{\artin}
\end{equation}
of functors of Artin local rings.  On the other hand, the functor
$\Deformations{C_{0}}$ is \emph{effectively pro-representable}: there
is deformation $( C'/S',f',j')$  with $S'\iso \spf \W k \psb{u}$, such
that the natural map 
\begin{equation} \label{eq:65}
     (\un{C'/S'})^{\artin} \rightarrow \Deformations{C_{0}}^{\artin}
\end{equation}
is an isomorphism (see for example \cite{DeligneRapoport}).  It follows that 
\[
     (\un{C'/S'})^{\noether} \iso \Deformations{C_{0}}^{\noether}.
\]
Combining the isomorphisms~\eqref{eq:64}
and~\eqref{eq:65} with the isomorphism 
\[
      \Deformations{\fmlgpof{C_{0}}}^{\noether} \iso \GS^{\noether}
\]
gives an isomorphism of formal schemes
\[
     S \xra{\iso} S'
\]
and, if $C$ is the elliptic curve over $S$ obtained from $C'/S'$ by
pull-back, an isomorphism
\[
    \CS^{\noether} \rightarrow \GS^{\noether}.
\]
\end{proof}

\begin{Example} \label{ex-supsing-char-2}
In characteristic $2$ the elliptic curve $C_{0}$ given by the 
Weierstrass equation 
\[
    y^{2} + y = x^{3}
\]
is supersingular (e.g. \cite{Si:AEC95}).  The universal deformation of
its formal group is a formal group $G$ over $S\iso \spf
\Z_{2}\psb{u_{1}}$.  It is well-known (e.g. by the Exact Functor
Theorem \cite{Landweber:Exa}) that there is a spectrum $\E$ with 
\[
   \GpOf{\E}/S_{\E} = G/S:
\]
it is a form of $E_{2}$.
The Serre-Tate Theorem endows $\E$ with the structure of an elliptic
spectrum: if $C/S$ is the universal deformation of $C_{0}$ to $S$, 
then there is a canonical isomorphism 
\[
 \GpOf{\E}= G\iso \fmlgpof{C}
\]
of formal groups over $S_{\E}$.
\end{Example}

\subsection{Descent for level structures on a  Serre-Tate curve} 
\label{sec:desc-level-struct-e-c}

We suppose that $C_{0}$ is a supersingular elliptic curve over a
perfect field $k$ of characteristic $p>0$, and that $C/S$ is the
universal deformation of $C_{0}$ provided by Theorem
\ref{t-th-serre-tate}.  Let $G/S$ be the universal deformation of the
formal group $\fmlgpof{C_{0}}$.  Since $\CS$ is a functor over $\Fgps$,
Definition 
\ref{def-descent-functors} provides a notion of descent data for level
structures on $\CS$.  

Theorem \ref{t-th-serre-tate} 
gives an isomorphism of formal groups $\fmlgpof{C}\iso G$, and so 
the descent
data \eqref{eq:50} for the 
Lubin-Tate formal groups give descent data  
\[
   d_{1}: \ulvl_{1} (\CS)
          \rightarrow
          \CS.
\]
Explicitly, suppose that 
\[
    A_{T} \xra{\ell} i^{*}\fmlgpof{C}
\]
is a level structure over a Noetherian local formal scheme $T$.  The
descent data provide an isogeny of formal groups 
\begin{equation}\label{eq:69}
    i^{*}\fmlgpof{C} \xra{f_{\ell}} \psi_{\ell}^{*}\fmlgpof{C}
\end{equation}
over $T$ with kernel $\divisorellA.$   

It is natural to ask for an isogeny of elliptic curves 
\[
        i^{*} C \xra{g_{\ell}} \psi_{\ell}^{*}C
\]
extending $f_{\ell}$.  This corresponds, in the language of
\S\ref{sec:desc-level-struct}, to replacing the functor  $\Fgps$ with
the functor $\EllCurves$ whose value on a ring $R$
is the set of elliptic curves $C/\spec R$.  
Thus we shall refer to 
descent data for level structures on $\CS$ together with isogenies
$g_{\ell}$ extending $f_{\ell}$ as \emph{descent data for level
structures on $\CS$ over $\EllCurves$}.

Over $T$ we have an isogeny of elliptic curves
\[
       i^{*} C \xra{\quot} C'
\]
with kernel $\divisorellA$, and 
a canonical isomorphism 
\begin{equation} \label{eq:68}
    \fmlgpof{C'}\iso\psi_{\ell}^{*} \fmlgpof{C}
\end{equation}
of formal groups over $T$, as in \S\ref{sec:level-struct-ellipt}.
Theorem \ref{t-th-serre-tate} implies that there is a unique
isomorphism of elliptic curves  
\[
     C'\iso \psi_{\ell}^{*}C
\]
extending \eqref{eq:68}; put another way, we have the following.

\begin{Corollary}\label{t-co-hinfty-ec}
The functor $\CS$
has descent data for level structures over $Ell$, whose
restriction to $\underline{\fmlgpof{C}/S} = \GS$ are the descent
data given by  Proposition \ref{t-pr-lubin-tate-hinfty}. 
In particular, for each level structure
\[
     A_{T} \xra{\ell} i^{*} \fmlgpof{C},
\]
there is a canonical isogeny of elliptic curves $g_{\ell}$ making the
diagram 
\[
\begin{CD}
i^{*}\fmlgpof{C} @>>> i^{*} C \\
@V f_{\ell} VV @VV g_{\ell} V \\
\psi_{\ell}^{*}\fmlgpof{C} @>>> \psi_{\ell}^{*} C
\end{CD}    
\]
commute. \qed
\end{Corollary}

\subsection{The cubical structure of an elliptic curve is compatible with
descent}
\label{sec:cubical-structure-descent}

The uniqueness of the cubical structure in Theorem \ref{t-th-abel}
implies the following.

\begin{Proposition} \label{t-prop-ell-curve-hinfty}
Let $C$ be an elliptic curve, and let $s (C/S)$ be the cubical
structure of Theorem \ref{t-th-abel}.  If $i^{*}C\to \psi^{*}C$ is an
isogeny, then  
\[
        \psi^{*}s (C/S) = \Tilde{\normq}i^{*}s (C/S).
\] \qed
\end{Proposition}


\section{The sigma orientation}
\label{sec:sigma-orientation}

Suppose that $\E$ is a homogeneous ring spectrum and let $G=G_{\E}$.
Let $V$ be the line bundle 
\[
   V = \prod_{j=1}^{k} (1-L_{i})
\]
over $(\cp)^{k}$.  In Lemma \ref{t-le-Theta-and-BUn} we observed that 
Proposition \ref{t-pr-line-v} specializes to give an isomorphism
\[
    t_{V}: \Line (V)\iso \Theta^{k} (\I_{G_{\E}} (\e)),
\]
and if 
\[
    g: \MU{2k} \to \E
\]
is an orientation, then the composition 
\[
(\cp^k)^V\to\MU{2k}\xra{g}\E
\]
represents a rigid section $s$ of $\Theta^{k} (\I_{G} (\e))$.  In fact
it is easily seen to be a $\Theta^{k}$-structure, that is a
$\pi_{0}\E$-valued point of $\uC^{k} (G;\I_{G} (\e))$.  
Similarly, if $g: \BU{2k}_{\plus} \to \E$ is a homotopy multiplicative
map, then the composite 
\[
  \cp^{k} \to \BU{2k} \to \E
\]
represents a $\Theta^{k}$-structure on the trivial line bundle
$\O_{G}$, and so a point of $\uC^{k} (G;\O_{G})$.  In 
\cite{AHS:ESWGTC} we proved 

\begin{Theorem} \label{t-th-ring-spectra-mu-to-e}
If $\E$ is a homogeneous spectrum and $k\leq 3$, then these
correspondences induce isomorphisms 
\begin{equation} \label{eq:46}
 \RingSpectra (\MU{2k},\E) \rightarrow \uC^{k} (G;\I_{G} (\e)) (\pi_{0}\E)
\end{equation}
and 
\[
\RingSpectra (\BU{2k}_{\plus},\E) \rightarrow \uC^{k} (G;\O_{G}) (\pi_{0}\E).
\] \qed
\end{Theorem}

Now suppose that $(\E,C,t)$ is an \emph{elliptic spectrum}: that is,
$\E$ is a homogeneous ring spectrum, $C$ is an elliptic curve over
$S_{\E}$, and $t$ is an isomorphism
\[
     t: \GpOf{\E}\iso \fmlgpof{C}
\]
of formal groups over $S_{\E}$.  Abel's Theorem \ref{t-th-abel} gives
a cubical structure $s (C/S)$ on $C$, which gives a  cubical
structure $t^{*}\fmlgpof{s} (C/S)$ on $\GpOf{\E}$.

\begin{Definition}\cite{AHS:ESWGTC} \label{def-2}
The \emph{sigma orientation} for $(\E,C,t)$ is the map of ring spectra 
\[
\sigma (\E,C,t): \MU{6} \rightarrow \E
\]
which corresponds to $t^{*}\fmlgpof{s} (C/S)$ under the
isomorphism~\eqref{eq:46}.
\end{Definition}

Now suppose that $\E$ is a homogeneous $\hinfty$ spectrum, with the
property that $\pi_{0}\E$ is an
admissible 
local ring
with perfect residue field of characteristic $p>0$.  Let $S=S_{\E}$.
Suppose that $(\E,C,t)$ is an elliptic spectrum.  In particular, the 
$G=\GpOf{\E}$ is 
of finite height.  By Theorem
\ref{t-th-hinfty-adds}, the $\hinfty$ structure on $\E$ gives descent
data for level structures on $G$.

\begin{Definition}\label{def-3}
An \emph{$\hinfty$ elliptic spectrum} is an elliptic spectrum
$(\E,C,t)$ whose underlying spectrum $\E$ is a homogeneous $\hinfty$
spectrum $\E$ as above, together with descent data for level
structures on $\CS$, considered as a functor over $\EllCurves$ as in
\S\ref{sec:desc-level-struct-e-c}, such that the diagram  of functors
over $\Fgps$  
\[
\begin{CD}
\ulvl_{1} (\CS) @> t >> \ulvl_{1} (\GS) \\
@V d_{1} VV @VV d_{1} V \\
\CS @> t >> \GS
\end{CD}
\]
commutes.
\end{Definition}

\begin{Proposition} \label{t-pr-sigma-hinfty-hinfty}
Let $(\E,C,t)$ be an $\hinfty$ elliptic spectrum, and suppose in
addition that $p$ is regular in $\pi_{0}\E$.  
Then the sigma orientation 
\[
    \musix \xra{\sigma_{(\E,C,t)}} \E
\]
is an $\hinfty$ map.
\end{Proposition}

\begin{proof}
By Proposition \ref{t-pr-norm-condition-suffices}, it suffices to show
that, for each level structure 
\[
    A_{T} \xra{\ell} i^{*} \fmlgpof{C},
\]
we have
\[
 \tilde \normq_{g_{\ell}} s (C/S_{\E}) = \psilestar s (C/S_{\E}),
\]
where $g_{\ell}$ is the isogeny of elliptic curves making the diagram 
\[
\begin{CD}
i^{*}\GpOf{\E}  @>>> i^{*} C \\
@V \psilge{\E} VV @VV g_{\ell} V \\
\psilestar \GpOf{\E}  @>>> \psilestar C
\end{CD}    
\]
commute.
Proposition \ref{t-prop-ell-curve-hinfty} gives the result.
\end{proof}

Now let $(\E,C,t)$ be the elliptic spectrum associated to the
universal deformation of a supersingular elliptic curve $C_{0}$ over a perfect
field $k$ of characteristic $p>0$.  For example, we may take 
$C_{0}$ to be the Weierstrass curve
\[
 y^{2} + y = x^{3}
\]
over $\GF{2}$ (Example~\ref{ex-supsing-char-2}).  Applying the
Proposition, Corollary \ref{t-co-hinfty-ec}, and
Proposition \ref{t-pr-d-1-top-is-d-1-frob} gives the 

\begin{Corollary} \label{t-co-sigma-to-E-2-hinfty}
The orientation 
\[
    \musix\xra{\sigma (\E,C,t)} \E
\]
is an $\hinfty$ map. \qed
\end{Corollary}

\appendix

\section{$\hinfty$-ring spectra}
\label{sec:hinfty-ring-spectra}

Given an integer $n\ge 0$, let $D_n:\Spectra U\to\Spectra U$ be the
functor
\begin{equation} \label{eq:28}
\E\mapsto \lien n\underset{\Sigma_n}{\wedge}\E^{(n)},
\end{equation}
where $\lien n=\lie(U^n,U)$ is the space of linear isometric embeddings
from $U^n$ to $U$.

An $E_\infty$ ring spectrum is a spectrum with maps 
\[
D_n(\E)\to\E,\quad n\ge 0,
\]
making the following diagrams commute:
\begin{equation} \label{eq:42}
\begin{CD}
\{1_U\}\rtimes \E@>>> D_1\E \\
@VVV @VVV \\
\E @= \E
\end{CD}
\qquad
\begin{CD}
D_n D_m\E @>>> D_{n+m}\E \\
@VVV @VVV \\
D_n \E @>>> \E.
\end{CD}
\end{equation}
An $H_\infty$ ring spectrum is a spectrum $\E$ together with maps
$D_n\E\to\E$ such that the diagrams~\eqref{eq:42} commute
up to homotopy.

The category of $\E_\infty$-ring spectra is naturally enriched over
topological spaces.  The \emph{space} of $E_\infty$-maps from $\E$ to
$\F$ is the
subspace of all maps consisting of those which make the diagrams
\[
\begin{CD}
D_n\E @>>> D_n\F \\
@VVV @VVV \\
\E @>>> \F
\end{CD}
\]
commute.  For a topological space $X$, the spectrum which
underlies the ``function object'' is simply the spectrum
$\E^{X_{\plus}}$.  
The spectrum which underlies $\E\otimes X$ is more difficult to
describe.  If $\E$ is only an $H_\infty$-ring spectrum, the spectrum
$\E^{X_{\plus}}$ is still $H_\infty$.

These remarks actually depend very little on the construction of the
functor $D_n$ and are mostly matters of pure category theory.  Indeed,
the map $D_n\E\to\E$ can be regarded as a natural transformation of
functors
\[
\Spectra U(\F,\E)\to\Spectra U(D_n\F,\E).
\]
Given a topological space $X$, we can
use~\eqref{eq:44} to define a transformation
\begin{equation} \label{eq:43}
\Spectra U(\F,\E^{X_{\plus}})\to\Spectra U(D_n\F,\E^{X_{\plus}}).
\end{equation}
as the composite
\begin{align*}
\Spectra U(\F,\E^{X_{\plus}}) &\iso\Spaces(X,\Spectra U(\F,\E)) \\
&\to\Spaces(X,\Spectra U(D_n\F,\E)) \\
&\iso\Spectra U(D_n\F,\E^{X_{\plus}}).
\end{align*}
A more subtle property is that the
transformation~\eqref{eq:43} is also given by
\begin{equation}\label{eq:45}
\begin{aligned}
\Spectra U(\F,\E^{X_{\plus}}) \iso\Spectra U(\F\wedge {X_{\plus}},\E))
&\xra{\phantom{\text{diag}}} \Spectra U(D_n(\F\wedge X_{\plus}),\E)) \\
& \xra{\text{diag}}\Spectra U(D_n(\F)\wedge X_{\plus}),\E) \\
&\iso\Spectra U(D_n\F,\E^{X_{\plus}}).
\end{aligned}
\end{equation}

An important property of the functors is summarized in the following
result of \cite{BMMS:Hin}.

\begin{Proposition}\label{t-pr-may-et-al}
There is a natural
weak equivalence
\[
\bigvee_{i+j=n} \lien2\wedge D_i(\E)\wedge D_j(\F)\to D_n(\E\vee\F).
\]
Furthermore, the $ij$-component of 
\splitpageyesno{
\begin{multline*}
D_n(\E) \xra{D_n(\nabla)} D_n(\E\vee\E) \\
\to \bigvee_{i+j=n} \lien2\wedge D_i(\E)\wedge D_j(\E) \\
\to
\prod_{i+j=n} \lien2\wedge D_i(\E)\wedge D_j(\E) 
\end{multline*}
}
{
\[
D_n(\E) \xra{D_n(\nabla)} D_n(\E\vee\E) 
\to \bigvee_{i+j=n} \lien2\wedge D_i(\E)\wedge D_j(\E)
\to
\prod_{i+j=n} \lien2\wedge D_i(\E)\wedge D_j(\E) 
\]
}
is the transfer map $\tr_{ij}$ with respect to the inclusion
$\Sigma_i\times\Sigma_j\subset\Sigma_n$.
\end{Proposition}

Note also that if $W$ is a virtual bundle of dimension $0$ over a
space $X$, then $D_A(X^W)$ is the Thom spectrum the virtual bundle 
$V_{\text{reg}}\otimes W$ over $D_A(X)$, where $V_{\text{reg}}$ is the
regular representation of $A^\ast$.

\section{Composition of operations}

\label{sec:comp-oper-deriv}

Let $\E$ be a homogeneous $\hinfty$ ring spectrum, such that 
$\pi_{0}\E$ is 
a 
local ring
with perfect residue field of characteristic $p>0$, and the formal group
$G=\GpOf{\E}$ is 
of finite height. 
In \S\ref{sec:algebr-geom-hinfty} we associate to
each map of formal schemes 
\[
i: T = \spf R \to S_{\E},
\]
finite abelian group $A$, and level structure
\[
\ell:A_{T}\to i^\ast G,
\]
a map of formal schemes $\psile: T\to S_{\E}$ and an isogeny
$\psilge{\E}:i^\ast G\to\psi_\ell^\ast G$ with kernel $\divisorellA$.  Theorem
\ref{t-th-hinfty-adds} asserts that these constitute descent data for
level structures on $G$.  Properties (1) and (3) of Definition
\ref{def-1} follow immediately from the construction.  In this section
we give a proof of (2): if 
\begin{equation} \label{eq:75}
 B \rightarrow A \xra{\pi}  C
\end{equation}
is a short exact sequence of finite abelian groups, then with
the notation
\begin{equation} \label{eq:79}
\begin{CD}
B@>>> A @>>> C \\
@V \ell' VV @V \ell VV @VV \ell'' V \\
i^\ast G @= i^\ast G @>> \psilpge{\E}> \psi_{\ell'}^{*} G,
\end{CD}
\end{equation}
we have 
\begin{align}
\psi_{\ell''} &  = \psi_{\ell} : T \rightarrow S  \label{eq:80}\\
\psilge{\E} &  = \psilpge{\E} \circ \psilppge{\E}:  i^{*}G \rightarrow
\psi_{\ell}^{*}G = \psi_{\ell''}^{*}G.\label{eq:81}
\end{align}
In fact we shall give a proof of~\eqref{eq:80}; the proof
of~\eqref{eq:81} is similar.

Elsewhere in this paper we use the notation $D_{n}$ for either
the functor~\eqref{eq:28} or for
the natural transformation 
\[
\Spectra U(\F,\E)\to\Spectra U(D_n\F,\E)
\]
associated to an $\hinfty$ ring spectrum $\E$ (and similarly for $D_{A}$).
In this section only, it is convenient to write $P_{n}$ or $P_{A}$ for
the natural transformation, reserving $D_{n}$ and $D_{A}$ to refer to
the functors~\eqref{eq:28} and~\eqref{eq-a-extended-power}.  Thus for 
\[
     f: \F\to \E,
\]
the diagram 
\[
\xymatrix{
 {D_{n}\F}
   \ar[r]^{P_{n}f}
   \ar[d]_{D_{n}f}
 &
{\E}
 \\
{D_{n}\E}
 \ar[ur]
}
\]
commutes; and if $R$ is a complete local ring and 
\[
    A_{\spf R} \xra{\ell} i^{*}G
\]
is a level structure, then by Definition~\ref{def-psile},
$\psile:\pi_{0}\E \rightarrow R$ is the composition 
\spmline{
\pi_0\E
 \xra{P_{A}}
\pi_0\Spectra U(D_A S^0,\E) =
\pi_0\E^{BA^\ast_\plus} \rightarrow
}
{
 \O ((BA^{\ast})_{\E}) \iso 
 \O (\uhom (A,G)) 
\xra{\chi_{\ell}}
R.
}

If $Z$ is a set, let $\Sigma_{Z}$ be the group of automorphisms of the
underlying set.  We let $Z$ act on itself by left multiplication, and
so we consider $Z$ to be a subgroup of $\Sigma_{Z}$.  If $T\subseteq
\Sigma_{Z}$ is a subgroup, and if $X$ is a group, then we write $T\wreath
X$ for the wreath product
\[
   T\wreath X = T \ltimes X^{Z}.
\]

Suppose that 
\begin{equation} \label{eq:78}
     X \rightarrow Y \xrightarrow{\pi} Z
\end{equation}
is a short exact sequence of abelian groups, and $Z$ is finite.  A
splitting 
\[
   s: Z \to Y
\]
of~\eqref{eq:78} as a sequence of sets determines a homomorphism 
\[
    g: Y \to Z\wreath X
\]
by the formula 
\[
    g (y) = (\pi (y), f_{y}),
\]
where 
\[
   f_{y}: Z \to X
\]
is the map of sets given by the formula 
\[
    f_{y} (z) = y + s (z) - s (\pi (y) + z).
\]
If $g'$ is another such homomorphism defined using a section $s'$,
then $g$ and $g'$ differ by conjugation by $(0,s-s')\in Z\wreath X$,
and so just the extension \eqref{eq:78} determines a homotopy class of maps 
\[
    BY \xrightarrow{Bg} B (Z\wreath X).
\]
If moreover $X$ is finite, then $s$ determines a homomorphism 
\[
    h: \Sigma_{Z}\wreath \Sigma_{X} \rightarrow \Sigma_{Y}
\]
by the formula 
\[
    h(\sigma,\tau) (s (z) + x) = s (\sigma (z)) + \tau (x)
\]
for $\sigma\in \Sigma_{Z}, \tau\in \Sigma_{X}, x\in X,$ and $z\in Z$.
Once again the resulting map 
\[
B (\Sigma_{Z}\wreath \Sigma_{X}) \rightarrow B\Sigma_{Y}
\]
is independent of the choice of section.  The maps $g$ and $h$ have been
defined so that the diagram 
\[
\begin{CD}
Z\wreath X @>>> \Sigma_{Z} \wreath \Sigma_{X} \\
@A g AA @VV h V  \\
Y @>>> \Sigma_{Y}
\end{CD}
\]
commutes.

Applying these observations to the dual of a short exact sequence
\[
   B \rightarrow  A \rightarrow C
\]
gives maps 
\[
    Bg_{\plus}: 
   D_{A}S^{0}  = BA^{*}_{\plus} \rightarrow  
                 B\left(  B^{*} \wreath BC^{*}\right)_{\plus}  
               =  D_{B}BC^{*}_{\plus} = D_{B}D_{C}S^{0} 
\]
 and 
\[
    Bh_{\plus}: D_{|B|}D_{|C|}S^{0} \rightarrow  D_{|A|}S^{0}.
\]

\begin{Lemma} \label{t-le-comp-ops}
In this situation, the diagram 
\[
\begin{CD}
\pi_{0}\E @> P_{A} >> \pi_{0}\E^{BA^{*}_{\plus}} \\
@V P_{C} VV @AA \pi_{0}\E^{Bg_{\plus}}A \\
\pi_{0}\E^{BC^{*}_{\plus}} @>> P_{B} > \pi_{0}\E^{D_{B^{*}}BC^{*}_{\plus}}
\end{CD}
\]
commutes.
\end{Lemma}

\begin{proof}
For $f: T\to E$, $P_{A}f$ is the composition 
\[
    D_{A}T \rightarrow D_{|A|} T \xra{D_{|A|}f} D_{|A|} \E \xra{} \E.
\]
Now let $f: S^{0} \to \E$, and consider the diagram 
\[
\xymatrix{
{D_{B}D_{C}S^{0}}
 \ar   `r[ddrrrr] [ddrrrr]^-{P_{B}P_{C}f}
 \ar[dr] \\
&
{D_{|B|}D_{C}S^{0}} 
 \ar[rr]_-{D_{|B|}P_{C}f}
 \ar[d]
& & 
{D_{|B|}\E}
 \ar[dr]
\\
&
{D_{|B|}D_{|C|}S^{0}} 
 \ar[rr]_-{D_{|B|} D_{|C|}f}  
 \ar[d]_{Bh_{\plus}}
& & 
{D_{|B|} D_{|C|} \E}
 \ar[d]
 \ar[u]
&
{\E} \\
&
{D_{|A|} S^{0}}
 \ar[rr]^-{D_{|A|} f}
& & 
{D_{|A|}\E}
 \ar[ur]
\\
D_{A}S^{0}
 \ar[uuuu]^{Bg_{\plus}}
 \ar[ur]
 \ar `r[uurrrr] [uurrrr]_{P_{A}f}
}
\]
The trapezoid on the left commutes by Lemma \ref{t-le-comp-ops}.  The
top inner rectangle is obtained by applying $D_{|B|}$ to the
definition of $P_{C}f$, and so the top outer composition is the
definition of $P_{B}P_{C}f$.  The lower inner rectangle commutes by
the naturality of $D$, and the lower outer composition is the definition
of $P_{A}$.  The right inner triangle is a case of the left diagram
of~\eqref{eq:42}, and so commutes because $\E$ is an $\hinfty$
spectrum.
\end{proof}

Now we turn to the situation of the diagram~\eqref{eq:79}.
Lemma \ref{t-le-comp-ops} implies that the top square in the diagram 
\begin{equation} \label{eq:77}
\xymatrix{
{\pi_{0}\E}
  \ar[d]_{P_{C}}
  \ar[r]^{P_{A}}
 \ar `u[r]
     `[drr]
      [drr]^{\psile}
&
{\pi_{0}\E^{BA^{*}_{+}}}
 \ar[dr]^{\chi_{\ell}}
\\
{\pi_{0}\E^{BC^{*}_{+}}}
 \ar[r]_{P_{B}}
&
{\pi_{0}\E^{D_{B^{*}}BC^{*}}_{+}}
 \ar[u]^{\pi_{0}\E^{Bg_{\plus}}} 
&
{R}
\\
{\pi_{0}\E}
 \ar[u]
 \ar[r]^{P_{B}}
 \ar `d[r]
     `[urr]
      [urr]_{\psilpe}
&
{\pi_{0}\E^{BB^{*}_{+}}}
 \ar[u]
 \ar[ur]_{\chi_{\ell'}}
}
\end{equation}
commutes.  The commutativity of the bottom left square is obvious; the 
two right-hand corners
commute by the definition of $\psile$~\eqref{def-psile}.

\begin{Lemma} \label{t-le-chiellpp}
If 
\[
    \chi_{\ell''} : \pi_{0}\E^{BC^{*}_{\plus}} \rightarrow R
\]
is the homomorphism classifying the homomorphism 
\[
    C \xra{\ell''} (\psi_{\ell'}^{\E})^{*} G,
\]
then 
\begin{equation} \label{eq:76}
   \chi_{\ell''}  = \chi_{\ell} \circ \pi_{0}\E^{Bg_{\plus}} \circ P_{B}.
\end{equation}
\end{Lemma}

\begin{proof}
By definition, $\chi_{\ell''}$ classifies the homomorphism 
\[
   C \xra{\ell} \psi_{\ell'}^{*} G.
\]
Let us temporarily write 
\[
    u = \chi_{\ell} \circ \pi_{0}\E^{Bg_{\plus}} \circ P_{B}.
\]
As in the proof of Lemma \ref{lem-psi-additive},  the axioms of an
$H_{\infty}$ structure together with 
Proposition \ref{t-pr-ra-kills-transfers} imply that $u$ is a
continuous ring  homomorphism.  
The commutativity of the diagram~\eqref{eq:77} shows that the
diagram 
\[
\xymatrix{
{\pi_{0}\E}
 \ar[d]
 \ar[dr]^{\psilpe}
\\
{\pi_{0}\E^{BC^{*}_{+}}}
 \ar[r]_-{u}
&
{R}
}
\]
commutes.  In view of the isomorphism 
\[
     (BC^{*})_{\E} \iso \uhom (C,G)
\]
of Proposition \ref{t-E-BA-ast}, $u$ classifies \emph{some} homomorphism 
\[
    C \xra{w} \psilpestar G,
\]
and it remains to show that $w=\ell''$.
To show that, it suffices to 
show that the diagram 
\[
\begin{CD}
A @> \ell >> i^{*} G \\
@VVV         @VV \psilge{\E} V \\
C @> w >> \psilpestar G
\end{CD}   
\]
commutes.  

Suppose that $a\in A$ with image $c\in C$, and consider the diagram 
\[
\xymatrix{
{\pi_{0}\E^{\cpplus}}
 \ar[r]^{P_{B}}
 \ar[d]_{\pi_{0}\E^{Bc_{\plus}}}
&
{\pi_{0}\E^{D_{B}\cp_{\plus}}}
 \ar[d]^{\pi_{0}\E^{D_{B}c_{\plus}}}
 \ar[r]^{\pi_{0}\E^{B\Delta_{\plus}}}
&
{\pi_{0}\E^{(BB^{*}\times \cp)_{\plus}}}
 \ar[d]^{\pi_{0}\E^{B (\pi\times a)_{\plus}}}
\\
{\pi_{0}\E^{BC^{*}_{\plus}}}
 \ar[r]^{P_{B}}
&
{\pi_{0}\E^{D_{B}BC^{*}_{\plus}}}
 \ar[r]^{\pi_{0}\E^{Bg_{\plus}}}
&
{\pi_{0}\E^{BA^{*}_{\plus}}}
 \ar[r]^{\chi_{\ell}}
&
{R.}
}
\]
The outer clockwise composition is the map of rings corresponding to
the point 
\[
     \spf R \xra{ \psilge{\E} (\ell (a)) } G,
\]
while the outer counterclockwise composition corresponds to the point  
\[
     \spf R \xra{ w (c)} G.
\]
It is clear that the left square commutes.  To see that the right
square commutes, observe that we have a commutative diagram 
\[
\begin{CD}
C^{*} @>>> A^{*} @>>> B^{*} \\
@V c VV       @VV \pi \times a V        @VV = V \\
\C^{\times} @>>> B^{*}\times \C^{\times} @>>> B^{*}. 
\end{CD}
\]
A choice of section 
\[
   s: B^{*} \to A^{*}
\]
gives the map $Bg : BA^{*} \rightarrow  D_{B}BC^{*}$.  It also gives 
a section 
\[
  s: B^{*} \to B^{*}\times \C^{\times},
\]
and so a map 
\[
   Bg : BB^{*} \times \cp \rightarrow D_{B}\cp
\]
such that the diagram 
\[
\begin{CD}
BA^{*} @> Bg >> D_{B}BC^{*} \\
@V B (\pi \times a) VV          @VV D_{B} Bc V \\
BB^{*}\times \cp @> Bg >> D_{B} \cp
\end{CD}
\]
commutes.  But the homotopy class of the map $Bg$ 
is independent of the choice of section, so
the diagram 
\[
\begin{CD}
BA^{*} @> Bg >> D_{B}BC^{*} \\
@V B (\pi \times a) VV          @VV D_{B} Bc V \\
BB^{*}\times \cp @> B\Delta >> D_{B} \cp
\end{CD}
\]
commutes.
\end{proof}

The commutativity of the diagram \eqref{eq:77} and Lemma
\ref{t-le-chiellpp} together imply that 
\begin{align*}
   \psile & = \chi_{\ell} \circ P_{A} \\
          & = \chi_{\ell} \circ \pi_{0}\E^{Bg_{\plus}} \circ P_{B}\circ P_{C} \\
          & = \chi_{\ell''} \circ P_{C}  \\
          & = \psilppe
\end{align*}
as required.


\end{document}